\input psfigb.sty
%
%
\catcode`@=11 
%
%

\def\b@lank{ }

\newif\if@simboli
\newif\if@riferimenti
\newif\if@bozze

\def\bozze{\@bozzetrue\font\tt@bozze=cmtt8}
\def\og@gi{\number\day\space\ifcase\month\or 
   gennaio\or febbraio\or marzo\or aprile\or maggio\or giugno\or 
   luglio\or agosto\or settembre\or ottobre\or novembre\or dicembre\fi
   \space\number\year}
\newcount\min@uti
\newcount\or@a
\newcount\ausil@iario
\min@uti=\number\time
\or@a=\number\time
\divide\or@a by 60
\ausil@iario=-\number\or@a
\multiply\ausil@iario by 60
\advance\min@uti by \number\ausil@iario
\def\ora@esecuzione{\the\or@a:\the\min@uti}  
\def\makefootline{\baselineskip=24pt\line{\the\footline}
    \if@bozze\vskip-10pt\tt@bozze
             \noindent \jobname\hfill\og@gi, ore \ora@esecuzione\fi}

\newwrite\file@simboli
\def\simboli{
    \immediate\write16{ !!! Genera il file \jobname.SMB }
    \@simbolitrue\immediate\openout\file@simboli=\jobname.smb}

\newwrite\file@ausiliario
\def\riferimentifuturi{
    \immediate\write16{ !!! Genera il file \jobname.AUX }
    \@riferimentitrue\openin1 \jobname.aux
    \ifeof1\relax\else\closein1\relax\input\jobname.aux\fi
    \immediate\openout\file@ausiliario=\jobname.aux}

\newcount\eq@num\global\eq@num=0
\newcount\sect@num\global\sect@num=0

\newif\if@ndoppia
\def\numerazionedoppia{\@ndoppiatrue\gdef\la@sezionecorrente{\the\sect@num}}

\def\se@indefinito#1{\expandafter\ifx\csname#1\endcsname\relax}
\def\spo@glia#1>{} 

\newif\if@primasezione
\@primasezionetrue

\def\s@ection#1\par{\immediate
    \write16{#1}\if@primasezione\global\@primasezionefalse\else\goodbreak
    \vskip\spaziosoprasez\fi\noindent
    {\sezfont #1}\nobreak\vskip\spaziosottosez\nobreak\noindent}
%

\font\sezfont=cmbx10

\def\sezpreset#1{\global\sect@num=#1
    \immediate\write16{ !!! sez-preset = #1 }   }

\def\spaziosoprasez{50pt plus 60pt}
\def\spaziosottosez{15pt}
\def\spaziotitsez{5truemm}

\def\sref#1{\se@indefinito{@s@#1}\immediate\write16{ ??? \string\sref{#1}
    non definita !!!}
    \expandafter\xdef\csname @s@#1\endcsname{??}\fi\csname @s@#1\endcsname}

\def\autosez#1#2\par{
    \global\advance\sect@num by 1\if@ndoppia\global\eq@num=0\fi
    \xdef\la@sezionecorrente{\the\sect@num}
    \def\usa@getta{1}\se@indefinito{@s@#1}\def\usa@getta{2}\fi
    \expandafter\ifx\csname @s@#1\endcsname\la@sezionecorrente\def
    \usa@getta{2}\fi
    \ifodd\usa@getta\immediate\write16
      { ??? possibili riferimenti errati a \string\sref{#1} !!!}\fi
    \expandafter\xdef\csname @s@#1\endcsname{\la@sezionecorrente}
    \immediate\write16{\la@sezionecorrente. #2}
    \if@simboli
      \immediate\write\file@simboli{ }\immediate\write\file@simboli{ }
      \immediate\write\file@simboli{  Sezione 
                                  \la@sezionecorrente :   sref.   #1}
      \immediate\write\file@simboli{ } \fi
    \if@riferimenti
      \immediate\write\file@ausiliario{\string\expandafter\string\edef
      \string\csname\b@lank @s@#1\string\endcsname{\la@sezionecorrente}}\fi
    \goodbreak\vskip\spaziosoprasez
    \noindent\if@bozze\llap{\tt@bozze#1\ }\fi
      {\sezfont\the\sect@num.\hskip\spaziotitsez #2}\par\nobreak
    \vskip\spaziosottosez\nobreak\noindent}

\def\semiautosez#1#2\par{
    \gdef\la@sezionecorrente{#1}\if@ndoppia\global\eq@num=0\fi
    \if@simboli
      \immediate\write\file@simboli{ }\immediate\write\file@simboli{ }
      \immediate\write\file@simboli{  Sezione ** : sref.
          \expandafter\spo@glia\meaning\la@sezionecorrente}
      \immediate\write\file@simboli{ }\fi
    \s@ection#2\par}


\def\eqpreset#1{\global\eq@num=#1
     \immediate\write16{ !!! eq-preset = #1 }     }

\def\eqref#1{\se@indefinito{@eq@#1}
    \immediate\write16{ ??? \string\eqref{#1} non definita !!!}
    \expandafter\xdef\csname @eq@#1\endcsname{??}
    \fi\csname @eq@#1\endcsname}

\def\eqlabel#1{\global\advance\eq@num by 1
    \if@ndoppia\xdef\il@numero{\la@sezionecorrente.\the\eq@num}
       \else\xdef\il@numero{\the\eq@num}\fi
    \def\usa@getta{1}\se@indefinito{@eq@#1}\def\usa@getta{2}\fi
    \expandafter\ifx\csname @eq@#1\endcsname\il@numero\def\usa@getta{2}\fi
    \ifodd\usa@getta\immediate\write16
       { ??? possibili riferimenti errati a \string\eqref{#1} !!!}\fi
    \expandafter\xdef\csname @eq@#1\endcsname{\il@numero}
    \if@ndoppia
       \def\usa@getta{\expandafter\spo@glia\meaning
       \la@sezionecorrente.\the\eq@num}
       \else\def\usa@getta{\the\eq@num}\fi
    \if@simboli
       \immediate\write\file@simboli{  Equazione 
            \usa@getta :  eqref.   #1}\fi
    \if@riferimenti
       \immediate\write\file@ausiliario{\string\expandafter\string\edef
       \string\csname\b@lank @eq@#1\string\endcsname{\usa@getta}}\fi}

\def\autoreqno#1{\eqlabel{#1}\eqno(\csname @eq@#1\endcsname)
       \if@bozze\rlap{\tt@bozze\ #1}\fi}
\def\autoleqno#1{\eqlabel{#1}\leqno\if@bozze\llap{\tt@bozze#1\ }
       \fi(\csname @eq@#1\endcsname)}
\def\eqrefp#1{(\eqref{#1})}
\def\numeriadestra{\let\autoeqno=\autoreqno}
\def\numaeriasinistra{\let\autoeqno=\autoleqno}
\numeriadestra

\newcount\cit@num\global\cit@num=0

\newwrite\file@bibliografia
\newif\if@bibliografia
\@bibliografiafalse

\def\lp@cite{[}
\def\rp@cite{]}
\def\trap@cite#1{\lp@cite #1\rp@cite}
\def\lp@bibl{[}
\def\rp@bibl{]}
\def\trap@bibl#1{\lp@bibl #1\rp@bibl}

\def\refe@renza#1{\if@bibliografia\immediate        
    \write\file@bibliografia{
    \string\item{\trap@bibl{\cref{#1}}}\string
    \bibl@ref{#1}\string\bibl@skip}\fi}

\def\ref@ridefinita#1{\if@bibliografia\immediate\write\file@bibliografia{ 
    \string\item{?? \trap@bibl{\cref{#1}}} ??? tentativo di ridefinire la 
      citazione #1 !!! \string\bibl@skip}\fi}

\def\bibl@ref#1{\se@indefinito{@ref@#1}\immediate
    \write16{ ??? biblitem #1 indefinito !!!}\expandafter\xdef
    \csname @ref@#1\endcsname{ ??}\fi\csname @ref@#1\endcsname}

\def\c@label#1{\global\advance\cit@num by 1\xdef            
   \la@citazione{\the\cit@num}\expandafter
   \xdef\csname @c@#1\endcsname{\la@citazione}}

\def\bibl@skip{\vskip 0truept}


\def\stileincite#1#2{\global\def\lp@cite{#1}\global
    \def\rp@cite{#2}}
\def\stileinbibl#1#2{\global\def\lp@bibl{#1}\global
    \def\rp@bibl{#2}}

\def\citpreset#1{\global\cit@num=#1
    \immediate\write16{ !!! cit-preset = #1 }    }

\def\autobibliografia{\global\@bibliografiatrue\immediate
    \write16{ !!! Genera il file \jobname.BIB}\immediate
    \openout\file@bibliografia=\jobname.bib}

\def\cref#1{\se@indefinito                  
   {@c@#1}\c@label{#1}\refe@renza{#1}\fi\csname @c@#1\endcsname}

\def\cite#1{\trap@cite{\cref{#1}}}                        
\def\upcite#1{$^{\,\trap@cite{\cref{#1}}}$}               
\def\ccite#1#2{\trap@cite{\cref{#1},\cref{#2}}}           
\def\upccite#1#2{$^{\,\trap@cite{\cref{#1},\cref{#2}}}$}  
\def\cccite#1#2#3{\trap@cite{\cref{#1},\cref{#2},\cref{#3}}}          
\def\upcccite#1#2#3{$^{\,\trap@cite{\cref{#1},\cref{#2},\cref{#3}}}$} 
\def\ncite#1#2{\trap@cite{\cref{#1}--\cref{#2}}}          
\def\upncite#1#2{$^{\,\trap@cite{\cref{#1}-\cref{#2}}}$}  

\def\clabel#1{\se@indefinito{@c@#1}\c@label           
    {#1}\refe@renza{#1}\else\c@label{#1}\ref@ridefinita{#1}\fi}

\def\biblskip#1{\def\bibl@skip{\vskip #1}}           

\def\insertbibliografia{\if@bibliografia             
    \immediate\write\file@bibliografia{ }
    \immediate\closeout\file@bibliografia
    \catcode`@=11\input\jobname.bib\catcode`@=12\fi}


\def\commento#1{\relax} 
\def\biblitem#1#2\par{\expandafter\xdef\csname @ref@#1\endcsname{#2}}


\catcode`@=12

\input nuovii.tex

\magnification 1200                  




 
\numerazionedoppia
\autobibliografia

\def\ID{{\bf 1}}
\def\RE{\hbox{\rm{Re}}\,}
\def\Id{\hbox{\rm{Id}}\,}
\def\ddx{{{\rm d}\over{\rm d} x}}
\def\ddtx{{{\rm d}\over{\rm d}\hat x}}
\def\ddz{{{\rm d}\over{\rm d} z}}
\def\thi{{\vartheta_\infty}}
\def\eps{\varepsilon}

\centerline{\bf MONODROMY OF CERTAIN PAINLEVE'--VI TRANSCENDENTS AND}
\vskip 0.3 cm
\centerline{\bf REFLECTION GROUPS}
\vskip 0.3 cm
\centerline{\bf B. Dubrovin and M. Mazzocco}
\vskip 0.2 cm
\centerline{\it International School of Advanced Studies, SISSA-ISAS, Trieste.}
\vskip 0.6 cm

\noindent{\bf Abstract.}\quad We study the global analytic properties of the 
solutions of a particular family of Painlev\'e VI equations with the 
parameters $\beta=\gamma=0$, $\delta={1\over2}$ and $\alpha$ arbitrary. 
We introduce a class of solutions having critical behaviour of algebraic 
type, and completely compute the structure of the analytic continuation of 
these solutions in terms of an auxiliary reflection group in the three 
dimensional space. The analytic continuation is given in terms of an 
action of the braid group on the triples of generators of the reflection 
group. This result is used to classify all the algebraic solutions of our 
Painlev\'e VI equation.

\vskip 0.6 cm\noindent{\sl SISSA preprint no. 149/97/FM, 28 November 1997, 
revised 8 June 1998.}

\vfill\eject

\semiautosez{0}{INTRODUCTION}

In this paper, we will study the structure of the analytic continuation 
of the solutions of the following differential equation
$$
\eqalign{
y_{xx}=&{1\over2}\left({1\over y}+{1\over y-1}+{1\over y-x}\right) y_x^2 -
\left({1\over x}+{1\over x-1}+{1\over y-x}\right)y_x\cr
&+{1\over2}{y(y-1)(y-x)\over x^2(x-1)^2}\left[(2\mu-1)^2+
{x(x-1)\over(y-x)^2}\right],\cr}\eqno{PVI\mu}
$$
in the complex plane, $\mu$ is an arbitrary complex parameter satisfying the 
condition $2\mu\not\in\interi$.

This is a particular case of the general Painlev\'e VI equation (see [Ince]) 
PVI$(\alpha,\beta,\gamma,\delta)$, that depends on four parameters 
$\alpha,\beta,\gamma,\delta$, specified by the following choice of the 
parameters:
$$
\alpha={(2\mu-1)^2\over2},\qquad \beta=\gamma=0\qquad \delta={1\over2}.
$$
The general solution $y(x;c_1,c_2)$ of PVI$(\alpha,\beta,\gamma,\delta)$ 
satisfies the following two important properties (see [Pain]):
\item{1)} The solution $y(x;c_1,c_2)$ can be analytically continued to a 
meromorphic function on the universal covering of 
$\overline\complessi\backslash\{0,1,\infty\}$. 
\item{2)} For generic values of the integration constants $c_1,c_2$ and of 
the parameters $\alpha,\beta,\gamma,\delta$, the solution $y(x;c_1,c_2)$ can 
not be expressed via elementary or classical transcendental functions.

\vskip 0.2 cm
The former claim is the so-called {\it Painlev\'e property}\/  of the equation 
PVI$(\alpha,\beta,\gamma,\delta)$, i.e. its solutions $y(x;c_1,c_2)$ may have 
complicated singularities only at the {\it critical points}\/ of the equation, 
$0,1,\infty$, the position of which does not depend on the choice of 
the particular solution (the so-called {\it fixed singularities}\/), and 
all the other singularities of the solution are poles. Positions of the poles
depend on the integration constants (the so-called 
{\it movable singularities}\/).

All the second order ordinary differential equations of the type:   
$$
y_{xx}={\cal R}(x,y,y_x),
$$
where ${\cal R}$ is rational in $y_x$ and meromorphic in $x$ and $y$ and 
satisfies the Painlev\'e property of absence of movable critical 
singularities, were classified by Painlev\'e and Gambier (see [Pain], 
and [Gamb]). Only six of these equations, which are given in the 
{\it Painlev\'e-Gambier list,}\/ satisfy the property 2), i.e. they can not 
be reduced to known differential equations for elementary and classical 
special functions. The solutions of these equations define some new 
functions, the so-called {\it Painlev\'e transcendents.}\/
PVI$(\alpha,\beta,\gamma,\delta)$ is the most general equation of 
Painlev\'e-Gambier list. Indeed all the others can be obtained from 
PVI$(\alpha,\beta,\gamma,\delta)$ by a confluence procedure (see [Ince] 
\S 14.4).

There are many physical applications of particular solutions of the 
Painlev\'e equations which we do not discuss here. We mention only the 
paper [Tod] where our PVI$\mu$ appears in the problem of the construction 
of self-dual Bianchi-type IX Einstein metrics, and the paper [Dub] where 
the same equation was used to classify the solutions of WDVV equation in 
$2D$-topological field theories.

The name of transcendents could be misleading; indeed, for some particular 
values of $(c_1,c_2,\alpha,\beta,\gamma,\delta)$, the solution 
$y(x;c_1,c_2)$ can be expressed via classical functions. For example 
Picard (see [Pic] and [Ok]) showed that the general solution of 
PVI$(0,0,0,{1\over2})$ can be expressed via elliptic functions, and, more 
recently, Hitchin [Hit] obtained the general solution of 
PVI$({1\over8},{1\over8},{1\over8},{3\over8})$ in terms of the Jacobi 
theta-functions (see also [Man]). Particular examples of classical 
solutions, that can be expressed via hypergeometric functions, of PVI were 
first constructed by Lukashevich [Luka]. A general approach to study the 
classical solutions of PVI was proposed by Okamoto (see [Ok1][Ok2]). One of 
the main tools of this approach is the symmetry group of PVI: the particular 
solutions are those being invariant with respect to some symmetry of PVI. 
The symmetries act in a non trivial way on the space of the parameters 
$(\alpha,\beta,\gamma,\delta)$. Okamoto described the fundamental region 
of the action of this symmetry group and showed that all the classical 
solutions known at that moment, fit into the boundary of this fundamental 
region.

The theory of the classical solutions of the Painlev\'e equations was 
developed by Umemura and Watanabe ([Um], [Um1], [Um2], [Um3], [Wat]); in 
particular, all the one-parameter families of classical solutions of PVI 
were classified in [Wat]. Watanabe also proved that, loosely speaking, 
all the other 
classical solutions of PVI (i.e. not belonging to the one-parameter families) 
can only be given by algebraic functions.

Examples of algebraic solutions were found in [Hit1], for 
PVI$({1\over8},-{1\over8},{1\over2 k^2},{1\over2}-{1\over2 k^2})$, for an 
arbitrary integer $k$. Other examples for PVI$\mu$ were constructed in 
[Dub]. They turn out to be related to the group of symmetries of the 
regular polyhedra in the three dimensional space. Other algebraic solutions
of PVI can be extracted from the recent paper [Seg].

The main aim of our work is to elaborate a tool to classify all the 
algebraic solutions of the Painlev\'e VI equation (for the other five 
Painlev\'e equations, algebraic solutions have been 
classified, see [Kit], [Wat1], [Mur] and [Mur1]). Our idea is very close 
to the main idea of the classical paper of Schwartz (see [Schw]) devoted 
to the classification of the algebraic solutions of the Gauss 
hypergeometric equation. Let $y(x;c_1,c_2)$ be a branch of a solution 
of PVI; its analytic continuation along any closed path $\gamma$ avoiding 
the singularities is a new branch $y(x;c^\gamma_1,c^\gamma_2)$ 
with new integration constants $c^\gamma_1,c^\gamma_2$. Since all the 
singularities of the solution on 
$\overline\complessi\backslash\{0,1,\infty\}$ are poles, the result 
of the analytic continuation depends only on the homotopy class of the 
loop $\gamma$ on the Riemann sphere with three punctures. As a consequence, 
the structure of the analytic continuation is described by an action of 
the fundamental group:
$$
\gamma\in\overline\complessi\backslash\{0,1,\infty\},\quad
\gamma:(c_1,c_2)\to(c^\gamma_1,c^\gamma_2).\autoeqno{in1}
$$
To classify all the algebraic solutions of Painlev\'e VI, all the finite 
orbits of this action must be classified.

Our problem differs from Schwartz's linear analogue, because \eqrefp{in1} is 
not a linear representation but a non-linear action of the fundamental group. 
It is also more involved than the problem of the classification of the 
algebraic solutions of the other five Painlev\'e equations, because the PVI 
is the only equation on the Painlev\'e-Gambier list having a non-abelian 
fundamental group of the complement of the critical locus.

Although the main idea seems to work for the general 
PVI$(\alpha,\beta,\gamma,\delta)$, we managed to completely describe 
the action \eqrefp{in1}, and to solve the problem of the classification 
of the algebraic solutions, only for the particular one-parameter family 
PVI$\mu$. Nevertheless, we have decided to publish these results separately, 
postponing the investigation of the general case to another paper 
(in an effort to keep the paper within a reasonable size, we also postpone 
the study of the resonant case $2\mu\in\interi$, see [Ma]). One of the 
motivations for the present publication is a nice geometrical interpretation 
of the structure of the analytic continuation \eqrefp{in1}, that seems to 
disappear in the general PVI equation.

We now outline the main results and describe the structure of the paper. 
Let us introduce a class of solutions of PVI$\mu$ a-priori containing all 
the algebraic solutions. We say that a branch of a solution $y(x;c_1,c_2)$ 
has {\it critical behaviour of algebraic type,}\/ if there exist three 
real numbers $l_0,l_1,l_\infty$ and three non-zero complex numbers 
$a_0,a_1,a_\infty$, such that
$$
y(x)= \left\{\eqalign{ a_0 x^{l_0} \left(1+{\cal O}(x^\eps)\right),
\qquad\hbox{as}\quad x\rightarrow 0,\cr
1-a_1(1-x)^{l_1} \left(1+{\cal O}((1-x)^\eps)\right),
\qquad\hbox{as}\quad x\rightarrow 1,\cr
 a_\infty x^{1-l_\infty} \left(1+{\cal O}(x^{-\eps})\right),
\qquad\hbox{as}\quad x\rightarrow \infty,\cr
}\right.\autoeqno{B1}
$$
where $\eps>0$ is small enough. We show that there exists a three-parameter 
family of solutions of PVI$\mu$ with critical behaviour of algebraic type, 
where $\mu$ itself is a function of  $l_0,l_1,l_\infty$. Of course, for an 
algebraic solution, the indices $l_0,l_1,l_\infty$ must be rational.

It turns out that the three-parameter family of solutions \eqrefp{B1} is 
closed under the analytic continuation \eqrefp{in1}, if and only if $\mu$ 
is real. One of our main results is the parameterization of the solutions 
\eqrefp{B1} by ordered triples of planes in the three dimensional 
{\it Euclidean space}\/(see Section 1.4). In particular, the indices 
$l_0,l_1,l_\infty$ are related to the angles $\pi r_0,\pi r_1,\pi r_\infty$ 
between the planes:
$$
l_i=\left\{\eqalign{
2r_i\quad\hbox{if}\quad 0<r_i\leq{1\over2}\cr
2-2r_i\quad\hbox{if}\quad{1\over2}\leq r_i<1\cr
}\right.
\qquad i=0,1,\infty,
$$
and the parameter $\mu$ is determined within the ambiguity 
$\mu\mapsto\pm\mu+n$, $n\in\interi$, by the equation:
$$
\sin^2\pi\mu = \cos^2\pi r_0+ \cos^2\pi r_1+ \cos^2\pi r_\infty+ 
2 \cos\pi r_0\cos\pi r_1\cos\pi r_\infty.
$$
This ambiguity and the one due to the reordering of the planes can be 
absorbed by the symmetries of PVI$\mu$ described in Section 1.2.

We compute the analytic continuation \eqrefp{in1} in terms of some 
elementary operations on the planes. This computation leads to prove that, 
for an algebraic solution of PVI$\mu$, the reflections in the planes must 
generate the symmetry group of a regular polyhedron in $\reali^3$. Another  
result of this paper is the classification of all algebraic 
solutions of PVI$\mu$. They are in one-to-one correspondence, modulo the 
symmetries of the equations described in Section 1.2, with the reciprocal 
pairs of the three-dimensional regular polyhedra and star-polyhedra (the 
description of the star-polyhedra can be found in [Cox]). The solutions 
corresponding to the regular tetrahedron, cube and icosahedron are the 
ones obtained in [Dub] using the theory of polynomial Frobenius manifolds. 
The solutions corresponding to the regular great 
icosahedron, and regular great dodecahedron are new. Our method not only
allows to classify the solutions, but also to obtain the explicit formulae, 
as we do in Section 2.4.

The main tool to obtain these results is the isomonodromy deformation 
method (see [Fuchs], [Sch] and [JMU], [ItN], [FlN]). The Painlev\'e VI is 
represented as the equation of isomonodromy deformation of the auxiliary 
Fuchsian system
$$
{{\rm d} Y\over{\rm d} z}=
\left({A_0\over z}+{A_1\over z-1}+{A_x\over z-x}\right) Y.
\autoeqno{in2}
$$
For PVI$\mu$, the $2\times2$ matrices $A_0,A_1,A_x$ are nilpotent and
$$
A_0+A_1+A_x=\pmatrix{-\mu&0\cr 0&\mu\cr}.
$$
The entries of the matrices $A_i$ are complicated expressions of $x,y,y_x$ 
and of some quadrature $\int R(x,y){\rm d}x$. The monodromy of 
\eqrefp{in2} remains constant if and only if $y=y(x)$ satisfies PVI. Thus, 
the solutions of PVI$\mu$ are parameterized by the monodromy data of the 
Fuchsian system \eqrefp{in2} (see Section 1.1). 
In section 1.2, we compute the structure of the analytic continuation in 
terms of the monodromy data. On this basis, in Section 1.3, we classify 
all the monodromy data of the algebraic solutions of PVI$\mu$. To this end, 
we classify all the rational solutions of certain trigonometric equations 
using the method of a paper by Gordan (see [Gor]). In Section 1.4, we 
parameterize the monodromy data of PVI$\mu$ by ordered triples of planes 
in the three-dimensional space, considered modulo rotations.
The structure of the analytic continuation of the solutions of PVI$\mu$ is 
reformulated in terms of a certain action of the braid group $B_3$ on the 
triples of planes. The group $G$ generated by the reflections with respect 
to the planes remains unchanged. For the algebraic solutions, the group 
$G$ turns out to coincide with the symmetry group of one of the regular 
polyhedra in the three-dimensional Euclidean space. We also give another 
proof, suggested by E. Vinberg, of this result, and we establish 
that the class of solutions of PVI$\mu$ parameterized by triples of planes 
in the three dimensional Euclidean space is invariant with respect to the 
analytic continuation.

In the second part of the paper, we identify this class of solutions of 
PVI$\mu$ with the class of solutions having critical behaviour of algebraic 
type \eqrefp{B1}. In Section 2.1, we prove that the solution $y(x)$ of the 
form \eqrefp{B1}, for a fixed value of $\mu$, is uniquely determined by its 
asymptotic behaviour near one of the critical points, i.e. by any of the 
pairs $(a_0,l_0)$, $(a_1,l_1)$,$(a_\infty,l_\infty)$. In particular, we 
prove that, for an algebraic solution of PVI$\mu$, the indices 
$l_0,l_1,l_\infty$ must satisfy:
$$
0<l_i\leq1,\qquad i=0,1,\infty.
$$
To derive the {\it connection formulae}\/ establishing the relations 
between these pairs, we use (see Section 2.3) the properly adapted method 
of Jimbo (see [Jim]). This method allows to express the monodromy data of 
the auxiliary Fuchsian system \eqrefp{in2} in terms of the parameters 
$(a_0,l_0)$, $(a_1,l_1)$ or $(a_\infty,l_\infty)$.  For convenience of 
the reader, and because of some differences between the assumptions of 
Jimbo's work and ours, we give a complete derivation of the connection 
formulae in Section 2.2. Using the results of the Sections 1.3 and 1.4, 
we complete the computation of the critical behaviour \eqrefp{B1} for all 
the branches of the analytic continuation of the solution. The result of 
this computation is used in Section 2.4 to obtain the explicit formulae 
for the algebraic solutions of PVI$\mu$.

\vskip 0.2 cm
\noindent{\bf Remark 0.1.}\quad The resulting classification of the algebraic 
solutions of PVI$\mu$ is in striking similarity to the Schwartz's 
classification (see [Schw]) of the algebraic solutions of the 
hypergeometric equation. According to Schwartz, the algebraic solutions 
of the hypergeometric equation, considered modulo contiguity 
transformations, are of fifteen types (the first type consists of an 
infinite sequence of solutions). The rows $(2-15)$ of Schwartz's 
list (see, for example, the table in Section 2.7.2 of [Bat]) correspond 
to the triples of generating reflections of the symmetry groups of regular 
polyhedra in the three-dimensional Euclidean space (we are grateful to E. 
Vinberg for bringing this point to our attention). The parameter 
$\lambda,\mu,\nu$ of the hypergeometric equation shown in the table are 
just the angles between the mirrors of the reflections, divided by $\pi$.

According to our classification, the algebraic solutions of PVI$\mu$, 
considered modulo symmetries, are in one-to-one correspondence to the 
classes of equivalence of the triples of generating reflections in the 
symmetry groups of regular polyhedra. The equivalence is defined by an 
action of the braid group $B_3$ on the triples and by orthogonal 
transformations. We find that in the groups 
$G=W(A_3)$ and $G=W(B_3)$, the symmetry groups of respectively the regular 
tetrahedron and of the cube or regular octahedron, there is only one 
equivalence class of triples of generating reflections; these are given  
respectively by the rows $(2,3)$ and by $(4,5)$ of Schwartz's table. 
In the group $W(H_3)$ of symmetries of regular icosahedron or regular 
dodecahedron, there are three equivalence classes of triples of  
reflections which are given  respectively by the rows $(6,8,13)$,  
$(11,14,15)$ and $(7,9,10,12)$ of the Schwartz's table and correspond 
to icosahedron, great icosahedron and great dodecahedron (or to their 
reciprocal pairs, see [Cox]). To establish the correspondence, we 
associate a {\it standard}\/ system of generating reflections to a regular 
polyhedron in the following way: let $H$ be the center of the polyhedron, 
$O$ the center of a face, $P$ a vertex of this face and $Q$ the center of 
an edge of the same face through the vertex $P$. Then the reflections with 
respect to the planes $HOP$, $HOQ$ and $HPQ$ are the standard system of 
generators. Our five algebraic solutions correspond to the classes of 
equivalence of the standard systems of generators obtained by this 
contruction applied to tetrahedron, cube, icosahedron, great icosahedron, 
great dodecahedron.

Summarizing, we see that the list of all the algebraic solutions of 
PVI$\mu$ is obtained by folding of the list of Schwartz modulo the action 
of the braid group. This relation between the algebraic solutions of 
PVI$\mu$ and the algebraic hypergeometric functions seems to be surprising 
also from the point of view of the results of Watanabe (see [Wat]) who 
classified all the one-parameter families of classical solutions of 
PVI$\mu$ (essentially, all of them are given by hypergeometric functions). 
Using these results, one can  easily check that our algebraic solutions 
do not belong to any of the one-parameter families of classical solutions 
of PVI$\mu$.

\vfill\eject
\noindent{\bf  Acknowledgments}
\vskip 0.3 cm
The authors are indebted to E. Vinberg for the elegant proof of theorem 1.8. 
We thank A. Akhmedov for a simple proof of the Algebraic Lemma, Section 1.4. 
We thank also R. Conte for drawing our attention to the classical work of 
Picard (see [Pic]) and V. Sokolov and F. Zanolin for useful discussions.

%





\semiautosez{1}{1. STRUCTURE OF ANALYTIC CONTINUATION AND ALGEBRAIC 
SOLUTIONS OF PVI${\bf \mu}$}

\noindent{\bf  1.1. Painlev\'e VI equation as isomonodromy deformation 
equation.}
\vskip 0.3 cm

In this section we show how the PVI$\mu$ equation can be reduced to the 
isomonodromy deformation equation of an auxiliary Fuchsian system 
(see [Sch], [JMU]); moreover we describe the parameterization, 
essentially due to Schlesinger (see. [Sch]), of the solutions of the PVI$\mu$ 
equation by the monodromy data of such Fuchsian system. 

\vskip 0.3 cm
\noindent{\bf 1.1.1. An auxiliary Fuchsian system and its monodromy 
data.}\quad In this subsection, we introduce an auxiliary Fuchsian system,
define its monodromy and connection matrices, and extabilish the 
correspondence between monodromy matrices and coefficients of the Fuchsian 
system for a given set of poles.

Let us consider the following Fuchsian system with four regular 
singularities at $u_1$, $u_2$, $u_3$ and $\infty$:
$$
\ddz Y={\cal A}(z)Y,
\qquad\qquad z\in\overline\complessi\backslash\{u_1,u_2,u_3,\infty\}
\autoeqno{N1}
$$
where ${\cal A}(z)$ is a matrix-valued function:
$$
{\cal A}(z)=
{{\cal A}_1\over z-u_1} + {{\cal A}_2\over z-u_2} + {{\cal A}_3\over z-u_3},
$$
${\cal A}_i$ being $2\times 2$ matrices independent on $z$, and 
$u_1,u_2,u_3$ being pairwise distinct complex numbers. We assume that 
the matrices ${\cal A}_i$ satisfy the following conditions:
$$
{\cal A}_i^2=0\quad\hbox{and}\quad
{\cal A}_\infty:= -{\cal A}_1-{\cal A}_2-{\cal A}_3=\pmatrix{
\mu & 0\cr 0 & -\mu\cr}.\autoeqno{N1.3}
$$
Indeed, we will see in the latter part of this section that this choice
corresponds to the particular case PVI$\mu$ of the Painlev\'e VI equation. 
In this paper we consider the {\it non-resonant}\/ case $2\mu\not\in\interi$. 

The solution $Y(z)$ of the system \eqrefp{N1} is a multi-valued analytic 
function in the punctured Riemann sphere 
$\complessi\backslash\{u_1,u_2,u_3\}$ and its multivaluedness 
is described by the so-called {\it monodromy matrices.}\/  

Let us briefly recall the definition of the monodromy matrices of the 
Fuchsian system \eqrefp{N1}. First, we fix a basis 
$\gamma_1,\gamma_2,\gamma_3$ of loops in the fundamental group, with base 
point at $\infty$, of the punctured Riemann sphere 
$$
\pi_1\left(\overline\complessi\backslash\{u_1,u_2,u_3,\infty\},\hat x\right),
$$ 
and  a fundamental matrix for the system \eqrefp{N1}.
To fix the basis of the loops, we first perform some cuts between the 
singularities, namely we cut three parallel segments $\pi_i$ between the 
point at infinity and each $u_i$; the segments $\pi_i$ are ordered 
according to the order of the points $u_1,u_2,u_3$, as in the figure 1. 
Take $\gamma_i$ to be a simple closed curve starting and finishing at 
infinity, going around $u_i$ in positive direction ($\gamma_i$ is oriented 
counter-clockwise, $u_i$ lies inside, while the other singular points lie 
outside) and not crossing the cuts $\pi_i$. Near $\infty$, we take every 
loop $\gamma_i$ close to the cut $\pi_i$ as in the figure 1.

\midinsert
\centerline{\psfig{file=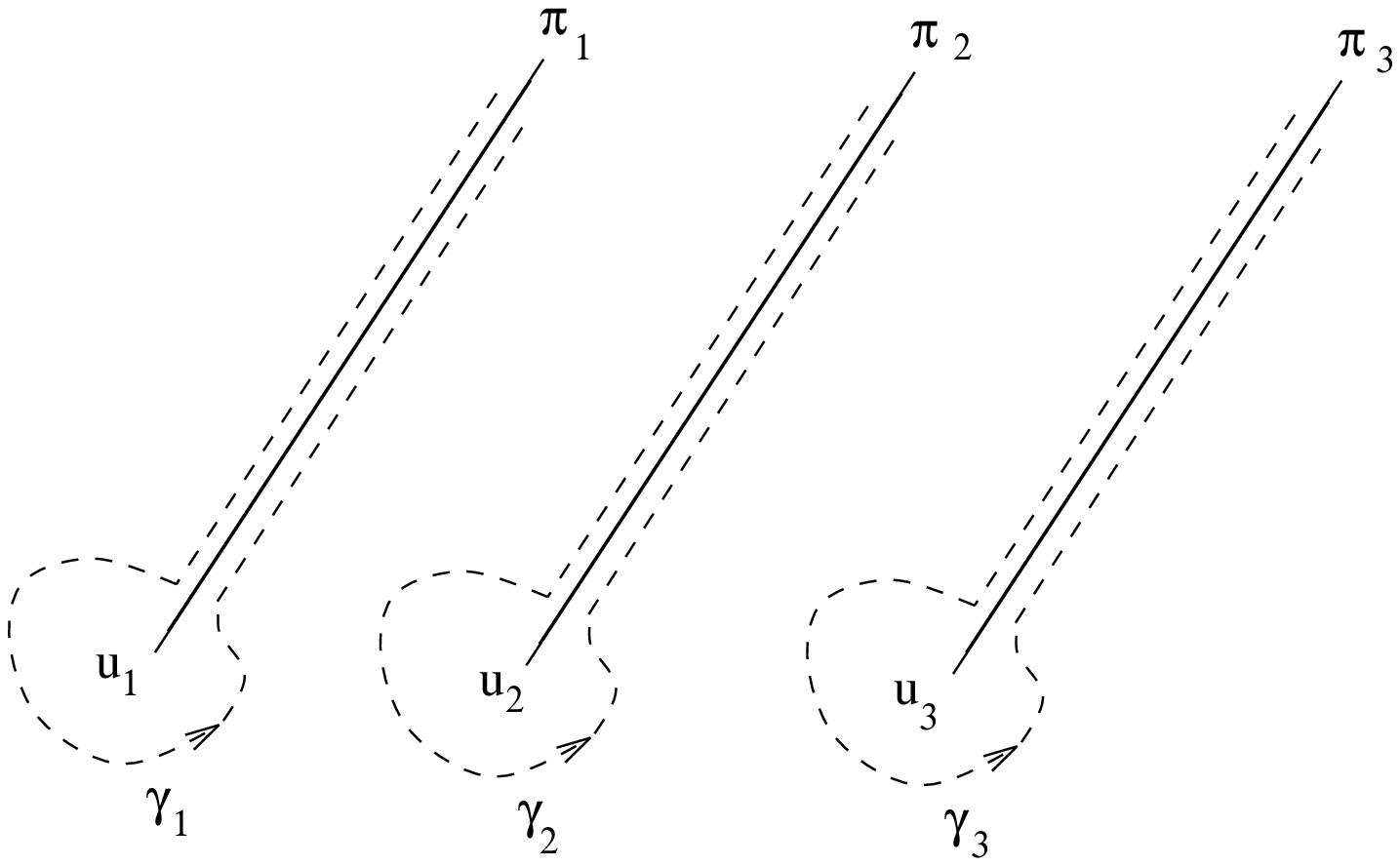,height=4cm}}
\vskip 0.6 cm
\centerline{{\bf Fig.1.} The cuts $\pi_i$ between the singularities $u_i$ 
and the oriented loops $\gamma_i$.}
\endinsert

\noindent Now, we fix the fundamental matrix $Y_\infty(z)$ of the system 
\eqrefp{N1} in such a way that
$$
Y_\infty(z)=\left(\ID +{\cal O}({1\over z})\right) 
\pmatrix{z^{-\mu}&0\cr 0& z^\mu\cr},\quad\hbox{as}\quad z\rightarrow\infty,
\autoeqno{N1.5}
$$
where $z^\mu:=e^{\mu\log z}$, with the choice of the principal branch of 
the logarithm, with the branch-cut along the common direction of the cuts
$\pi_1,\pi_2,\pi_3$. Such a fundamental matrix $Y_\infty(z)$ exists and, due 
to the non-resonance condition, it is uniquely determined. It can be 
analytically continued to an analytic function on the universal covering 
of $\overline\complessi\backslash\{u_1,u_2,u_3,\infty\}$. For any element 
$\gamma\in\pi_1
\left(\overline\complessi\backslash\{u_1,u_2,u_3,\infty\},\infty\right)$
we denote the result of the analytic continuation of $Y_\infty(z)$ along 
the loop $\gamma$ by $\gamma[Y_\infty(z)]$. Since $\gamma[Y_\infty(z)]$ 
and $Y_\infty(z)$ are two fundamental matrices in the neighborhood of 
infinity, they must be related by the following relation:
$$
\gamma[Y_\infty(z)]= Y_\infty(z) M_\gamma
$$
for some constant invertible $2\times 2$ matrix $M_\gamma$ depending only 
on the homotopy class of $\gamma$. Particularly, the matrix 
$M_\infty:=M_{\gamma_\infty}$, $\gamma_\infty$ being a simple loop around 
infinity in the clock-wise direction, is given by:
$$
M_\infty=\pmatrix{
\exp (2 i \pi \mu) & 0 \cr 0 &\exp (- 2 i  \pi \mu) \cr}.\autoeqno{N9}
$$
The resulting {\it monodromy representation}\/ is an anti-homomorphism:
$$
\matrix{
\pi_1\left(\overline\complessi\backslash\{u_1,u_2,u_3,\infty\},\infty\right)
&\rightarrow & SL_2(\complessi)\cr
\gamma&\mapsto & M_\gamma\cr}\autoeqno{N2.2}
$$
$$
M_{\gamma\tilde\gamma}= M_{\tilde\gamma}  M_{\gamma}.\autoeqno{N2}
$$ 
The images $M_i:=M_{\gamma_i}$ of the generators $\gamma_i$, $i=1,2,3$ of 
the fundamental group, are called {\it the monodromy matrices}\/ of 
the Fuchsian system \eqrefp{N1}. They generate the {\it monodromy group of 
the system,}\/ i.e. the image of the representation \eqrefp{N2.2}.
Moreover, due to the fact that, in our 
particular case, the ${\cal A}_i$ are nilpotent, satisfy the following 
relations:
$$
\det(M_i)=1, \quad {\rm Tr}(M_i)=2,\quad\hbox{for}\quad i=1,2,3,\autoeqno{N5}
$$
with $M_i=\ID$ if and only if ${\cal A}_i=0$. Moreover, since the loop 
$(\gamma_1 \gamma_2 \gamma_3)^{-1}$ is homotopic to $\gamma_\infty$, 
the following relation holds:
$$
M_\infty M_3 M_2 M_1=\ID.\autoeqno{N6}
$$
A simultaneous conjugation ${\cal A}_i\mapsto D^{-1} {\cal A}_i D$, 
$i=1,2,3$ of the coefficients ${\cal A}_i$ of the Fuchsian system 
\eqrefp{N1} by a diagonal matrix $D$, implies the same conjugation of the 
monodromy matrices $M_\gamma\mapsto D^{-1} M_\gamma D$, for any 
$\gamma\in\pi_1\left(\overline\complessi\backslash\{u_1,u_2,u_3,\infty\},
\infty\right)$.

We now recall the definition of the {\it connection matrices.}\/ Let us 
assume that $M_i\neq \ID$, or equivalently ${\cal A}_i\neq 0$, for every 
$i=1,2,3$.
We choose the fundamental matrices $Y_i(z)$ of the system \eqrefp{N1}, such 
that:
$$
Y_i=G_i\left(\ID +{\cal O}(z-u_i)\right)(z-u_i)^{J},\quad\hbox{as}\quad 
z\rightarrow u_i, \autoeqno{N6.1}
$$
where $J$ is the Jordan normal form of ${\cal A}_i$, namely 
$J=\pmatrix{0&1\cr 0 & 0\cr}$, the invertible matrix $G_i$ is defined by 
${\cal A}_i=G_i J G_i^{-1}$, and the choice of the branch of $\log(z-u_i)$ 
needed in the definition of
$$
(z-u_i)^{J}=\pmatrix{1&\log(z-u_i)\cr 0 &1\cr}
$$
is similar to the one above.
The fundamental matrix  $Y_i(z)$ is uniquely determined up to the ambiguity:
$$
Y_i(z)\mapsto Y_i(z) R_i,
$$
where $R_i$ is any matrix commuting with $J$.

Continuing, along, say, the right-hand-side of the cut $\pi_i$, the solution 
$Y_\infty$ to a neighborhood of $u_i$, we obtain another fundamental matrix 
around $u_i$, that must be related to $Y_i(z)$ by:
$$
Y_\infty(z)=Y_i(z)C_i,\autoeqno{N3}
$$
for some invertible matrix $C_i$. The matrices $C_1,C_2, C_3$  are called 
{\it connection matrices,}\/ and are related to the monodromy matrices as 
follows:
$$
M_i=C_i^{-1} \exp(2\pi i J) C_i,\qquad i=1,2,3.\autoeqno{N4}
$$

\proclaim Lemma 1.1. Given three matrices $M_1$, $M_2$, $M_3$, 
$M_i\neq \ID$ for every $i=1,2,3$, satisfying the relations \eqrefp{N5} 
and \eqrefp{N6}, then
\item{i)} there exist three matrices $C_1,C_2, C_3$ satisfying the 
\eqrefp{N4}. Moreover they are uniquely determined by the matrices 
$M_1$, $M_2$, $M_3$, up to the ambiguity $C_i\mapsto R_i^{-1} C_i$,
where $R_i J=J R_i$, for $i=1,2,3$. 
\item{ii)} If the matrices $M_1$, $M_2$, $M_3$ are the monodromy matrices 
of a Fuchsian system of the form \eqrefp{N1}, then any triple $C_1,C_2,C_3$
satisfying \eqrefp{N4} can be realized as the connection matrices of 
the Fuchsian system itself. 

\noindent Proof. i) By the \eqrefp{N5}, the monodromy 
matrices have all the eigenvalues equal to one; moreover they can be 
reduced to the Jordan normal form because $M_i\neq \ID$. Namely there 
exists a matrix $\tilde C_i$ such that:
$$
M_i = {\tilde C_i}^{-1} \pmatrix{1&1\cr 0 &1\cr} \tilde C_i.
$$ 
Taking
$$
C_i =  \pmatrix{1&2\pi i\cr 0 &1\cr} \tilde C_i
$$
we obtain the needed matrix.
Two such matrices $C_i$ and $C'_i$ give the same matrix $M_i$ if and only 
if $C_i^{-1} C'_i$ commutes with $J$, namely if and only if they are 
related by $C_i = R_i^{-1} C'_i$. ii) Let us now assume that 
$C'_1,C'_2, C'_3$ are the connection matrices of a Fuchsian system 
of the form \eqrefp{N1}, with monodromy matrices $M_1$, $M_2$, $M_3$; 
id est $Y_\infty(z)=Y'_i(z)C'_i$, $i=1,2,3$, for some choice of the 
solutions $Y'_1$, $Y'_2$ and $Y'_3$ of the form \eqrefp{N6.1}. We have 
$$
M_i=(C'_i)^{-1} \exp(2\pi i J)C'_i=
C_i^{-1} \exp(2\pi i J)C_i, \quad i=1,2,3.
$$
So the matrices $R_i=C'_i C_i^{-1}$ must commute with $J$ and $C_1$, $C_2$, 
$C_3$ are the connection matrices with respect to the new solutions 
$Y_i(z)=Y'_i(z)R_i$. {\hfill QED}
\vskip 0.3 cm

Now, we state the result about the correspondence between monodromy 
data and coefficients of the Fuchsian system, for a given set of poles:

\proclaim Lemma 1.2. Two Fuchsian systems \eqrefp{N1} with the same 
poles $u_1$, $u_2$ and $u_3$, and the same value of $\mu$, coincide 
if and only if they have the same monodromy matrices $M_1$, $M_2$, 
$M_3$, with respect to the same basis of the loops $\gamma_1$, 
$\gamma_2$ and $\gamma_3$.

\noindent Proof. Let $Y_\infty^{(1)}(z)$ and 
$Y_\infty^{(2)}(z)$ be the fundamental matrices of the form 
\eqrefp{N1.5} of the two Fuchsian systems. Let us consider the 
following matrix:
$$
Y(z):= Y_\infty^{(2)}(z)Y_\infty^{(1)}(z)^{-1}.
$$
$Y(z)$ is an analytic function around infinity: 
$$
Y(z)=1+{\cal O}\left({1\over z}\right),\quad\hbox{as}\, z\rightarrow\infty.
$$
Since the monodromy matrices coincide, $Y(z)$ is a single valued function 
on $\overline\complessi\backslash\{u_1,u_2,u_3\}$. Let us prove that $Y(z)$ 
is analytic also at the points $u_i$. 
Due to Lemma 1.1, we can choose the fundamental matrices $Y_i^{(1)}(z)$ 
and $Y_i^{(2)}(z)$ in such a way that
$$
Y_\infty^{(1),(2)}(z)= Y_i^{(1),(2)}(z)C_i \quad i=1,2,3.
$$
with the same connection matrices $C_i$. Then near the point $u_i$,
$$
Y(z)=G_i^{(2)}\left(\ID+{\cal O}(z-u_i) \right)
\left[G_i^{(1)}\left(\ID+{\cal O}(z-u_i) \right)\right]^{-1}.
$$
This proves that $Y(z)$ is an analytic function on all 
$\overline\complessi$ and then, by the Liouville theorem
$$
Y(z)=\ID,
$$
and the two Fuchsian systems must coincide.

\proclaim Corollary 1.1. Two Fuchsian systems \eqrefp{N1} with the same 
poles $u_1$, $u_2$ and $u_3$, and the same value of $\mu$, are conjugated
$$
{\cal A}_i^{(1)}=D^{-1} {\cal A}_i^{(2)} D,\qquad i=1,2,3,
$$
with a diagonal matrix $D$, if and only if their monodromy matrices 
$M_i^{(1)}$ and $M_i^{(2)}$, with respect to the same basis of the 
loops $\gamma_1$, $\gamma_2$ and $\gamma_3$, are conjugated:
$$
M_i^{(1)}=D^{-1} M_i^{(2)} D,\qquad i=1,2,3.
$$

\vskip 0.3 cm
\noindent{\bf 1.1.2. The isomonodromy deformations of the Fuchsian system 
\eqrefp{N1} and the Painlev\'e equation PVI$\mu$.}\quad
We now want to deform the poles of the Fuchsian system keeping the 
monodromy fixed. The theory of these deformations is described by the 
following two results:

\proclaim Theorem 1.1. Let $M_1$, $M_2$, $M_3$ be the monodromy matrices of 
the Fuchsian system:
$$
\ddz Y^0=
\left({{\cal A}^0_1\over z-u^0_1}+{{\cal A}^0_2\over z-u^0_2}+
{{\cal A}^0_3\over z-u^0_3}\right)Y^0,
\autoeqno{N7}
$$
of the above form \eqrefp{N1.3}, with pairwise distinct poles $u_i^0$, and 
with respect to some basis $\gamma_1,\gamma_2,\gamma_3$ of the loops 
in $\pi_1\left(\overline\complessi\backslash\{u^0_1,u^0_2,u^0_3,\infty\},
\infty\right)$. Then there exists a neighborhood $U\subset\complessi^3$ 
of the point $u^0=(u^0_1,u^0_2,u^0_3)$ such that, for any 
$u= (u_1,u_2,u_3)\in U$, there exists a unique triple ${\cal A}_1(u)$, 
${\cal A}_2(u)$, ${\cal A}_3(u)$ of analytic matrix valued functions 
such that:
$$
{\cal A}_i(u^0)={\cal A}_i^0,\quad i=1,2,3,
$$
and the monodromy matrices of the Fuchsian system
$$
\ddz Y= A(z;u) Y=
\left({{\cal A}_1(u)\over z-u_1}+{{\cal A}_2(u)\over z-u_2}+
{{\cal A}_3(u)\over z-u_3}\right)Y, \autoeqno{N8}
$$
with respect to the same basis\footnote{${}^{1}$}{Observe that the 
basis $\gamma_1,\gamma_2,\gamma_3$ of 
$\pi_1\left(\overline\complessi\backslash\{u_1,u_2,u_3,\infty\},
\infty\right)$ varies continuously with small variations of $u_1,u_2,u_3$. 
This new basis is homotopic to the initial one, so we can identify them.}
$\gamma_1,\gamma_2,\gamma_3$ of the loops, coincide with the given $M_1$, 
$M_2$, $M_3$.
The matrices ${\cal A}_i(u)$ are the solutions of the Cauchy problem with 
the initial data ${\cal A}_i^0$ for the following Schlesinger equations:
$$
{\partial\over\partial u_j}{\cal A}_i= 
{[{\cal A}_i,{\cal A}_j]\over u_i-u_j},\quad
{\partial\over\partial u_i}{\cal A}_i= 
-\sum_{j\neq i}{[{\cal A}_i,{\cal A}_j]\over u_i-u_j}. 
\autoeqno{N10}
$$
The solution $Y_\infty^0(z)$ of \eqrefp{N7} of the form \eqrefp{N1.5} 
can be uniquely continued, for $z\neq u_i$ $i=1,2,3$, to an analytic function
$$
Y_\infty(z,u),\quad u\in U,
$$
such that
$$
Y_\infty(z,u^0)=Y_\infty^0(z).
$$
This continuation is the local solution of the Cauchy problem with 
the initial data $Y_\infty^0$ for the following system that is 
compatible to the system \eqrefp{N8}:
$$
{\partial\over\partial u_i} Y = -{{\cal A}_i(u)\over z-u_i} Y.
$$
Moreover the functions ${\cal A}_i(u)$ and $Y_\infty(z,u)$ can be 
continued analytically to global meromorphic functions on the universal 
coverings of
$$
\complessi^3\backslash\{diags\}:=
\left\{(u_1,u_2,u_3)\in\complessi^3\,|\,u_i\neq u_j\,\hbox{for}\, 
i\neq j\right\},
$$
and
$$
\left\{(z,u_1,u_2,u_3)\in\complessi^4\,|\,u_i\neq u_j\,\hbox{for}\, 
i\neq j\,\hbox{and}\,z\neq u_i,\, i=1,2,3\right\},
$$
respectively.

The proof can be found, for example, in [Mal], [Miwa], 
[Sib]. We recall the theorem of solvability of the inverse problem of the 
monodromy (see [Dek]):

\proclaim Theorem 1.2. Given three arbitrary matrices, satisfying 
\eqrefp{N5} and \eqrefp{N6}, with $M_\infty$ of the form \eqrefp{N9}, and 
given a point $u^0=(u^0_1,u^0_2,u^0_3)\in \complessi^3\backslash\{diags\}$, 
for any neighborhood $U$ of $u^0$, there exist $(u_1,u_2,u_3)\in U$ and a 
Fuchsian system of the form \eqrefp{N1}, with the given monodromy matrices, 
the given $\mu$ and with poles in $u_1,u_2,u_3$.

\noindent{\bf Remark 1.1.}\quad Fuchsian systems of the form \eqrefp{N1}, 
with coefficients ${\cal A}_i$ satisfying \eqrefp{N1.3}, depend on four 
parameters, one of them being $\mu$. The triples of the monodromy matrices 
satisfying \eqrefp{N5} and \eqrefp{N6}, with $M_\infty$ 
of the form \eqrefp{N9}, depend on four parameters too. 
Loosely speaking, Theorems 1.1 and 1.2 claim that, not only the monodromy 
matrices are first integrals for the equations of isomonodromy deformation 
\eqrefp{N10}, but they provide a full system of first integrals for such 
equations. We denote ${\cal A}(u_1,u_2,u_3;M_1,M_2,M_3)$ the solution of the 
Schlesinger equations locally uniquelly determined by the triple of monodromy 
matrices $(M_1,M_2,M_3)$.
\vskip 0.2 cm

All the above arguments remain valid for a general $2\times 2$ Fuchsian 
system, provided the non-resonancy condition of the eigenvalues of 
${\cal A}_i$ and ${\cal A}_\infty$.
\vskip 0.2 cm

\noindent{\bf Remark 1.2.}\quad We observe that the isomonodromy deformations 
equations preserve the connection matrices $C_i$ too. This follows from 
Lemma 1.1.
\vskip 0.3 cm

\noindent{\bf 1.1.3. Reduction to the PVI$\mu$ equation.}\quad
Let us now explain, following [JMU], how to reduce the Schlesinger equations 
\eqrefp{N10} to the PVI$\mu$ equation.
The Schlesinger equations are invariant with respect to the gauge 
transformations of the form:
$$
{\cal A}_i\mapsto D^{-1} {\cal A}_i D,\quad i=1,2,3,\quad\hbox{for any $D$ 
diagonal matrix}.
$$
First of all we have to factor out such gauge transformations; to this aim, 
we introduce two coordinates $(p,q)$ on the quotient of the space of the 
matrices satisfying \eqrefp{N1.3} with respect to the equivalence relation
$$
{\cal A}_i\sim D^{-1} {\cal A}_i D,\quad i=1,2,3,\quad\hbox{for any $D$ 
diagonal matrix}.
\autoeqno{N12}
$$
The coordinates $(p,q)$ are defined as follows: $q$ is the root of the 
following linear equation:
$$
[{\cal A}(q;u_1,u_2,u_3)]_{12}=0, 
$$
and $p$ is given by:
$$
p=[{\cal A}(q;u_1,u_2,u_3)]_{11},
$$
where ${\cal A}(z;u_1,u_2,u_3)$ is given in \eqrefp{N8}. The matrices 
${\cal A}_i$ are expressed rationally in terms of the coordinates 
$(p,q)$ and an auxiliary coordinate $k$, coming from the gauge freedom 
\eqrefp{N12}:
$$
\eqalign{
\left({\cal A}_i\right)_{11}&=-\left({\cal A}_i\right)_{22}
={q-u_i\over2\mu P'(u_i)}\left[P(q)p^2+2\mu {P(q)\over q-u_i} p+
\mu^2(q+2u_i-\sum_ju_j)\right],\cr
\left({\cal A}_i\right)_{12}&=-\mu k {q-u_i\over P'(u_i)},\cr
\left({\cal A}_i\right)_{21}&= k^{-1} {q-u_i\over4\mu^3 P'(u_i)}
\left[P(q)p^2+2\mu {P(q)\over q-u_i} p
+\mu^2(q+2u_i-\sum_ju_j)\right]^2,\cr}
\autoeqno{N12.5}
$$
for $i=1,2,3$, where $P(z) = (z-u_1)(z-u_2)(z-u_3)$ and 
$P'(z)={{\rm d}P\over{\rm d}z}$. 
The Schlesinger equations in these coordinates reduce to:
$$
\eqalign{
{\partial q\over\partial u_i} &=
{P(q)\over P'(u_i)}\left[2 p + {1\over q-u_i}\right]\cr
{\partial p\over\partial u_i} &=
-{P'(q) p^2 +(2q+u_i-\sum_ju_j)p+\mu(1-\mu)\over P'(u_i)},
\cr}\autoeqno{N13}
$$
for $i=1,2,3$. The system of the {\it reduced Schlesinger equations}\/ 
\eqrefp{N13} is invariant under the transformations of the form
$$
u_i\mapsto a u_i + b,\qquad
q\mapsto a q + b,\qquad
p\mapsto {p\over a}, \qquad\forall a,b\in\complessi,\quad a\neq0.
$$
We introduce the following new invariant variables:
$$
\eqalign{
x&={u_2-u_1\over u_3-u_1},\cr
y&={q-u_1\over u_3-u_1};\cr}\autoeqno{N13.5}
$$
the system \eqrefp{N13}, expressed in the these new variables, reduces to 
the PVI$\mu$ equation for $y(x)$.

\vskip 0.3 cm
\noindent{\bf Remark 1.3.}\quad The system \eqrefp{N13} admits the following 
{\it singular solutions}\/ (see [Ok1] and [Wat]):
$$
q\equiv u_i\quad\hbox{for some $i$},
$$
and $p$, in the variable $x$, can be expressed via Gauss hypergeometric 
functions (see [Ok1]). Moreover the monodromy group of the system 
\eqrefp{N1} reduces to the monodromy group of the Gauss hypergeometric 
equation, namely the following lemma holds true:

\proclaim Lemma 1.3. The solutions of the full Schlesinger equations, 
corresponding to the solution $q\equiv u_i$, for some $i$, have the form:
$$
A_i(u)\equiv0,\quad\hbox{and for}\, j\neq i\quad 
A_j(u)=D(u)^{-1} A_j^0 D(u),
$$
where $D(u)$ is a diagonal matrix depending on $u$, and $A_j^0$ is a 
constant matrix. The monodromy matrix $M_i$ of the corresponding Fuchsian 
system turns out to be the identity. Conversely, if one of the monodromy 
matrices $M_i$ is the identity, $M_i=\ID$, then the solution of 
\eqrefp{N13} is degenerate.

\noindent Proof. The matrix $A_i$, for $q\equiv u_i$, is identically $0$, 
thanks to \eqrefp{N12.5}. Having $A_i\equiv0$, $M_i$ is $\ID$. Conversely, 
if $M_i=\ID$, then $A_i\equiv0$. Solving the Schlesinger equations 
\eqrefp{N8}, we obtain $q\equiv u_i$, and the equation for $p$ is reduced to 
a Gauss hypergeometric equation. {\hfill QED}
\vskip 0.3 cm

The singular solutions do not give any solution of the PVI$\mu$ equation. 
All the other solutions do, via \eqrefp{N13.5}.
Conversely, starting from any solution $y(x)$ of PVI$\mu$, we arrive at the 
solution:
$$
\eqalign{
q&=(u_3-u_1)y\left({u_2-u_1\over u_3-u_1}\right)+u_1\cr
p&={P'(u_2)\over 2 P(q)} y'\left({u_2-u_1\over u_3-u_1}\right)-{1\over 2} 
{1\over q-u_2} \cr}
$$
of the reduced Schlesinger equations \eqrefp{N13}. To obtain a solution of 
the full Schlesinger equations, the function $k$ must be given by a quadrature:
$$
{\partial k\over\partial u_i}=(2\mu-1){q-u_i\over P'(u_i)}.
$$

\vskip 0.2 cm\noindent
We conclude this section summarizing all the above results in the following:

\proclaim Theorem 1.3. The branches of solutions of the PVI$\mu$ equation 
near a given point $x_0\in\overline\complessi\backslash\{0,1,\infty\}$, 
are in one-to-one correspondence with the triples of the monodromy matrices 
$M_1$, $M_2$, $M_3$ satisfying \eqrefp{N5} and \eqrefp{N6}, with $M_\infty$
of the form \eqrefp{N9}, none of them 
being equal to $\ID$, considered modulo diagonal conjugations.

\vskip 0.3 cm
\noindent{\bf Remark 1.4.}\quad A triple of $2\times2$ matrices 
$M_1,M_2,M_3\in SL(2;\complessi)$, considered modulo conjugations, is a point 
$\rho$ of the space of representations 
$$
\rho:F_3\to SL(2;\complessi)
$$
of the free group $F_3$ with three generators $\gamma_1,\gamma_2,\gamma_3$,
specified by
$$
M_i=\rho(\gamma_i),\quad i=1,2,3.
$$
In the general case, i.e. with the matrices ${\cal A}_i$ and 
${\cal A}_\infty$ not necassarly of the form \eqrefp{N1.3}, the corresponding 
solution $(p,q)$ of the reduced Schlesinger equations will be denoted
$$
p=p(u_1,u_2,u_3;\rho),\qquad q=q(u_1,u_2,u_3;\rho).
$$
It is locally uniquelly specified by the representation $\rho$, provided the 
non-resonancy condition of the eigenvalues of ${\cal A}_i$ and 
${\cal A}_\infty$.

\vskip 0.3 cm
\noindent{\bf  1.2. The structure of the analytic continuation.}
\vskip 0.3 cm
We parameterized branches of the solutions of PVI$\mu$ by triples of 
monodromy matrices. Now we show how do these parameters change with a 
change of the branch in the process of analytic continuation of the 
solutions along a path in $\overline\complessi\backslash\{0,1,\infty\}$. 
Recall that, as it follows from Theorem 1.1, the solutions of PVI$\mu$, 
defined in a neighborhood of a given point 
$x_0\in\overline\complessi\backslash\{0,1,\infty\}$, can be 
analytically continued to a meromorphic function on the universal 
covering of $\overline\complessi\backslash\{0,1,\infty\}$ (the above 
mentioned Painlev\'e Property). The fundamental group 
$\pi_1\left(\overline\complessi\backslash\{0,1,\infty\}\right)$ is 
non-abelian. As a consequence, the global structure of the analytic 
continuation of the solutions of PVI is more involved than that of the 
other Painlev\'e equations. In fact the solutions of PI,..., PV have 
at most two critical singularities and the corresponding fundamental group 
is abelian.

As a first step we introduce a parameterization of the monodromy matrices.

\vskip 0.3 cm
\noindent{\bf 1.2.1. The parameterization of the monodromy data}\quad Let 
${\cal M}_1$, ${\cal M}_2$ and ${\cal M}_3$ be three linear operators 
${\cal M}_i:\complessi^2\rightarrow\complessi^2$ satisfying \eqrefp{N5}. 
We introduce for them a parameterization which will be useful for 
studying the analytic continuation of the solutions of the PVI$\mu$ equation.

\proclaim Lemma 1.4. If ${\cal M}_1$, ${\cal M}_2$ are such that
$$
{\rm Tr}({\cal M}_1{\cal M}_2)\neq 2,
$$
then there exists a basis in $\complessi^2$ such that, in this basis, the 
matrices of ${\cal M}_1$, ${\cal M}_2$ have the form:
$$
M_1=\pmatrix{1&-x_1\cr 0&1\cr},\qquad M_2=\pmatrix{1&0\cr x_1&1\cr},
\autoeqno{N14}
$$
where $x_1=\sqrt{2-{\rm Tr}({\cal M}_1{\cal M}_2)}$; when ${\cal M}_1$, 
${\cal M}_2$ are such that ${\rm Tr}({\cal M}_1{\cal M}_2)= 2$, they have 
a common eigenvector, and then there exists a basis in $\complessi^2$ such 
that, in this basis, the matrices $M_1$, $M_2$ are both upper-triangular.

\noindent Proof. Due to the \eqrefp{N5}, there exist two 
vectors $e_1$ and $e_2$ such that
$$
{\cal M}_1 e_1= e_1,\qquad {\cal M}_2 e_2= e_2.
$$
We now prove that these two vectors are linearly dependent if and only if 
${\rm Tr}({\cal M}_1{\cal M}_2)= 2$. In fact if the two vectors are linearly 
dependent, then we can find a linear independent vector $e_2'$ such that, 
in the basis $(e_1,e_2')$ the matrices of $M_1$, $M_2$ have the form:
$$
M_1=\pmatrix{1&\lambda_1\cr 0&1\cr},\qquad M_2=
\pmatrix{1&\lambda_2\cr 0&1\cr},
$$
so ${\rm Tr}(M_1 M_2)= 2$. Conversely, in the basis $(e_1,e'_2)$ the 
matrix $M_1$ has the form $M_1=\pmatrix{1&\lambda_1\cr 0&1\cr}$ and, 
requiring that 
$$
{\rm Tr}({\cal M}_1{\cal M}_2)= 2,\quad {\rm eigenv}(M_2)=1,
$$ 
also the matrix $M_2$ must have the above form 
$M_2=\pmatrix{1&\lambda_2\cr 0&1\cr}$. Then, the two vectors $e_1$ and 
$e_2$ are linearly dependent.
As a consequence, if ${\rm Tr}({\cal M}_1{\cal M}_2)\neq2$, the two vectors 
$e_1$ and $e_2$ are linearly independent, and in the basis $(e_1,e_2)$ the 
matrices of $M_1$, $M_2$ have the form:
$$
M_1=\pmatrix{1&\lambda_1\cr 0&1\cr},\qquad M_2
=\pmatrix{1&0\cr\lambda_2&1\cr},
$$
with ${\rm Tr}(M_1M_2)= 2+\lambda_1\lambda_2$. Rescaling the basic 
vectors $(e_1,e_2)$, we obtain the \eqrefp{N14}. 
{\hfill QED}

\proclaim Lemma 1.5. Let ${\cal M}_1$, ${\cal M}_2$, ${\cal M}_3$ satisfy 
also the condition \eqrefp{N6} with ${\cal M}_\infty$ given by \eqrefp{N9}, 
and $2\mu\not\in\interi$. Then the following statements are true:
\item{i)} If two of the following numbers
$$
{\rm Tr}({\cal M}_1{\cal M}_2),\quad {\rm Tr}({\cal M}_1{\cal M}_3),
\quad {\rm Tr}({\cal M}_3{\cal M}_2)
$$
are equal to $2$, then one of the matrices of $M_i$ is equal to one.
\item{ii)} If ${\rm Tr}({\cal M}_1{\cal M}_2)\neq 2$, then  there exists a 
basis in $\complessi^2$ such that, in this basis, the matrices $M_1$, 
$M_2$ and $M_3$ have the form
$$
M_1=\pmatrix{1&-x_1\cr 0&1\cr},\quad M_2=\pmatrix{1&0\cr x_1&1\cr},
\quad M_3=\pmatrix{1+ {x_2 x_3\over x_1}& -{x_2^2\over x_1}\cr 
{x_3^2\over x_1}&1-{x_2 x_3\over x_1}\cr},
\autoeqno{N15}
$$
where
$$
{\rm Tr}(M_1M_2)=2-x_1^2,\quad {\rm Tr}(M_3M_2)=2-x_2^2,\quad 
{\rm Tr}(M_1M_3)=2-x_3^2,
$$
and
$$
x_1^2+x_2^2+x_3^2-x_1 x_2 x_3=4\sin^2\pi\mu.\autoeqno{N16}
$$
\item{iii)} If two triples of matrices $M_1$, $M_2$, $M_3$ and $M'_1$, 
$M'_2$, $M'_3$ satisfying \eqrefp{N6}, with none of them equal to $\ID$, 
have the form \eqrefp{N15} with the parameters $(x_1,x_2,x_3)$ and 
$(x'_1,x'_2,x'_3)$ respectively, then these triples are conjugated 
$$
M_i=T^{-1} M_i' T
$$
with some invertible matrix $T$ if and only if the triple 
$(x'_1,x'_2,x'_3)$ is equal to the triple $(x_1,x_2,x_3)$, up to the 
change of the sign of two of the coordinates.

\noindent Proof.
\item{i)} Let us assume that
$$
{\rm Tr}({\cal M}_1{\cal M}_2)=2,\quad {\rm Tr}({\cal M}_1{\cal M}_3)=2.
$$
Let $e_1$ and $e_3$ be the common eigenvectors of ${\cal M}_1$, 
${\cal M}_2$ and ${\cal M}_1$, ${\cal M}_3$ respectively, (see Lemma 1.4).
If ${\cal M}_1\neq \ID$, then the eigenvectors $e_1$ and $e_3$ coincide.
Then we can find a linear independent vector $e_2'$ such 
that, in the basis $(e_1,e_2')$ the matrices of $M_1$, $M_2$, $M_3$ all 
have the form
$$
M_i=\pmatrix{1&\lambda_i\cr 0&1\cr},\quad i=1,2,3.
$$
Then
$$
{\rm Tr}(M_3 M_2 M_1)1{\rm Tr}(M_\infty)=2.
$$
This contradicts the assumption $2\mu\not\in\interi$.
\item{ii)} Let us choose the basis such that, according to Lemma 1.4, 
the matrices $M_1$, $M_2$ have the form \eqrefp{N14}. Solving the 
equations
$$
{\rm Tr}(M_3M_2)=2-x_2^2,\quad 
{\rm Tr}(M_1M_3)=2-x_3^2,
$$
we arrive at the formula \eqrefp{N15}. The \eqrefp{N16} is obtained by 
straightforward computations from
$$
{\rm Tr}(M_3 M_2 M_1)=2\cos2\pi\mu.
$$
\item{iii)} The two triples of matrices $M_1$, $M_2$, $M_3$ and $M'_1$, 
$M'_2$, $M'_3$ are conjugated 
$$
M_i=T^{-1} M_i' T
$$
with some invertible matrix $T$ if and only if they are the matrices of 
the same operators ${\cal M}_1$, ${\cal M}_2$, ${\cal M}_3$, written in 
different bases. Since the traces do not depend on the choice of the basis, 
then
$$
x_i^2={x_i'}^2,\quad i=1,2,3.
$$
According to the proof of Lemma 1.4, the basis $(e_1,e_2)$ is uniquely 
determined up to changes of sign. A change of sign $e_1\mapsto - e_1$ 
corresponds to the change of sign $x_1\mapsto - x_1$; then the form of the 
matrix $M_3$ is preserved if and only if we change one of the signs of 
$x_2$ or $x_3$.
{\hfill QED}
\vskip 0.3 cm

\noindent{\bf Remark 1.5.}\quad The matrices \eqrefp{N15} have a simple 
geometrical meaning. Let us consider the three-dimensional linear space 
with a basis $(e_1,e_2,e_3)$ and with a skew-symmetric bilinear form 
$\{\cdot,\cdot\}$ such that 
$$
\{e_1,e_2\}=x_1,\quad \{e_1,e_3\}=x_3,\quad \{e_2,e_3\}=x_2.
$$
Let us consider the reflections $R_1,R_2,R_3$ in this space, with respect 
to the hyperplanes skew-orthogonal to the basic vectors:
$$
R_i(x)=x-\{e_i,x\}e_i,\qquad i=1,2,3.
$$
The reflections have a one-dimensional invariant subspace, namely the 
kernel of the bilinear form. The matrices of the reflections acting on 
the quotient are the \eqrefp{N15}.

\vskip 0.3 cm
\noindent{\bf Definition.}\quad
A triple $(x_1,x_2,x_3)$ is called {\it admissible}\/ if it has at most 
one coordinate equal to zero. Two such triples are called 
{\it equivalent}\/ if they are equal up to the change of two signs 
of the coordinates.
\vskip 0.3 cm

Observe that for an admissible triple $(x_1,x_2,x_3)$ none of the 
matrices \eqrefp{N15} is equal to the identity. So the admissible 
triples correspond to the non-singular solutions of the reduced 
Schlesinger equations \eqrefp{N13}. Moreover, two equivalent triples 
generate the same solution. We can summarize the above results in the 
following:

\proclaim Theorem 1.4. The branches of solutions of the PVI$\mu$ 
equation near a given point 
$x_0\in\overline\complessi\backslash\{0,1,\infty\}$ are in one-to-one 
correspondence with the equivalence classes of the admissible triples 
satisfying \eqrefp{N16}.

\noindent Proof. Starting from a solution of PVI$\mu$ 
we obtain the monodromy matrices satisfying \eqrefp{N5}. None of them is 
equal to the identity. So the canonical form \eqrefp{N15} of $M_1,M_2,M_3$ 
is determined uniquely up to a choice of the admissible triple 
$(x_1,x_2,x_3)$ within the equivalence class. Conversely, given an 
admissible triple $(x_1,x_2,x_3)$ satisfying \eqrefp{N16}, we obtain 
the matrices  $M_1,M_2,M_3$ of the form \eqrefp{N15}. 
The matrix $M_3 M_2 M_1$ is diagonalizable with the eigenvalues 
$\exp(\pm 2 \pi i \mu)$ (here we use the non-resonance condition 
$2\mu\not\in\interi$). Reducing this matrix to the diagonal form
$$
M_3M_2M_1=
T^{-1}\pmatrix{\exp( 2 \pi i \mu)&0\cr 0& \exp(-2 \pi i \mu)\cr} T
$$
we obtain the monodromy matrices $T M_i T^{-1}$ satisfying \eqrefp{N5} 
and thus specifying a branch of the solution of PVI$\mu$.

\vskip 0.3 cm
\noindent{\bf 1.2.2. Monodromy data and symmetries of PVI$\mu$.}\quad
The Painlev\'e VI equation possesses a rich family of symmetries, i. e. 
transformations of the dependent and independent variables $(y,x)$, and 
also of the parameters, that preserve the shape of the equation. The 
theory of these symmetries, and its applications to the construction of 
particular solutions, was developed in [Ok].
Here we list the symmetries which preserve our PVI$\mu$ and compute 
their action on the monodromy data.

First of all we observe that the trivial symmetry $\mu\mapsto 1-\mu$ 
preserves the Painleve' equation, i.e. $\hbox{PVI}\mu=\hbox{PVI}(1-\mu)$, 
so it maps the solutions $y(x)$ in themselves.

Then we consider the permutations of the poles $u_1,u_2,u_3$ which generate 
the action of the symmetric group $S_3$ on the solutions $y(x)$. In 
particular the involution 
$$
i_1: u_2\leftrightarrow u_3,
$$
produces the transformation
$$
x\mapsto{1\over x},\qquad y\mapsto {y\over x},\autoeqno{s1}
$$
and 
$$
i_2: u_1\leftrightarrow u_3,
$$
produces the transformation
$$
x\mapsto 1 - x,\qquad y\mapsto 1 -  y.\autoeqno{s2}
$$
Both these transformations clearly preserve the equation PVI$\mu$.

Let us compute the action of these symmetries on the monodromy data. The 
only thing that changes is the basis in the fundamental group 
$\pi_1(\overline\complessi\backslash\{u_1,u_2,u_3,\infty\})$. In fact, the 
cuts $\pi_1,\pi_2,\pi_3$ along which we take our basis 
$\gamma_1,\gamma_2,\gamma_3$, are ordered according to the order of the 
poles. Applying the transformation $i_1$ we then arrive at the new basis 
$\gamma'_1,\gamma'_2,\gamma'_3$ shown in figure 2.
\midinsert
\centerline{\psfig{file=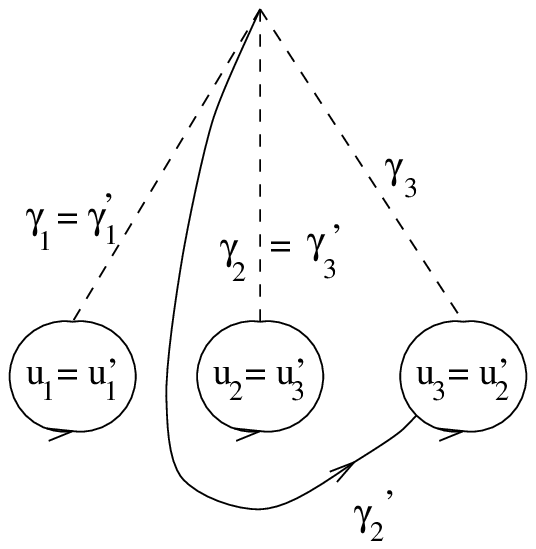,height=4cm}}
\vskip 0.6 cm
\centerline{{\bf Fig.2.} The new basis $\gamma'_1,\gamma'_2,\gamma'_3$ 
obtained by the action of $i_1$.}
\endinsert
%
%
%
\noindent This new basis has the following form
$$
\gamma'_1=\gamma_1,
\quad \gamma'_2=\gamma_2\gamma_3\gamma_2^{-1},
\quad \gamma'_3=\gamma_2.
$$
As a consequence the new monodromy matrices are
$$
 M'_1= M_1,
\quad  M'_2= M_2^{-1} M_3 M_2,
\quad  M'_3= M_2.
$$
For the second transformation $i_2$, the basis of the new loops is shown 
in figure 3. 
\midinsert
\centerline{\psfig{file=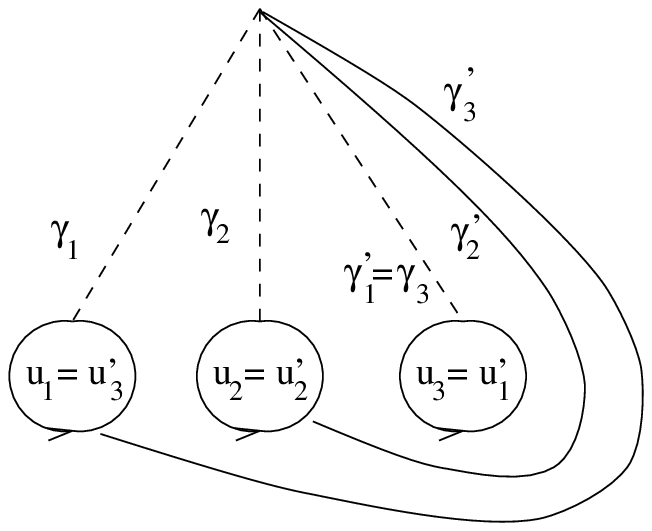,height=4cm}}
\vskip 0.7 cm
\centerline{{\bf Fig.3.} The new basis $\gamma'_1,\gamma'_2,\gamma'_3$ 
obtained by the action of $i_2$.}
\endinsert
%
%
%
\noindent It has the following form
$$
\gamma'_1=\gamma_3,
\quad \gamma'_2=\gamma_3^{-1}\gamma_2\gamma_3,
\quad \gamma'_3=\gamma_3^{-1}\gamma_2^{-1}\gamma_1\gamma_2\gamma_3.
$$
The new monodromy matrices are
$$
 M'_1= M_3,
\quad  M'_2= M_3 M_2^{-1} M_3,
\quad  M'_3= M_3 M_2 M_1 M_2^{-1} M_3^{-1}.
$$

\proclaim Lemma 1.6. In the coordinates $(x_1,x_2,x_3)$ on the space of 
the monodromy matrices, the action of the symmetries $i_1,\,i_2$ is given 
by the formulae
$$
i_1:(x_1,x_2,x_3)\mapsto (x_3-x_1 x_2,-x_2,x_1),\qquad
i_2:(x_1,x_2,x_3)\mapsto (-x_2,-x_1,x_1 x_2-x_3).
$$

\noindent The proof is straightforward.

\vskip 0.2 cm

The last symmetry is more complicated because it changes the value of the 
parameter $\mu$, i.e. $\mu\mapsto-\mu$, or equivalently $\mu\mapsto 1+\mu$  
as it follows form the fact that $\hbox{PVI}(-\mu)=\hbox{PVI}(1+\mu)$.
This simmetry comes from the following simultaneous conjugation of the 
coefficients of the Fuchsian system:
$$
A_i\rightarrow \Sigma A_i \Sigma,
$$
where
$$
\Sigma=\Sigma^{-1}=\pmatrix{0&1\cr 1&0\cr}.
$$
Indeed, 
$$ 
\Sigma A_\infty \Sigma=-A_\infty.
$$
Using the parameterization \eqrefp{N12.5} of the matrices $A_1,A_2,A_3$ by 
the coordinates $(p,q)$, we arrive at the following

\proclaim Lemma 1.7. The formula
$$
\tilde y = y {\left(p_0 (y')^2+p_1 y'+p_2\right)^2
\over q_0(y')^4+ q_1(y')^3+
 q_2(y')^2+ q_3 y'+ q_4},\autoeqno{B10}
$$
where
$$
\eqalign{
p_0&=x^2(x-1)^2,\cr
p_1& = 2 x(x-1)(y-1)[2\mu(y-x)-y]\cr
p_2& =y(y-1)[y(y-1)-4\mu(y-1)(y-x)+4\mu^2(y-x)(y-x-1)]\cr
q_0& =x^4(x-1)^4\cr
q_1& =-4 x^3(x-1)^3 y(y-1)\cr
q_2& =2 x^2(x-1)^2 y(y-1)[3y(y-1)+4\mu^2(y-x)(1+x-3y)]\cr
q_3& =4x(x-1)y^2(y-1)^2[-y(y-1)-16\mu^3(y-x)^2+4\mu^2(y-x)(3y-x-1)]\cr
q_4& =y^2(y-1)^2\big\{y^2(y-1)^2+64\mu^3y(y-1)(y-x)^2
-8\mu^2y(y-1)(y-x)(3y-x-1)+\cr
&+16\mu^4(y-x)^2[(x-1)^2+y(2+2x-3y)]\big\}\cr},\autoeqno{B11}
$$
transforms solutions of PVI$\mu$ to solutions of PVI$(-\mu)$. The class of 
equivalence of the monodromy data $(x_1,x_2,x_3)$ does not change under such 
a symmetry.

\noindent Proof. The new monodromy matrices 
$M'_1,M'_2,M'_3$ have the form
$$
M'_i=\Sigma M_i \Sigma,\quad i=1,2,3.
$$
Then, the canonical form \eqrefp{N15} of the monodromy operators does not 
change. {\hfill QED}
\vskip 0.3 cm

Other symmetries are superpositions of \eqrefp{B10} with the trivial one 
$\mu\to1-\mu$. Using these symmetries, one can transform PVI$_\mu$ to 
PVI$_{\mu'}$, with $\mu'=\pm\mu+n$ for an arbitrary integer $n$.

\vskip 0.3 cm
\noindent{\bf Remark 1.6.}\quad One can show that the above symmetries, and 
their superpositions, exhaust all the birational transformations preserving 
our one-parameter family of PVI equations. We will not do it here (see [Ok]). 
It is important, however, that these symmetries preserve the class of 
algebraic solutions of PVI$\mu$. We will classify all the algebraic 
solutions modulo the above symmetries. 

\vskip 0.3 cm
\noindent{\bf Remark 1.7.}\quad It is not difficult to show that the 
denominator of the formula \eqrefp{B10} does not vanish identically 
for any solution of PVI$\mu$, with $2\mu\not\in\interi$. Indeed, 
eliminating $y_{xx}$ and $y_x$ form the system 
$$
\eqalign{
y_{xx}=&{1\over2}\left({1\over y}+{1\over y-1}+{1\over y-x}\right) 
y_x^2 -\left({1\over x}+{1\over x-1}+{1\over y-x}\right)y_x\cr
&+{1\over2}{y(y-1)(y-x)\over x^2(x-1)^2}\left[(2\mu-1)^2+
{x(x-1)\over(y-x)^2}\right],\cr}
$$
$$
\eqalign{
Q(y_x,y,x,\mu)=&0,\cr
{{\rm d}\over {\rm d} x} Q(y_x,y,x,\mu)=&0,\cr}
$$
where $Q$ is the denominator, the resultant equation 
$$
(2\mu+1)^4\mu^{16}\left[x(x-1)^2 \right]^4\left[y(y-1)(y-x) \right]^4
$$
never vanishes.

\vskip 0.3 cm
\noindent{\bf 1.2.3. The analytic continuation of the solutions of PVI$\mu$ 
and the braid group $B_3$.}\quad In this subsection, we describe the 
procedure of the analytic continuation towards an action of the braid 
group on the admissible triples $(x_1,x_2,x_3)$ parameterizing the 
branches of the solutions of PVI$\mu$.

According to Theorem 1.1, any solution of the Schlesinger equations 
can be continued analytically from a point $(u_1^0,u_2^0,u_3^0)$ to 
another point $(u_1^1,u_2^1,u_3^1)$ along a path 
$$
(u_1(t),u_2(t),u_3(t))\in\complessi^3\backslash\{diags\},
\qquad 0\leq t\leq 1,
$$
with 
$$
u_i(0)=u_i^0,\quad\hbox{and}\quad u_i(1)=u_i^1,
$$
provided that the end-points are not the poles of the solution. The result 
of the analytic continuation depends only on the homotopy class of the 
path in $\complessi^3\backslash\{diags\}$. Particularly, to find all the 
branches of a solution near a given point $u^0=(u_1^0,u_2^0,u_3^0)$ one 
has to compute the results of the analytic continuation along any homotopy 
class of closed loops in $\complessi^3\backslash\{diags\}$ with the 
beginning and the end at the point $u^0=(u_1^0,u_2^0,u_3^0)$.
Let 
$$
\beta\in\pi_1\left(\complessi^3\backslash\{diags\};\, u^0\right)
$$
be an arbitrary loop. Any solution of the Schlesinger equations near the point 
$u^0=(u_1^0,u_2^0,u_3^0)$, is uniquely determined by the monodromy matrices 
$M_1$, $M_2$ and $M_3$, computed in the basis $\gamma_1$, $\gamma_2$, 
$\gamma_3$. 
Continuing analytically this solution along the loop $\beta$, we arrive at 
another branch of the same solution near $u^0$. This new branch is 
specified, according to Theorem 1.3, by some new monodromy matrices 
$M_1^\beta$, $M_2^\beta$ and $M_3^\beta$, computed in the same basis 
$\gamma_1$, $\gamma_2$, $\gamma_3$. Our nearest goal is to compute these 
new matrices for any loop 
$\beta\in\pi_1\left(\complessi^3\backslash\{diags\};u^0\right)$.

The fundamental group 
$\pi_1\left(\complessi^3\backslash\{diags\};\, u^0\right)$ is isomorphic 
to the pure (or unpermuted) braid group, $P_3$ with three strings 
(see [Bir]); this is a subgroup of the full braid group $B_3$. The full 
braid group is isomorphic to the fundamental group of the same space 
where the permutations are allowed:
$$
B_3\simeq\pi_1\left(\complessi^3\backslash\{diags\}\big/S_3;\, u^0\right),
$$
$S_3$ being the symmetric group acting by permutations of the coordinates 
$(u_1,u_2,u_3)$. Any loop in $B_3$ has the form
$$
(u_1(t),u_2(t),u_3(t))\in\complessi^3\backslash\{diags\},
\qquad 0\leq t\leq 1,
$$
with 
$$
u_i(0)=u_i^0,\quad u_i(1)=u^0_{p(i)},
$$
where $p$ is a permutation of $\{1,2,3\}$.
The elements of the subgroup $P_3$ of pure braids are specified by 
the condition $p={\rm id}$.

To simplify the computations we extend the procedure of the analytic 
continuation to the full braid group
$$
M_1,\,M_2,\,M_3\mapsto M_1^\beta,\,M_2^\beta,\,M_3^\beta, 
\quad \beta\in B_3=
\pi_1\left(\complessi^3\backslash\{diags\}\big/S_3;\, u^0\right).
$$
For a generic braid $\beta\in B_3$, the new monodromy matrices describe the 
superposition of the analytic continuation and of the permutation 
$$
u_i\mapsto u_{p(i)},\quad {\cal A}_i\mapsto {\cal A}_{p(i)}.\autoeqno{N20}
$$
The braid group $B_3$ admits a presentation with generators $\beta_1$ and 
$\beta_2$ and the defining relation
$$
\beta_1 \beta_2 \beta_1=\beta_2 \beta_1\beta_2.
$$
The generators $\beta_1$ and $\beta_2$ are shown in the figure 4.

\midinsert
\centerline{\psfig{file=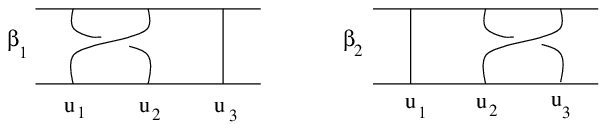,height=2cm}}
\vskip 0.5 cm
\centerline{{\bf Fig.4.} The generators of the braid group $B_3$.}
\endinsert

\proclaim Lemma 1.8. For the generators $\beta_1$, $\beta_2$ shown in the 
figure 4, the matrices $M_i^\beta$ have the following form:
$$
M_1^{\beta_1}=M_2,\quad M_2^{\beta_1}=M_2 M_1 M_2^{-1},
\quad M_3^{\beta_1}= M_3,\autoeqno{D2}
$$
$$
M_1^{\beta_2}=M_1,\quad M_2^{\beta_2}=M_3,
\quad M_3^{\beta_2}= M_3 M_2 M_3^{-1}.\autoeqno{D3}
$$

\noindent Proof. 
Changing the positions of the points $u_1$ and $u_2$ by the braid $\beta_1$, 
the basis of the loops will be deformed into the new basis 
$\gamma_1',\gamma_2',\gamma_3'$ shown in the figure 5.

\midinsert
\centerline{\psfig{file=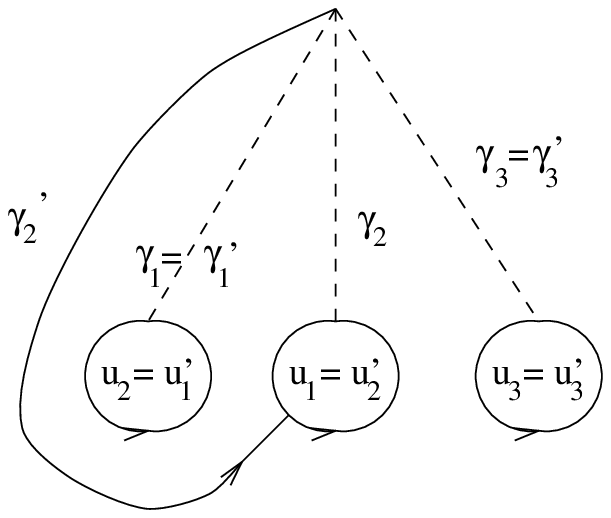,height=4cm}}
\vskip 0.6 cm
\centerline{{\bf Fig.5.} The new loops $\gamma'_i$ obtained under the 
action of the braid $\beta_1$.}
\endinsert

Thanks to the fact that we deal with isomonodromy deformations, the 
monodromy matrices $M_i'$ of the system \eqrefp{N1} with respect to the 
new basis $\gamma_1',\gamma_2',\gamma_3'$ are the same $M_i$, up to the 
reordering:
$$
M_1'=M_2,\quad M_2'=M_1,\quad M_3'=M_3.\autoeqno{D5}
$$
We want to compute the monodromy matrices with respect to the old basis 
$\gamma_1,\gamma_2,\gamma_3$. To this aim we notice the following obvious 
relation in the fundamental group:
$$
\gamma_1=\gamma'_1,\quad \gamma_2=(\gamma_1')^{-1}\gamma_2'\gamma_1',
\quad \gamma_3=\gamma_3'.
$$
Using this relations and the \eqrefp{D5}, we immediately obtain the 
\eqrefp{D2}. Similarly the deformation of the basis of the fundamental 
group corresponding to the braid $\beta_2$ is shown in the figure 6.
\midinsert
\centerline{\psfig{file=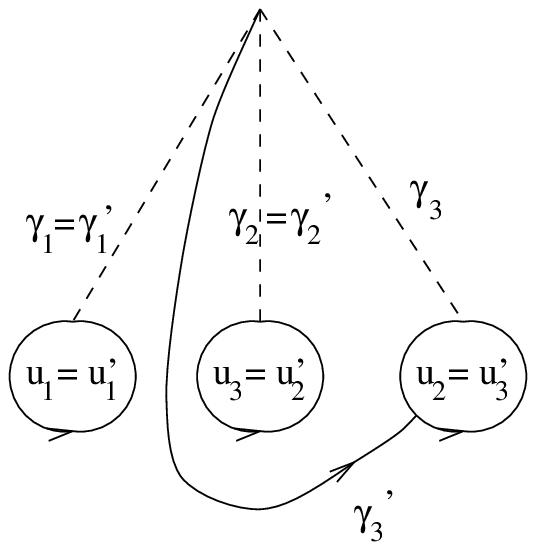,height=4cm}}
\vskip 0.5 cm
\centerline{{\bf Fig.6.} The new loops $\gamma'_i$ obtained under the 
action of the braid $\beta_2$.}
\endinsert

Here we have the permutation
$$
M_1'=M_1\quad 
M_2'=M_3',\quad
M_3'=M_2,
$$
and the relations in the fundamental group:
$$
\gamma_1=\gamma'_1,\quad \gamma_2=\gamma_2',
\quad \gamma_3=(\gamma_2')^{-1}\gamma_3'\gamma_2'.
$$
From this we obtain the \eqrefp{D3} and the lemma is proved.{\hfill QED}
\vskip 0.3 cm

The action \eqrefp{D2}, \eqrefp{D3} of the braid group on the triples of 
monodromy matrices commutes with the diagonal conjugation of them. As a 
consequence this action not only describes the structure of the analytic 
continuation of the solutions of the Schlesinger equations \eqrefp{N10}, 
but also of the reduced ones \eqrefp{N13}. Moreover, the class of the 
singular solutions is closed under this analytic continuation. In fact 
if some of the matrices $M_i$ is equal to $\ID$ then for any $\beta$ there 
is a $j$ such that $M_j^\beta=\ID$. As a consequence the following lemma 
holds true:

\proclaim Lemma 1.9. The structure of the analytic continuation of the 
solutions of the PVI$\mu$ equation is determined by the action \eqrefp{D2}, 
\eqrefp{D3} of the braid group on the triples of monodromy matrices.

Our next step is to rewrite the action \eqrefp{D2}, \eqrefp{D3} of the 
braid group in the coordinates $(x_1,x_2,x_3)$ in the space of the monodromy 
data. This is given by the following

\proclaim Lemma 1.10. In the coordinates $(x_1,x_2,x_3)$, the action 
\eqrefp{D2}, \eqrefp{D3} of the braid group is given by the formulae:
$$
\eqalign{
\beta_1:(x_1,x_2,x_3)&\mapsto(-x_1,x_3-x_1x_2,x_2),\cr
\beta_2:(x_1,x_2,x_3)&\mapsto(x_3,-x_2,x_1-x_2 x_3).\cr}\autoeqno{BD}
$$

\noindent Proof. The above formulae are obtained by 
straightforward computations from \eqrefp{D2}, \eqrefp{D3} by means 
of the parameterization of the monodromy matrices \eqrefp{N15}.

\vskip 0.3 cm
We can summarize the results of this section in the following:

\proclaim Theorem 1.5. The structure of the analytic continuation of the 
solutions of the PVI$\mu$ equation is determined by the action \eqrefp{BD} 
of the braid group on the triples $(x_1,x_2,x_3)$.

\vskip 0.3 cm
\noindent{\bf Remark 1.8.}\quad It is easy to see that the braid 
$(\beta_1\beta_2)^3$ acts trivially on the monodromy data. This braid is the 
generator of the center of $B_3$ (see [Bir]). The quotient 
$$
B_3/center\simeq PSL(2;\interi)
$$
coincides with the mapping class group of the complex plane with three 
punctures [Bir].
Also in the general case, the structure of analytic continuation of solutions
of PVI equation is described by the following natural action 
$\rho\to\rho^\beta$ of the mapping class group on the representation space 
(see remark 1.4)
$$
\rho^\beta(\gamma)=\rho(\beta^{-1}_\star(\gamma))\autoeqno{ro}
$$
where
$$
\gamma\in F_3\simeq
\pi_1\left(\overline\complessi\backslash\{u_1,u_2,u_3,\infty\},\infty\right),
$$
$$
\beta:\overline\complessi\backslash\{u_1,u_2,u_3,\infty\}
\to \overline\complessi\backslash\{u_1,u_2,u_3,\infty\},\qquad 
\beta(\infty)=\infty
$$
is a homeomorphism, and
$$
\rho:F_3\to SL(2;\complessi).
$$
Our action \eqrefp{BD} is obtained restricting \eqrefp{ro} onto the subspace
of representations of the form \eqrefp{N5}. The problem of selection of 
algebraic solutions of Painlev\'e VI (see below) with generic values of the 
parameters $\alpha,\beta,\gamma,\delta$ can be reduced to the classification
of finite orbits of the action \eqrefp{ro}.

\vskip 0.3 cm
\noindent{\bf  1.3. Monodromy data and algebraic solutions of the PVI$\mu$ 
equation.}

\vskip 0.3 cm
\noindent{\bf 1.3.1. A preliminary discussion on the algebraic solutions of 
the PVI$\mu$ equation and their monodromy data.}\quad Here we state some 
necessary condition for the triples $(x_1,x_2,x_3)$ to generate the 
algebraic solutions.

\vskip 0.3 cm
\noindent {\bf Definition.} A solution $y(x)$ is called {\it algebraic}\/ 
if there exists a polynomial in two variables such that 
$$
F(x,y(x))\equiv0.
$$
\vskip 0.3 cm

If $y(x)$ is an algebraic solution then the correspondent solution 
$p({\bf u})$, $q({\bf u})$, ${\bf u}=(u_1,u_2,u_3)$ of the reduced 
Schlesinger equations 
\eqrefp{N13} is also algebraic. According to Theorem 1.1, the solutions
of the reduced Schlesinger equations \eqrefp{N13} can ramify only on 
the diagonals $u_1=u_2$, $u_1=u_3$, $u_3=u_2$. Analogously the 
ramification points of $y(x)$ are allowed to lie only at $0,\,1,\,\infty$. 

We now characterize the monodromy data such that the correspondent solution 
of the PVI$\mu$ equation is algebraic.

\proclaim Lemma 1.11. A necessary and sufficient condition for a solution of 
PVI$\mu$ to be algebraic is that the correspondent monodromy matrices, 
defined modulo diagonal conjugations, have a finite orbit under the action 
of the braid group \eqrefp{D2}, \eqrefp{D3}.

\noindent Proof. By definition, any algebraic function 
has a finite number of branches. Allowing also the permutations 
\eqrefp{N20}, we still obtain a finite number of values for $M_1^\beta$, 
$M_2^\beta$ and $M_3^\beta$,
$\beta\in B_3$ up to diagonal conjugations.{\hfill QED}
\vskip 0.3 cm

\proclaim Corollary 1.2. An admissible triple $(x_1,x_2,x_3)$ specifies an 
algebraic solution of PVI$\mu$, with $2\mu\not\in\interi$, if and only if 
it satisfies \eqrefp{N16} and its orbit, under the action \eqrefp{BD} of 
the braid group, is finite.

\noindent{\bf Remark 1.9.}\quad We stress that the action \eqrefp{BD} 
preserves the relation \eqrefp{N16}.

\vskip 0.2 cm
In this way, the problem of the classification of all the algebraic 
solutions of the PVI$\mu$ reduces to the problem of the classification of 
all the finite orbits of the action \eqrefp{BD} under the braid group in 
the three dimensional space (see [Dub], appendix F). Here we give a simple 
necessary condition for a triple $(x_1,x_2,x_3)$ to belong to a finite orbit.

\proclaim Lemma 1.12. Let $(x_1,x_2,x_3)$ be a triple belonging to a finite 
orbit. Then:
$$
x_i=-2\cos\pi r_i,\quad r_i\in{\bf Q},\quad 0\leq r_i\leq 1,
\quad i=1,2,3.
\autoeqno{N17}
$$
Here ${\bf Q}$ is the set of rational numbers.

\noindent Proof. Let us prove the statement for, say, the 
coordinate $x_1$. Consider the transformation
$$
\beta_1^2:(x_1,x_2,x_3)\mapsto(x_1,x_2+x_1 x_3-x_1^2 x_2,x_3-x_1 x_2),
$$
as a linear map on the plane $(x_2,x_3)$. This linear map preserves the 
quadratic form
$$
x_2^2+x_3^2-x_1x_2x_3.
$$
If $x_1=2$, we put $r_1=1$; otherwise we reduce the quadratic form to the 
principal axes, introducing the new coordinates
$$
\tilde x_2={\sqrt{2+x_1}\over 2}(x_2-x_3),\qquad 
\tilde x_3={\sqrt{2+x_1}\over 2}(x_2+x_3).
$$
In these new coordinates the preserved quadratic form becomes a sum of 
squares and the transformation $\beta_1^2$ is a rotation by the angle 
$\pi+2 \alpha$, where $\alpha$ is such that $x_1=-2\cos\alpha$. To have a 
finite orbit of $(\tilde x_2,\tilde x_2)$ under the iterations of 
$\beta_1^2$, the angle $\alpha$ must be a rational multiple of $\pi$. In 
this way the statement for $x_1$ is proved. 
To prove it for $x_2$ and $x_3$ we have to consider the iterations of 
$\beta^2_2$ and $\beta_2^{-1} \beta_1^2\beta_2$ respectively. {\hfill QED}
\vskip 0.3 cm

\noindent{\bf Remark 1.10.}\quad Thanks to the above lemma, for the finite 
orbits of the braid group, it is equivalent to deal with the triples 
$(x_1,x_2,x_3)$, or with the {\it triangles}\/ with angles 
$(\pi r_1,\pi r_2,\pi r_3)$, with $x_i=-2\cos\pi r_i$ and 
$0\leq r_i\leq 1$ (we may assume, changing if necessary two of the 
signs, that at most one of the $x_i$ is positive). Observe that the 
quantity 
$$
x_1^2+x_2^2+x_3^2-x_1 x_2 x_3 -4
$$ 
is greater than 0 if and only if the triangle $(r_1,r_2,r_3)$ is 
hyperbolic, namely $\sum r_i<1$; it is equal to 0, if and only if 
the triangle $(r_1,r_2,r_3)$ is flat, namely $\sum r_i=1$, and it is 
less than 0 if and only if the triangle $(r_1,r_2,r_3)$ is spherical, 
namely $\sum r_i>1$. Thanks to \eqrefp{N16}, a flat triangle gives a 
resonant value of $\mu$, and it is thus forbidden.
\vskip 0.3 cm

\noindent{\bf 1.3.2. Classification of the triples $(x_1,x_2,x_3)$ 
corresponding to the algebraic solutions.}\quad
We deal with the classification of all the finite orbits of the triples 
$(x_1,x_2,x_3)$ of the form \eqrefp{N17}, with at 
most\footnote{${}^{3}$}{This corresponds to the fact that we deal 
only with admissible triples.} one $r_i$ being equal to ${1\over 2}$. 
According to Lemma 1.12, any point of these $B_3$-orbits must have 
the same form \eqrefp{N17}. This condition is crucial in the classification.
\vskip 0.3 cm
\noindent{\bf Definition.} We say that an admissible triple $(x_1,x_2,x_3)$ 
is {\it good}\/ if for any braid $\beta\in B_3$ one has
$$
\beta(x_1,x_2,x_3)=
(-2\cos\pi r_1^\beta,-2\cos\pi r_2^\beta,-2\cos\pi r_3^\beta),
$$
with some rational numbers $0\leq r_i^\beta\leq 1$.

\proclaim Theorem 1.6. Any good triple belongs to the orbit of one of the 
following five
$$
\left(-2\cos{\pi\over2},-2\cos{\pi\over3},-2\cos{\pi\over3}\right),
\autoeqno{cl1}
$$
$$
\left(-2\cos{\pi\over2},-2\cos{\pi\over3},-2\cos{\pi\over4}\right),
\autoeqno{cl2}
$$
$$
\left(-2\cos{\pi\over2},-2\cos{\pi\over3},-2\cos{\pi\over5}\right),
\autoeqno{cl3}
$$
$$
\left(-2\cos{\pi\over2},-2\cos{\pi\over3},-2\cos{2\pi\over5}\right),
\autoeqno{cl4}
$$
$$
\left(-2\cos{\pi\over2},-2\cos{\pi\over5},-2\cos{2\pi\over5}\right).
\autoeqno{cl5}
$$
All these orbits are finite and pairwise distinct. They contain all the 
permutations of the triples \eqrefp{cl1}, \eqrefp{cl2}, \eqrefp{cl3}, 
\eqrefp{cl4} and \eqrefp{cl5}, and also the triples
$$
\left(2\cos{\pi\over3},2\cos{\pi\over3},2\cos{\pi\over3}\right),
\eqno{(\eqref{cl1}')}
$$
$$
\left(-2\cos{2\pi\over3},-2\cos{\pi\over4},-2\cos{\pi\over4}\right),
\eqno{(\eqref{cl2}')}
$$
$$
\left(-2\cos{2\pi\over3},-2\cos{\pi\over5},-2\cos{\pi\over5}\right)\qquad
\left(-2\cos{4\pi\over5},-2\cos{4\pi\over5},-2\cos{4\pi\over5}\right),
\eqno{(\eqref{cl3}')}
$$
$$
\left(-2\cos{2\pi\over3},-2\cos{2\pi\over5},-2\cos{2\pi\over5}\right),\qquad
\left(-2\cos{2\pi\over5},-2\cos{2\pi\over5},-2\cos{2\pi\over5}\right),
\eqno{(\eqref{cl4}')}
$$
$$
\eqalign{
\left(-2\cos{3\pi\over5},-2\cos{\pi\over3},-2\cos{\pi\over5}\right),\,
&\left(-2\cos{2\pi\over5},-2\cos{\pi\over3},-2\cos{\pi\over3}\right),\cr
\left(-2\cos{2\pi\over3},\right.
&\left. -2\cos{\pi\over3},-2\cos{\pi\over5}\right),\cr}
\eqno{(\eqref{cl5}')}
$$
respectively, together with all their permutations.

\proclaim Corollary 1.3. There are five finite orbits of the action 
\eqrefp{BD} of the braid group on the space of the admissible triples
$(x_1,x_2,x_3)$ satisfying
$$
x_1^2+x_2^2+x_3^2-x_1 x_2 x_3\neq4.
$$
The lengths of the orbits \eqrefp{cl1}, \eqrefp{cl2}, \eqrefp{cl3}, 
\eqrefp{cl4} and \eqrefp{cl5}, are equal to $4$, $9$, $10$, $10$ and $18$ 
respectively. 

\noindent{\bf Remark 1.11.}\quad The action of the pure braid group $P_3$ 
on the above orbits gives the same orbits for any of them but \eqrefp{cl2}.
The orbit \eqrefp{cl2}, under the action of the pure braid group $P_3$, 
splits into three different orbits of three points. So the $P_3$-orbit 
\eqrefp{cl1} has four points, the three $P_3$-orbits \eqrefp{cl2} 
have three points each, \eqrefp{cl3} and \eqrefp{cl4} have ten points 
each and  \eqrefp{cl5} has eighteen points. These orbits give rise to all 
the algebraic solutions of the PVI$\mu$ equation, for $\mu$ is given by 
\eqrefp{N16}. The number of the points of each orbit with respect to the 
action of $P_3$ coincides with the number of the branches of the 
correspondent algebraic solution.

\vskip 0.3 cm
\noindent
{\bf Proof of Theorem 1.6.}\quad 
The braid group acting on the classes of triples $(x_1,x_2,x_3)$, is 
generated by the braid $\beta_1$ and by the cyclic permutation:
$$
(x_1,x_2,x_3)\mapsto(x_3,x_1,x_2).
$$
As a consequence it suffices to study the operator:
$$
(x_i,x_j,x_k)\mapsto (-x_i,x_j, x_k-x_i x_j),
$$
up to cyclic permutations.
This transformation works on the triangles with angles 
$\pi r_i$, $\pi r_j$, $\pi r_k$ as follows:
$$
(r_i,r_j,r_k)\mapsto (1-r_i,r_j, r'_k),\autoeqno{cl5.5}
$$
where $r'_k$ is such that:
$$
\cos\pi r_k' = \cos\pi r_k+ 2\cos\pi r_i \cos\pi r_j.\autoeqno{cl6}
$$
The first step is to classify all the rational triples $(r_i,r_j,r_k)$ 
such that $r'_k$, defined by \eqrefp{cl6} is a rational number, 
$1>r'_k>0$, for every choice of $i\neq j\neq k \neq i$, $i,j,k=1,2,3$. 
Equivalently we want to classify all the rational solutions of the 
following equation:
$$
\cos\pi r_k+ \cos\pi (r_i+ r_j) + \cos\pi(r_i - r_j)+ \cos\pi(1- r'_k)=0,
$$
or all the rational quadruples $(\varphi_1,\varphi_2,\varphi_3,\varphi_4)$ 
such that:
$$
\cos 2 \pi \varphi_1+ \cos 2\pi \varphi_2+\cos 2 \pi \varphi_3
+ \cos 2\pi \varphi_4 = 0,\autoeqno{cl7}
$$
where the $\varphi_i$ are related with the $r_i$ by the following 
relations:
$$
\varphi_1=r_k/2,  \quad \varphi_2={r_i+ r_j\over 2},\quad 
\varphi_3={|r_i- r_j|\over 2},\quad \varphi_4={|1- r'_k|\over 2}.
\autoeqno{cl6.5}
$$
Such a classification is given by the following:

\proclaim Lemma 1.13. The only rational solutions 
$(\varphi_1,\varphi_2,\varphi_3,\varphi_4)$, $0\leq\varphi_i<1$, considered 
up to permutations  and up to transformations $\varphi_i\to1-\varphi_i$, of 
the equation \eqrefp{cl7} consist of the following non--trivial solutions:
$$
\left({1\over 30},{11\over30},{2\over5},{1\over6}\right)\eqno{(a)}
$$
$$
\left({7\over 30},{17\over30},{1\over5},{1\over6}\right)\eqno{(b)}
$$
$$
\left({1\over7},{2\over7},{3\over7},{1\over6}\right)\eqno{(c)}
$$
and of the following ``trivial'' ones, of three types:
\noindent\item{(d):} $\cos 2 \pi \varphi_4=0$. The solutions obtained in [Cro]
have the form
$$
(d.1)\, : \,
\left({1\over3},{1\over10},{3\over10},{1\over4}\right),\;
(d.2)\, : \,
\left(\varphi,\varphi+{1\over3},\varphi+{2\over3},{1\over4}\right),\;
(d.3)\, : \,
\left({1\over4},\varphi,|\varphi-{1\over2}|,{1\over4}\right),
$$
where $\varphi$ is any rational number $0\leq\varphi<1$.
\item{(e):} $\cos 2 \pi \varphi_4=1$. The solutions obtained in [Gor]
have the form
$$
(e.1)\, : \,
\left({1\over3},{1\over4},{1\over3},0\right),\qquad
(e.2)\, : \,
\left({1\over2},\varphi,|\varphi-{1\over2}|,0\right),\qquad
(e.3)\, : \,
\left({1\over3},{1\over5},{2\over5},0\right),
$$
where $\varphi$ is any rational number $0\leq\varphi<1$. 
\item{(f):} 
$\cos 2 \pi \varphi_1+ \cos 2\pi \varphi_2=0, \quad
\cos 2 \pi \varphi_3 +\cos 2 \pi \varphi_4=0$. The solutions are obvious
$$
\varphi_2=|1/2- \varphi_1|,\qquad \varphi_4=|1/2- \varphi_3|,
$$
where $\varphi_1$, $\varphi_3$ are two arbitrary rational numbers 
$0\leq\varphi_i<1$.

\noindent Proof.\quad
We follow the idea of Gordan [Gor] (see also [Cro]). In this proof we 
use the same notations as in [Cro], except for the $\varphi_i$ which there 
are called $r_i$. Let us recall the notations. Let $\varphi_k={n_k\over d_k}$ 
where $d_k,\, n_k$ are either positive coprime integers,  $d_k>n_k$, or 
$n_k=0$. Let $p$ be the largest prime which is a divisor of 
$d_1,\, d_2,\, d_3$, or $d_4$ and let $\delta_k, \, l_k, \, c_k, \nu_k$ be 
the integers such that
$$
d_k=\delta_k p^{l_k}\quad\hbox{and}\quad n_k=c_k \delta_k + \nu_k p^{l_k},
$$
where $\delta_k $ is prime to $p$, $0\leq c_k<p^{l_k}$, $c_k=0$ if $l_k=0$, 
but otherwise $c_k$ is prime to $p$. So
$$
\varphi_k={\nu_k\over\delta_k}+{c_k\over p^{l_k}}=f_k+{c_k\over p^{l_k}}.
$$
We assume that $l_1\geq l_2 \geq l_3 \geq l_4$ and define the function:
$$
g_k(x)=\left\{
\matrix{
&{1\over2}\left[e^{2\pi i f_k}x^{c_kp^{l_1-l_k}}+ 
e^{-2\pi i f_k}x^{p^{l_1}-c_kp^{l_1-l_k}}\right] &
\quad\hbox{if}\quad c_k\neq0\cr
&\cos 2\pi \varphi_k&\quad\hbox{if}\quad c_k=0\cr
}\right.
$$
and, in our case: 
$$
U(x)=\sum_1^4g_k(x).
$$
As in [Cro], 
$g_k\left(\exp\left({2\pi i\over p^{l_1}}\right)\right)=\cos2\pi\varphi_k$ 
and $U\left(\exp\left({2\pi i\over p^{l_1}}\right)\right)=0$.
Let us introduce the polynomial 
$$
P(x)=1+x^{p^{l_1-1}}+x^{2 p^{l_1-1}}\cdots x^{(p-1)p^{l_1-1}}.
$$
This is the minimal polynomial of $\exp\left({2\pi i\over p^{l_1}}\right)$ 
with coefficients in ${\bf Q}$, that is such that i) 
$P\left(\exp\left({2\pi i\over p^{l_1}}\right)\right)=0$ and ii) $P(x)$ is 
irreducible in the ring of polynomials with rational coefficients. A stronger 
result was proved by Kronecker (see [Kr]): the polynomial $P(x)$ remains 
irreducible over any extension of the form 
${\bf Q}(\zeta_1,\cdots,\zeta_n)$, where $\zeta_i$ is a root of the unity 
of the order coprime with $p$. As a consequence, the following lemma holds 
true (see [Gor])

\proclaim Lemma 1.14. If we express the polynomial $U(x)$ as a sum of 
polynomials $U_t(x)$, 
$$
U(x)=\sum_{t=0}^{p^{l_1-1}-1} U_t(x),
$$ 
where $U_t(x)$ contains those terms of $U(x)$ of the form $b x^c$ with 
$c=t\, {\rm mod}\left(p^{l_1-1}\right)$, then every $U_t(x)$ is divisible 
by $P(x)$.

We now apply this lemma in our case. The indices of the powers of $x$ are:
$$
c_1, \; p^{l_1} - c_1, \; c_2 p^{l_1-l_2}, \; p^{l_1}-c_2 p^{l_1-l_2},
\; c_3 p^{l_1-l_3}, \; p^{l_1}-c_3 p^{l_1-l_3},
\; c_4 p^{l_1-l_4}, \; p^{l_1}-c_4 p^{l_1-l_4}.
$$
If all the following conditions are satisfied:
$$
l_1,l_2,l_3>1,\quad l_1>l_2,l_3,l_4,\quad l_2>l_3,l_4,\quad l_3>l_4, 
\quad l_4>0,
$$
then there are no indices equal to each other 
${\rm mod}\left(p^{l_1-1}\right)$ and there is no solution of \eqrefp{cl7}. 
So we have to study the cases in which one of them is violated.

\noindent 1): $l_1=1\geq l_2\geq l_3\geq l_4$. In this case, since the 
degree of $U(x)$ is less than $p$, and the degree of $P(x)$ is $p-1$, being 
$U(x)$ divisible by $P(x)$, we must have $U(x)=m P(x)$, for some constant 
$m$. There are four possibilities:
\item{1.1)}: $l_1=l_2=l_3=l_4=1$ then $U(0)=0$ and P(0)=1. Then $m=0$ and 
$U(x)\equiv0$; moreover if the sum of two (three) terms representing two 
(three) of the functions $g_k$ vanishes, then the sum of the two (three) 
functions vanishes. As a consequence there are only the following 
possibilities:
\item{1.1.1)}: $g_i=-g_j$ and $g_k=-g_l$ for some distinct 
$i,j,k,l=1,\cdots4$. This gives rise to the trivial case (f).
\item{1.1.2)}: $g_l=0$ for some $l=1,\cdots4$; this is the trivial 
case (d).
\item{1.1.3)}: $U(x)$ contains only two powers of $x$. If $b_1,\cdots,b_4$ 
are the coefficients if one of the powers $x^c$, then:
$$
b_1+b_2+b_3+b_4=0,\quad\hbox{and}\quad 
{1\over b_1} +{1\over b_2} +{1\over b_3} +{1\over b_4} =0,
$$
namely $b_1,\cdots, b_4$ are the solutions of the following biquadratic 
equation:
$$
z^4+(b_1 b_2+ b_1 b_3+b_1 b_4+b_2 b_3+b_2 b_4+b_3 b_4)z^2+b_1 b_2 b_3 b_4=0.
$$
As a consequence $b_i+b_j=0$, $b_l+b_k=0$, ${1\over b_i} +{1\over b_j}=0$ and 
${1\over b_l} +{1\over b_k}=0$, for some distinct $i,j,k,l=1,\cdots,4$. 
Then this case reduces to the trivial case (f).
\item{1.2)} $l_1=l_2=l_3=1$, $l_4=0$; then $U(0)=\cos2\pi\varphi_4$ and then
$U(x)=\cos2\pi\varphi_4 P(x)$,
where $P(x)$ is a polynomial with $p$ powers of $x$. Since in $U$ we have at 
most $7$ powers and $p$ must be prime, then $p$ can only be equal to 
$2,3,5,7$.
\item{1.2.1)} Case $p=2$. Since $p$ is the largest prime in 
$d_1,\cdots, d_4$, we must have $d_1=d_2=d_3=d_4=2$ and $\delta_k=1$. 
Then $\nu_k=0$, $c_k=1$ and this provides no solution.
\item{1.2.2)} Case $p=3$. In this case there are the two following 
possibilities:
$$
{1\over 2} e^{2\pi i f_1}+{1\over 2} e^{2\pi i f_2}+
{1\over 2} e^{2\pi i f_3}=
\cos2\pi\varphi_4= 
{1\over 2} e^{-2\pi i f_1}+{1\over 2} e^{-2\pi i f_2}+
{1\over 2} e^{-2\pi i f_3}
$$
or
$$ 
{1\over 2} e^{-2\pi i f_1}+{1\over 2} e^{2\pi i f_2}+
{1\over 2} e^{2\pi i f_3}=
\cos2\pi\varphi_4= 
{1\over 2} e^{2\pi i f_1}+{1\over 2} e^{-2\pi i f_2}+
{1\over 2} e^{-2\pi i f_3}.
$$
In both the case one can show that there are no solutions. In fact, for 
example, in the first case one has to solve the following equations:
$$
2\cos2\pi\varphi_4=\cos2\pi f_1+\cos2\pi f_2+\cos2\pi f_3,\qquad
\sin2\pi f_1+\sin2\pi f_2+\sin2\pi f_3=0.
$$
Using the classification of all the possible rational solution (d.1), (d.2), 
(d.3) of the case (d), one can show that there are no solutions.
\item{1.2.3)} Case $p=5$. In this case we have:
$$
{1\over2} e^{2\pi i f_k} ={1\over2} e^{-2\pi i f_k}=
{1\over2} e^{2\pi i f_i} +{1\over2} e^{\pm 2\pi i f_j}= \cos2\pi\varphi_4,
$$
for some distinct $i,j,k=1, 2, 3$.
Then $f_k$ is $0$ or ${1\over2}$ and $\varphi_4={1\over6}$ or 
$\varphi_4={1\over3}$ respectively. Following the same computations of 
[Cro] we obtain the two solutions (a) and (b).
\item{1.2.4)} Case $p=7$. In this case we have:
$$
{1\over2} e^{2\pi i f_1} ={1\over2} e^{-2\pi i f_1}=
{1\over2} e^{2\pi i f_2} ={1\over2} e^{-2\pi i f_2}=
{1\over2} e^{2\pi i f_3} ={1\over2} e^{-2\pi i f_3}=  
\cos2\pi\varphi_4,
$$
which has the following solutions:
$$
f_1=f_2=f_3=0\quad\hbox{and}\quad \varphi_4={1\over6}\quad\hbox{or}\quad 
f_1=f_2=f_3={1\over2}\quad\hbox{and}\quad \varphi_4={1\over3}.
$$
This gives the solution (c).
\item{1.3)} $l_1=l_2=1$ and $l_3=l_4=0$. Then 
$U(x)=(\cos 2\pi \varphi_3+ \cos 2\pi \varphi_4)P(x)$; again in $U$ we 
have at most $5$ powers and then $p=2,3,5$. The case $p=2$ is treated as 
in [Cro];\item{1.3.1)}: In the case $p=3$ either
$$
{1\over2} e^{2\pi i f_1} +{1\over2} e^{2\pi i f_2}=
{1\over2} e^{-2\pi i f_1} + {1\over2} e^{-2\pi i f_2}=  
\cos2\pi\varphi_3+\cos2\pi\varphi_4,
$$
or:
$$
{1\over2} e^{2\pi i f_1} +{1\over2} e^{-2\pi i f_2}=
{1\over2} e^{-2\pi i f_1} + {1\over2} e^{2\pi i f_2}=  
\cos2\pi\varphi_3+\cos2\pi\varphi_4.
$$
In the former case, for $f_1=f_2$, with 
$\cos2\pi f_2=\cos2\pi\varphi_3+\cos2\pi\varphi_4$ and this gives again the 
solution (b). The latter case is equivalent.
\item{1.3.2)}: In the case $p=5$ one has:
$$
{1\over2} e^{2\pi i f_1} ={1\over2} e^{-2\pi i f_1}=
{1\over2} e^{2\pi i f_2} ={1\over2} e^{-2\pi i f_2}=
\cos2\pi\varphi_3+\cos2\pi\varphi_4,
$$
which gives $f_1=f_2=0$ or $f_1=f_2={1\over2}$. We treat the former case 
(the latter is equivalent); then 
$\cos2\pi\varphi_3+\cos2\pi\varphi_4={1\over2}$ and
we can show that this case reduces to the trivial solutions (d) and (e).
\item{1.4)} $l_1=1$ and $l_2 = l_3=l_4=0$. In this case, as in [Cro], there 
is no solution, but the trivial one (d).
\item{2)} $l_1\geq2$, $l_1\geq l_2,l_3,l_4$. This case can be treated as the 
analogous one in [Cro].
This concludes the proof of Lemma 1.13.{\hfill QED}
\vskip 0.3 cm

We now use the above lemma to classify all the triangles which correspond
to good triples. Every quadruple generates twelve triangles. In 
fact, given a solution $(\varphi_1,\cdots,\varphi_4)$ we have six ways to 
choose the pair $(\varphi_i,\varphi_j)$ such that
$$
\cos 2 \pi \varphi_i+ \cos 2\pi \varphi_j= 2 \cos\pi(\varphi_i+\varphi_j)
\cos\pi(\varphi_i-\varphi_j).
$$
Chosen the pair $(\varphi_i,\varphi_j)$, we have two ways for choosing 
$\varphi_k$, in order to have the triangle
$$
\left(2\varphi_k, \varphi_i+\varphi_j, |\varphi_i-\varphi_j| \right).
\autoeqno{st}
$$
The remaining $\varphi_l$ is, by definition, such that the above triangle 
is mapped, by the braid \eqrefp{cl5.5}, to:
$$
\left( |\varphi_i-\varphi_j|,|1- \varphi_i-\varphi_j|,|1-2\varphi_l|\right).
$$
Let us analyze all the triangles generated by the solutions of the equation 
\eqrefp{cl7}, and keep the {\it good}\/ ones, namely the ones for which the 
new $r_k'$, given by \eqrefp{cl6}, is rational for every $i,j,k$, cyclic 
permutation of $1,2,3$.

In order to do this, observe that if there exists a permutation $p$ such that
the triple $(r_{p(1)},r_{p(2)},r_{p(3)})$ gives via \eqrefp{cl6.5} values of 
$\varphi_1,\varphi_2,\varphi_3$ such that there is not any rational 
$\varphi_4$ such that $\varphi_1,\varphi_2,\varphi_3,\varphi_4$ satisfy 
\eqrefp{cl7}, then $(r_1,r_2,r_3)$ is not a good triple. In fact, every 
permutation $p$ is generated by ciclic permutations and the permutation 
$p_{23}:(r_1,r_2,r_3)\to(r_1,r_3,r_2)$. Cyclic permutations are elements of 
the braid group, so the statement is obvious for them. For $p_{23}$, the 
statement is a trivial consequence of the fact that the triples 
$(r_1,r_2,r_3)$ and $(r_1,r_3,r_2)$ give via \eqrefp{cl6.5} the same values 
of $\varphi_1,\varphi_2,\varphi_3$. 

So we will exclude all the triangles $(r_1,r_2,r_3)$ for which there exists 
at least a permutation that gives rise to values of 
$(\varphi_1,\varphi_2,\varphi_3)$ for which rational solutions
$\varphi_4$ of \eqrefp{cl7} do not exist.

\vskip 0.3 cm
\noindent Solution (a). Using \eqrefp{st}, we obtain the triangles
$$
\left({1\over15}, {1\over30},{23\over30}\right),\quad
\left({1\over15}, {1\over5},{8\over15}\right),\quad
\left({1\over15}, {7\over30},{17\over30}\right),\quad
\left({4\over15}, {7\over30},{13\over30}\right),
$$
$$
\left({11\over15}, {11\over30},{13\over30}\right),\quad
\left({11\over15}, {2\over15},{1\over5}\right),\quad
\left({4\over5}, {2\over15},{1\over5}\right),\quad
\left({1\over5}, {1\over5},{7\over15}\right),
$$
$$
\left({1\over3}, {11\over30},{13\over30}\right),\quad
\left({2\over3}, {1\over30},{7\over30}\right),\quad
\left({1\over3}, {1\over5},{3\over5}\right),\quad
\left({1\over3}, {1\over3},{2\over5}\right).
$$
The last two points 
$$
\left({1\over3}, {1\over5},{3\over5}\right)\qquad
\left({1\over3}, {1\over3},{2\over5}\right)\autoeqno{cl8}
$$ 
belong to the orbit \eqrefp{cl5}.
The above values suitably permuted, except the \eqrefp{cl8}, give rise via 
\eqrefp{cl6.5}, to the following values of $(\varphi_1,\varphi_2,\varphi_3)$
(written in the same order as the correspondent generating triangles)
$$
\left({1\over60}, {5\over12},{7\over20}\right),\quad
\left({1\over10}, {3\over10},{7\over30}\right),\quad
\left({7\over60}, {19\over60},{1\over4}\right),\quad
\left({7\over60}, {7\over20},{1\over12}\right),\quad
\left({11\over60}, {7\over12},{3\over20}\right),
$$
$$
\left({1\over10}, {13\over30},{3\over10}\right),\quad
\left({1\over15}, {1\over2},{3\over10}\right),\quad
\left(0, {1\over5},{7\over30}\right),\quad
\left({1\over20}, {11\over60},{23\over60}\right),\quad
\left({1\over60}, {9\over20},{13\over60}\right).
$$
there isn't any rational number $\varphi_4$ such that any of the quadruples 
build with these triples and  $\varphi_4$ is in the class described by 
Lemma 1.13.
\vskip 0.3 cm
\noindent Solution (b). Using \eqrefp{st}, the triangles are
$$
\left({7\over15}, {11\over30},{23\over30}\right),\quad
\left({7\over15}, {2\over5},{11\over15}\right),\quad
\left({7\over15}, {1\over30},{11\over30}\right),\quad
\left({2\over15}, {1\over30},{19\over30}\right),
$$
$$
\left({2\over15}, {1\over30},{17\over30}\right),\quad
\left({2\over15}, {1\over15},{3\over5}\right),\quad
\left({2\over5}, {1\over15},{2\over5}\right),\quad
\left({2\over5}, {11\over15},{2\over5}\right),
$$
$$
\left({1\over3}, {1\over30},{13\over30}\right),\quad
\left({1\over3}, {11\over30},{23\over30}\right),\quad
\left({2\over5}, {1\over3},{4\over5}\right),\quad
\left({1\over3}, {1\over3},{4\over5}\right).
$$
The last two points are equivalent to:
$$
\left({1\over3}, {1\over5},{3\over5}\right)\qquad
\left({2\over3}, {1\over3},{1\over5}\right)
\autoeqno{cl9}
$$ 
of the orbit \eqrefp{cl5}. As before one can show that if $(r_1,r_2,r_3)$ is 
one of the above values, except the \eqrefp{cl9}, then there exists a 
permutation such that the $r'_k$ defined by \eqrefp{cl6} is no-more rational. 
In fact we obtain for example the following values of 
$(\varphi_1,\varphi_2,\varphi_3)$, which don't fall in the values obtained 
in Lemma 1.13:
$$
\left({11\over60},{37\over60},{3\over20}\right),\quad
\left({1\over5},{3\over5},{2\over15}\right),\quad
\left({1\over60},{5\over12},{1\over20}\right),\quad
\left({1\over60},{3\over4},{23\over60}\right),\quad
\left({1\over60},{47\over60},{7\over20}\right),
$$
$$
\left({1\over30},{23\over30},{11\over30}\right),\quad
\left({1\over30},{2\over5},0\right),\quad
\left({11\over30},{2\over5},0\right),\quad
\left({1\over60},{23\over60},{1\over20}\right),\quad
\left({11\over60},{11\over20},{13\over60}\right).
$$

\vskip 0.3 cm

\noindent Solution (c).
$$
\left({2\over7}, {1\over7},{5\over7}\right),\;
\left({2\over7}, {5\over42},{19\over42}\right),\;
\left({2\over7}, {11\over42},{25\over42}\right),\;
\left({4\over7}, {11\over42},{25\over42}\right),\;
\left({4\over7}, {2\over7},{4\over7}\right),\;
\left({4\over7}, {1\over42},{13\over42}\right),
$$
$$
\left({1\over7}, {1\over42},{29\over42}\right),\;
\left({1\over7}, {5\over42},{23\over42}\right),\;
\left({1\over7}, {1\over7},{4\over7}\right),\;
\left({1\over3}, {1\over7 },{ 3\over7 }\right),\;
\left({1\over3}, {2\over7},{4\over7}\right),\;
\left({1\over3}, {1\over7 },{ 5\over7 }\right).
$$
As before one can show that if $(r_1,r_2,r_3)$ is one of the above values  
then there exists a permutation such that the $r'_k$ defined by \eqrefp{cl6} 
is no more rational. In fact we obtain for example the following values of 
$(\varphi_1,\varphi_2,\varphi_3)$, which are not included in the values 
described by Lemma 1.13:
$$
\left({1\over14},{1\over2},{3\over14}\right),\quad
\left({5\over84},{31\over84},{1\over12}\right),\quad
\left({11\over84},{37\over84},{13\over84}\right),\quad
\left({11\over84},{7\over12},{1\over84}\right),
$$
$$
\left({1\over7},{4\over7},0\right),\quad
\left({1\over84},{37\over84},{11\over84}\right),\quad
\left({1\over84},{5\over12},{23\over84}\right),\quad
\left({5\over84},{29\over84},{17\over84}\right),
$$
$$
\left({1\over7},{2\over7},0\right),\quad
\left({1\over14},{8\over21},{1\over21}\right),\quad
\left({1\over7},{19\over42},{5\over42}\right),\quad
\left({1\over14},{11\over21},{4\over21}\right).
$$
\vskip 0.3 cm
\noindent Solution (d.1).\quad
$$
\left({2\over3}, {3\over20},{3\over20}\right),\quad
\left({1\over3}, {1\over20},{9\over20}\right),\quad
\left({1\over5}, {1\over20},{1\over20}\right),\quad
\left({1\over5}, {1\over30},{19\over30}\right),
$$
$$
\left({1\over5}, {1\over12},{7\over12}\right),\quad
\left({2\over5}, {1\over12},{5\over12}\right),\quad
\left({3\over5}, {3\over20},{7\over20}\right),\quad
\left({3\over5}, {7\over30},{13\over30}\right),
$$
$$
\left({1\over2}, {7\over30},{13\over30}\right),\quad
\left({1\over2}, {1\over30},{11\over30}\right),\quad
\left({2\over3}, {1\over5},{2\over5}\right),\quad
\left({1\over2}, {1\over5},{2\over5}\right).
$$
The last two points are equivalent to:
$$
\left({1\over3}, {1\over5},{3\over5}\right)\qquad
\left({1\over2}, {2\over5},{1\over5}\right)
\autoeqno{cl10}
$$ 
of the orbit \eqrefp{cl5}, which is now complete. We can again exclude all 
the other values of $(r_1,r_2,r_3)$, with the same trick as above.
\vskip 0.3 cm
\noindent Solution (d.2).\quad
In this case we have that any of the triangles generated is equivalent to 
one of the following:
$$
\left(|1-2 \varphi|, {|1-4\varphi|\over4},{1+4\varphi\over4} \right),\quad
\left({1\over2}, {|1-4\varphi|\over4},{|3-4\varphi|\over4} \right),\quad
\left({1\over2}, {|1-4\varphi|\over4},{1+4\varphi\over4} \right),\quad
$$
$$
\left(2 \varphi, {|1-4\varphi|\over4},{|3-4\varphi|\over4} \right),\quad
\left({1\over2},{1\over2}, {|1-4\varphi|\over2}\right),
$$
where $\varphi$ is an arbitrary rational number.
The last triangle is forbidden because it has two right angles, and the 
first four ones are all equivalent to a flat triangle, so they are again 
forbidden because they give rise to a an half-integer value of $\mu$.
\vskip 0.3 cm
\noindent Solution (d.3).\quad
The generated triangles are the following:
$$
\left({1\over2}, {1\over3},2\varphi+{1\over3} \right),\quad
\left(|{2\over3}-2\varphi|,{1\over3},2\varphi+{1\over3} \right),\quad
\left({1\over2},{2\over3},{2\over3}+2\varphi\right),\quad
\left({2\over3}+2\varphi,{2\over3},{2\over3}+2\varphi\right),
$$
$$
\left({2\over3}+2\varphi,{|1-4\varphi|\over4},{1+4\varphi\over4}\right),
\quad
\left(|{2\over3}-2\varphi|,{|1-4\varphi|\over4},{1+4\varphi\over4}\right),
\quad
\left({1\over2}, {1\over3},1+2\varphi\right),
$$
$$
\left(2\varphi,{1\over3},1+2\varphi\right),\quad
\left(|{2\over3}-2\varphi|,\varphi+{1\over12},\varphi+{7\over12}\right),\quad
\left(2\varphi,\varphi+{1\over12},\varphi+{7\over12}\right),
$$
$$
\left({2\over3}+2\varphi,\varphi+{5\over12},\varphi+{11\over12}\right),\quad
\left(2\varphi,\varphi+{5\over12},\varphi+{11\over12}\right).
$$
This case must be studied carefully because we have to classify the allowed 
values of the rational variable $\varphi$ in order that, applying the 
transformation \eqrefp{cl5.5}, we obtain always rational values. 

Let us analyze the first triangle. It is mapped by \eqrefp{cl5.5} to a 
triangle equivalent to the second:
$$
\left({1\over2}, {1\over3},2\varphi+{1\over3} \right)\mapsto
\left({1\over3}+2\varphi,{2\over3},2\varphi+{1\over3} \right)\sim
\left(|{2\over3}-2\varphi|,{1\over3},2\varphi+{1\over3} \right).
\autoeqno{cl11}
$$
If we apply the braid \eqrefp{cl5.5}, with $r_i={2\pi\over3}$, 
$r_j=r_k=2\varphi+{1\over3}$, we have to solve:
$$
\cos{2\pi\over3}+2\cos^2\pi(2\varphi+{1\over3}) = \cos\pi r_k',
$$
or, equivalently, for $\varphi$ and for the new $\varphi_k$:
$$
\cos{2\pi\over3}+\cos2\pi(2\varphi+{1\over3})+1+ \cos2\pi\varphi_k=0.
$$
We classify the values of the allowed $\varphi$ using the Lemma 1.13 in the 
case (e). We have six possibilities for $\varphi$:

\item{i)} if $\varphi_k={1\over2}$, then $\varphi={11\over24}$. In this case 
we obtain, from \eqrefp{cl11}, all the points of the orbit \eqrefp{cl2}.

\item{ii)} if $\varphi_k={1\over4}$, then $\varphi=0$. In this case we
obtain, from \eqrefp{cl11}, all the points of the orbit \eqrefp{cl1}.

\item{iii)} if $\varphi_k={1\over2}$, then $\varphi={1\over4}$. In this case 
we obtain, from \eqrefp{cl11}, the following points:
$$
\left({1\over2}, {1\over3},{5\over6}\right),\qquad
\left({1\over6},{2\over3},{1\over6} \right). 
$$
We exclude them because there exists a permutation such that the $r'_k$ 
defined by \eqrefp{cl6} is no-more rational.

\item{iv)} if $\varphi_k$ is free to vary, then 
$2\varphi+{1\over3}={1\over2}$. In this case we obtain, from \eqrefp{cl11}, 
the forbidden point:
$$
\left({1\over2}, {1\over3},{1\over2}\right),
$$
so we have to exclude it.

\item{v)} if $\varphi_k={1\over5}$, then $\varphi={1\over30}$, and we 
obtain, from \eqrefp{cl11}, the following two points of the orbit 
\eqrefp{cl4}:
$$
\left({1\over2}, {1\over3},{2\over5}\right),\qquad
\left({2\over5},{2\over3},{2\over5} \right). 
$$

\item{vi)} if $\varphi_k={2\over5}$, then $\varphi={13\over30}$, and we 
obtain the following two points of the orbit \eqrefp{cl3}:
$$
\left({1\over2}, {1\over3},{1\over5}\right),\qquad
\left({1\over5},{2\over3},{1\over5} \right). 
$$

In the same way we can study all the other triangles and show that we don't 
obtain any other value but the ones described in Theorem 1.6. 
\vskip 0.3 cm
\noindent Solution (e.1).\quad
The generated triangles are the following:
$$
\left(0,0,{2\over3}\right),\,
\left(0,{1\over12},{7\over12}\right),\,
\left(0,{1\over2},{2\over3}\right),\,
\left({2\over3},{1\over12},{7\over12}\right)\,
\left({2\over3},{1\over3},{1\over3}\right),\,
\left({1\over2},{1\over3},{1\over3}\right),\,
\left({2\over3},{1\over4},{1\over4}\right);
$$
the first four are forbidden because there exists a permutation such that the 
$r'_k$ defined by \eqrefp{cl6} is no-more rational. So we have to exclude 
them. The fifth and the sixth are points of the orbit \eqrefp{cl1} and the 
last of \eqrefp{cl2}. 
\vskip 0.3 cm
\noindent Solution (e.2).\quad
The generated triangles are the following:
$$
\left(2\varphi,{1\over2},{1\over2} \right),\quad
\left(|1-2\varphi|, {1\over2},{1\over2}\right),\quad
\left({1\over2},{1\over2},{|1+4\varphi|\over2}\right),\quad
\left(1,\varphi,\varphi\right),
$$
$$
\left(1,{|1-2\varphi|\over2},{|1-2\varphi|\over2}\right),\quad
\left(0,{1\over2},{|1-4\varphi|\over2}\right),\quad
\left(|1-2\varphi|,\varphi,\varphi\right),
\left(0,{|1-2\varphi|\over2},{1+2\varphi\over2}\right),
$$
$$
\left(0,\varphi,1-\varphi\right),\;
\left(2\varphi,{|1-2\varphi|\over2},{|1-2\varphi|\over2}\right),\;
\left(|1-2\varphi|,{|1-2\varphi|\over2},{1+2\varphi\over2}\right),\;
\left(2\varphi,\varphi,1-\varphi\right),
$$
They are all forbidden, the first three because they have two right angles, 
the next three ones, because we can prove that necessarily
$\varphi={1\over2}$, then the first has two right angles, the second 
gives $|\cos\pi r_k'|=3$ and the last gives $|\cos\pi r_k'|=2$; all the 
others because they are equivalent to a flat triangle. 

\vskip 0.3 cm
\noindent Solution (e.3).\quad
The generated triangles are the following:
$$
\left(0,{2\over15},{8\over15}\right),\;
\left(0,{1\over15},{11\over15}\right),\;
\left(0,{1\over5},{3\over5}\right),\;
\left({2\over3},{1\over5},{3\over5}\right),\;
\left({2\over5},{1\over15},{11\over15}\right),\;
\left({4\over5},{2\over15},{8\over15}\right),
$$
$$
\left({4\over5},{1\over5},{1\over5}\right),\;
\left({2\over3},{1\over5},{1\over5}\right),\;
\left({2\over3},{2\over5},{2\over5}\right),\;
\left({2\over5},{2\over5},{2\over5}\right),\;
\left({2\over5},{1\over3},{1\over3}\right),\;
\left({1\over5},{2\over3},{1\over3}\right).
$$
We exclude the first six because we can show that there exists a 
permutation such that the $r'_k$ defined by \eqrefp{cl6} is no-more 
rational. The seventh and the eighth give two points of the orbit 
\eqrefp{cl4}, the ninth and tenth two points of the orbit \eqrefp{cl6} 
and the last two, two points of \eqrefp{cl6}.
\vskip 0.3 cm
\noindent Solution (f)\quad
We have obtained all the points of all the orbits of Theorem 1.6. To show 
that there are no other points we still have to examine the case (f). In 
this case all the obtained triangles are equivalent to the following:
$$
\eqalign{
\left({\varphi_1\over2},{|2-\varphi_1+\varphi_3|\over4}, 
{|2-\varphi_1-\varphi_3|\over4}\right),&
\left({\varphi_1\over2},{1\over2},{|1-\varphi_1|\over2}\right),\cr
\left({\varphi_1\over2},
{|4-\varphi_1-\varphi_3|\over4},\right.
&\left.{|\varphi_1-\varphi_3|\over4}\right).\cr}
\autoeqno{cl16}
$$
Applying the transformation \eqrefp{cl5.5} to the above triangles, we find 
that we have to solve for $\varphi_1$, $\varphi_3$ and for the new $\varphi$ 
obtained from the \eqrefp{cl6}, the following three equations respectively:
$$
\cos\pi{|2-\varphi_3|\over2}+\cos\pi\varphi_1+\cos2\pi\varphi+1=0,
$$
$$
\cos\pi{|1+\varphi_1-\varphi_3|\over2}+
\cos\pi{|1-\varphi_1-\varphi_3|\over2}+\cos2\pi\varphi=0,
$$
$$
2\cos\pi{\varphi_1+\varphi_3\over4}+\cos\pi{|3\varphi_1-\varphi_3|\over4}
+\cos2\pi\varphi=0.
$$
We again can use Lemma 1.13 to prove that we don't obtain any new point. 
Let us show this in for the first triangle
$$
\left({\varphi_1\over2},{|2-\varphi_1+\varphi_3|\over4}, 
{|2-\varphi_1-\varphi_3|\over4}\right).\autoeqno{cl13}
$$
We have to solve the equation
$$
\cos\pi{|2-\varphi_3|\over2}+\cos\pi\varphi_1+\cos2\pi\varphi+1=0.
\autoeqno{cl14}
$$
Using Lemma 1.13, the possible values for 
$({|2-\varphi_3|\over2},\varphi_1,\varphi)$ are the (e.1), (e,2) and (e.3). 
Consider the case (e.1), then the possible solutions for the pair 
$(\varphi_1,\varphi_3)$, are
$$
(\varphi_1,\varphi_3)= ({1\over 2},{2\over 3}),\quad
({2\over 3},1),\quad
({2\over 3},{2\over 3}).
$$
Let us substitute these solutions in \eqrefp{cl13}; we obtain the triangles
$$
({1\over4},{13\over 24},{5\over 24}),\quad
({1\over 3},{7\over12},{1\over12}),\quad
({1\over 3},{1\over2},{1\over6}),
$$
which are all flat, and thus forbidden. Let us now consider the case (e.2).
In this case we obtain two possibilities: either $\varphi_3=0$ and 
$\varphi_1$ is a free parameter, or $\varphi_1=1$ and $\varphi_1$ is 
a free parameter. In both the cases the triangle \eqrefp{cl13} is flat, 
and thus forbidden. Let us now consider the last case (e.3). The possible 
solutions for the pair $(\varphi_1,\varphi_3)$, are
$$
(\varphi_1,\varphi_3)=({2\over 5},{2\over 3}),\quad
({2\over 3},{6\over5}),\quad
({2\over 5},{2\over 5}),\quad
({4\over 5},{6\over5}),\quad
({4\over 5},{2\over 3}),\quad
({2\over 3},{2\over5}).
$$
Substituting these values in \eqrefp{cl13}, we obtain all flat angles. We 
can repeat the same proof for the other two triangles in \eqrefp{cl16}.
In this way we conclude the proof of the theorem.{\hfill QED}
\vskip 0.3 cm

\noindent{\bf 1.4. Monodromy data and reflection groups.}
\vskip 0.3 cm
We reformulate here the above parameterization of the monodromy data by 
classes of equivalence of triples $(x_1,x_2,x_3)$ in a geometric way.
Let us consider a three-dimensional space $V$ with a basis $(e_1,e_2,e_3)$ 
and with a symmetric bilinear form $(\cdot,\cdot)$ given, in this basis, by 
the matrix
$$
g:=\pmatrix{2&x_1&x_3\cr
x_1&2&x_2\cr
x_3&x_2&2\cr}\autoeqno{N29}
$$
namely
$$
(e_i,e_i)=2,\quad\hbox{for}\quad i=1,2,3,\quad\hbox{and}\quad
(e_1,e_2)=x_1,\quad (e_2,e_3)=x_2,\quad (e_1,e_3)=x_3.
$$
Observe that the bilinear form \eqrefp{N29} does not degenerate. Indeed,
$$
\det g = 8- 2(x_1^2+x_2^2+x_3^2-x_1 x_2 x_3)= 8\cos^2\pi \mu\neq 0,
$$
due to the non-resonance assumption $2\mu\not\in\interi$. The three planes 
$p_1,p_2,p_3$ orthogonal to the basic vectors $(e_1,e_2,e_3)$ possess the 
following properties:
\item{1)} The normal vectors to these planes are non-isotropic 
(i.e. $(e_i,e_i)\neq0$).
\item{2)} None of the planes is orthogonal to the other two.

Conversely, a three-dimensional space $V$ with a non-degenerate 
symmetric bilinear form $(\cdot,\cdot)$ and with an ordered triple of 
planes satisfying the above conditions, uniquely determines the matrix 
$g$ of the form \eqrefp{N29}, and then the monodromy data of a solution 
of PVI$\mu$.

We define three reflections $R_1,R_2,R_3$ with respect to the three planes 
$(p_1,p_2,p_3)$:
$$
R_i:\eqalign{
V&\rightarrow V\cr
x&\mapsto x-(e_i,x)e_i\cr}
\qquad i=1,2,3.
$$
These reflection have the following matrix representation in the basis 
$(e_1,e_2,e_3)$:
$$
R_1=\pmatrix{-1&-x_1&-x_3\cr 0&1&0\cr 0&0&1\cr},\quad
R_2=\pmatrix{1&0&0\cr -x_1&-1&-x_2\cr 0&0&1\cr},\quad
R_3=\pmatrix{1&0&0\cr 0&1&0\cr -x_3&-x_2&-1\cr}.\autoeqno{N20.5}
$$
Let us consider the group $G\subset{\cal O}(V,(\cdot,\cdot))$ of the linear 
transformations of $V$, generated by the three reflections $R_1$, $R_2$, 
$R_3$. 
The matrix $g$ will be called the {\it Gram matrix}\/ of the reflection 
group $G$. It determines the subgroup $G\subset{\cal O}(V,(\cdot,\cdot))$
uniquely. We observe that, for an admissible triple, the group $G$ is 
irreducible, namely there are no non-trivial subspaces of $V$ which are 
invariant with respect to all the transformations of $G$. 

We conclude that the branches of the solutions of  PVI$\mu$ can be 
parameterized by groups $G\subset{\cal O}(3)$ with a marked ordered system of 
generating reflections $R_1,R_2,R_3$. Let us describe what happens with the 
triples of generators under the analytic continuation of the solution.

We define an action of the braid group $B_3$ on the systems of generators 
$R_1,R_2,R_3$ of the reflection group $G$:
$$
\eqalign{
\beta_1:(R_1,R_2,R_3)\mapsto 
(R_1,R_2,R_3)^{\beta_1}:=&(R_2,R_2R_1R_2,R_3),\cr
\beta_2:(R_1,R_2,R_3)\mapsto 
(R_1,R_2,R_3)^{\beta_2}:=&(R_1,R_3, R_3R_2R_3),\cr}\autoeqno{B19}
$$
where $\beta_{1,2}$ are the standard generators of the braid group. Observe 
that the groups generated by the reflections $(R_1,R_2,R_3)$ and 
$(R_1,R_2,R_3)^{\beta}$ coincide for any $\beta\in B_3$. In particular the 
following lemma holds true:

\proclaim Lemma 1.15. For any braid $\beta\in B_3$, the transformations 
$\beta(R_1,R_2,R_3)$ are reflections with respect to some planes orthogonal 
to some new basic vectors $(e^\beta_1,e^\beta_2,e^\beta_3)$. The Gram matrix 
with respect to the basis $(e^\beta_1,e^\beta_2,e^\beta_3)$ has the form:
$$
(e^\beta_i,e^\beta_i)=2,\quad i=1,2,3,\quad
(e^\beta_1,e^\beta_2)=x^\beta_1,\quad (e^\beta_2,e^\beta_3)=x^\beta_2,
\quad (e^\beta_1,e^\beta_3)=x^\beta_3,
$$
where $(x^\beta_1,x^\beta_2,x^\beta_3)=\beta(x_1,x_2,x_3)$.

\noindent Proof. It is sufficient to check the statement 
for the generators $\beta_{1,2}$. For $\beta=\beta_1$:
$$
e^{\beta_1}_1=e_2,\quad  e^{\beta_1}_2=e_1-x_1 e_2 ,\quad  
e^{\beta_1}_3=e_3,
$$
for $\beta=\beta_2$:
$$
e^{\beta_2}_1=e_1,\quad e^{\beta_2}_2=e_3,\quad e^{\beta_2}_3=e_2-x_2 e_3.  
$$
Computing the Gram matrix we prove the lemma. {\hfill QED}

\vskip 0.3 cm
\noindent{\bf 1.4.1. Reflection groups and algebraic solutions}\quad Let us 
figure out what are the reflection groups corresponding to the finite orbits 
classified in Theorem 1.6.

\proclaim Theorem 1.7. The orbit \eqrefp{cl1} corresponds to the group 
$W(A_3)$ of symmetries of regular tetrahedron, the orbit \eqrefp{cl2}  
corresponds to the group $W(B_3)$ of symmetries of the regular octahedron,
the orbits \eqrefp{cl3}, \eqrefp{cl4}, \eqrefp{cl5} correspond to different
choices of a system of generating reflections in the group $W(H_3)$ of 
symmetries of icosahedron.

\noindent Proof. It is sufficient to find one point in each of the orbits
\eqrefp{cl1}, \eqrefp{cl2}, \eqrefp{cl3}, \eqrefp{cl4}, \eqrefp{cl5} that 
corresponds to a triple of symmetry planes of a regular polyhedron. To this
end, we associate to a regular polyhedron a standard triple of symmetry 
planes using the following construction. Let $0$ be the center of the 
polyhedron. Take a face of the polyhedron and denote $H$ the center of this 
face, $P$ a vertex and $Q$ the center of an edge of the face passing through 
the vertex $P$. The standard triple consists of the symmetry planes trough
the points $OPQ$, $OQH$, $OHP$ respectively. Let us compute the angles 
between the planes of each regular 
polyhedron. It is convinient to use the Schl\"afli symbol $\{p,q\}$ for 
regular polyhedra (see [Cox]). In these notations, the face of the regular 
polyhedron $\{p,q\}$ is a regular $p$-gon, the vertex figure is a regular 
$q$-gon. We immediately see that the angles between the planes of the 
standard triple are
$$
\hbox{between}\quad OPQ\quad \hbox{and}\quad OQH \qquad {\pi\over 2},
$$
$$
\hbox{between}\quad OQH\quad \hbox{and}\quad OHP \qquad {\pi\over p},
$$
$$
\hbox{between}\quad OHP\quad  \hbox{and}\quad OPQ \qquad  {\pi\over q}.
$$
So, for the tetrahedron $\{3,3\}$ we obtain the angles 
$\left({\pi\over 2},{\pi\over 3},{\pi\over 3}\right)$, for the octahedron
$\{3,4\}$ the angles $\left({\pi\over 2},{\pi\over 3},{\pi\over 4}\right)$,
for the icosahedron $\{3,5\}$ the angles 
$\left({\pi\over 2},{\pi\over 3},{\pi\over 5}\right)$. In this way, we 
obtain the triples \eqrefp{cl1}, \eqrefp{cl2}, \eqrefp{cl3}. The reciprocal
polyherda (i.e. cube $\{4,3\}$ and dodecaherdon $\{5,3\}$) give the same 
angles up to permutations. As we already know, the permuted triples of 
\eqrefp{cl2} or \eqrefp{cl3} belong to the same orbit. So, the standard 
triples of cube and dodecahedron are $B_3$-equivalent to those of octahedron 
and icosahedron respectively.

To obtain the last two orbits \eqrefp{cl4} and \eqrefp{cl5}, we apply the 
above construction of the standard triple to great icosahedron and 
dodecahedron respectively. These non convex regular polyhedra both have
icosaherdal symmetry (see [Cox]). Their Schl\"afli symbols are 
$\left\{3,{5\over2}\right\}$ and $\left\{5,{5\over2}\right\}$ respectively.
This means that the faces of these polyhedra are regular triangles or
pentagons, but the vertex figures are pentagrams. The above computation
gives the triples \eqrefp{cl4} and \eqrefp{cl5}. Again we need not consider 
the reciprocal stellated polyhedra. Theorem 1.7 is proved.{\hfill QED}

\vskip 0.3 cm
\noindent{\bf 1.4.2. Classification of the monodromy data, second proof.}\quad
We present here another proof of Theorem 1.6, based on 
the idea suggested by E. B. Vinberg. We start with the following:

\proclaim Algebraic Lemma. Let $(x,y,z)$ be an admissible triple of real 
numbers, satisfying the inequalities:
$$
x^2+y^2+z^2-x y z >4,\autoeqno{N21}
$$
and
$$
|x|,\, |y|,\, |z|\, \leq 2.\autoeqno{N22}
$$
Then there exists a braid $\beta\in B_3$ such that the absolute value of some 
of the coordinates of $\beta(x,y,z)$ is strictly greater than $2$.

Before proving the lemma, we observe that we can assume, without loss of 
generality, that all the coordinates of $(x,y,z)$ are non-zero; in fact, 
for any admissible triple, there exists a braid $\beta\in B_3$ such that 
all the coordinates of $\beta(x,y,z)$ are non-zero. Let us denote $b_x$, 
$b_y$ and $b_z$ the following braids:
$$
\matrix{
b_x:=&\beta_2, & b_x(x,y,z)=(z,-x,x-y z),\cr
b_y:=&\beta_2^{-1}\beta_1\beta_2, & b_y(x,y,z)=(-y+x z,-x,-z),\cr
b_z:=&\beta_1, & b_z(x,y,z)=(-x,z-x y,y).\cr}
$$

\proclaim Lemma 1.16. Let $(x,y,z)$ be a triple of non-zero real numbers, 
satisfying 
$$
0<|z|,\, |x|,\, |y|\leq 2\autoeqno{N22.5}
$$
and
$$
x^2+y^2+z^2-x y z = 4+ c^2,\quad c>0. \autoeqno{N23}
$$
Denote $(x',y',z'):=\beta(x,y,z)$, where
$$
\beta=\left\{\eqalign{
b_x\quad\hbox{if}\quad & |x|\leq |y|,\, |z|,\cr
b_y\quad\hbox{if}\quad & |y|\leq |x|,\, |z|,\cr
b_z\quad\hbox{if}\quad & |z|\leq |x|,\, |y|.\cr}
\right.
$$
Then:
$$
\min\{ |x'|,\, |y'|,\, |z'|\}\geq\min\{ |x|,\, |y|,\, |z|\}\autoeqno{N24}
$$
and
$$
|x'|+|y'|+ |z'|\geq|x|+|y|+ |z|+ \min\{z^2, 2 c\}.\autoeqno{N25}
$$

\noindent Proof. Let us prove the lemma in the case where 
$|z|\leq |x|,\, |y|$ and $\beta=b_z$. The other cases can be proved in the 
same way. If the signs of $z$ and of $x y$ are opposite then
$$
|y'|=|z|+|x y|\geq |z|+ z^2,\quad |x'|=|x|, \quad |z'|=|z|
$$
and \eqrefp{N24}, \eqrefp{N25} are proved.
Let us suppose that the signs of $z$ and of $x y$ are the same. Changing the 
triple $(x,y,z)$ to an equivalent one, we can assume that all the coordinates 
are positive. If we prove now that
$$
2 z+ 2 c\leq x y,\autoeqno{N25.5}
$$
where $c$ is given in \eqrefp{N23}, we have that $|y'|=|x y-z|\geq z + 2c$ 
and the lemma is proved. To prove \eqrefp{N25.5} we find the constrained 
minimum of the function $x y$ on the domain $D$ defined by the conditions 
\eqrefp{N22.5} and \eqrefp{N23}. The Lagrange function
$$
F(x,y,z):= x y -\lambda\left(x^2+y^2+z^2-x y z \right),
$$
has the local maximum at
$$
x= y = \sqrt{4+c^2-z^2\over 2-z},
$$
and no minimum in the interior of $D$.
It remains to study the values of the function $x y$ on the boundary of $D$.
If, say, $z=y$ then the positive root $x$ of the equation
$$
x^2+2 z^2-x z^2=4 + c^2
$$
is greater than $2$. So the boundaries $z=y$ and $z=x$ are not reached for 
$(x,y,z)\in D$, and then $|z|< |x|,\, |y|$. It remains the last boundary 
to be studied. If, say, $y=2$, we find $x=z\pm c$. Since $x\geq z$, then 
$x=z+c$ and $x y=2(z+c)$; this is the minimum of the function $x y$.
{\hfill QED}
\vskip 0.3 cm

\noindent Proof of Algebraic Lemma. 
As observed above we can always reduce to the case where all the 
coordinates $(x,y,z)$ are non-zero. Put:
$$
\Delta(x,y,z):= \min \left\{ 
x^2,\, y^2, \, z^2,\, 2\sqrt{x^2+y^2+z^2-x y z-4}
\right\}.
$$
Using Lemma 1.16, we can build a braid $b_1$ such that the coordinates:
$$
(x_1,y_1,z_1):=b_1(x,y,z)
$$
satisfy the inequalities
$$
\min\{ |x_1|,\, |y_1|,\, |z_1|\}\geq\min\{ |x|,\, |y|,\, |z|\}
\quad
|x_1|+|y_1|+ |z_1|\geq|x|+|y|+ |z|+ \Delta(x,y,z).\autoeqno{N26}
$$
Since the quantity $x^2+y^2+z^2-x y z-4$ is preserved by the action of the 
braid group, we obtain:
$$
\Delta(x_1,y_1,z_1)\geq\Delta(x,y,z).
$$
If the absolute value of some of the coordinates $(x_1,y_1,z_1)$ is greater 
than 2, the lemma in proved. Otherwise we apply again the construction of 
Lemma 1.16 to the triple $(x_1,y_1,z_1)$. In this way we obtain a 
sequence of braids $b_1,\, b_2,\, b_3 \cdots$ such that the corresponding 
triples
$$
(x_{k+1},y_{k+1},z_{k+1}):=b_{k+1}(x_k,y_k,z_k)
$$ 
satisfy
$$
\Delta(x_{k+1},y_{k+1},z_{k+1})\geq\Delta(x_k,y_k,z_k).
$$
Iterating the inequality \eqrefp{N24}, we obtain that
$$
|x_k|+|y_k|+|z_k|\geq |x|+|y|+ |z|+ k \Delta(x,y,z).
$$
Hence, in a finite number of steps we build a triple such that the absolute 
value of at least one of the coordinates in greater than 2. This concludes 
the proof of Algebraic Lemma.{\hfill QED}
\vskip 0.3 cm

\proclaim Corollary 1.4. For an algebraic solution of PVI$\mu$, specified by 
an admissible triple $x_i=-2\cos2\pi r_i$, the value of $\mu$ must be real, 
the strict inequalities
$$
|x_i|<2,\quad i=1,2,3,\autoeqno{28.5}
$$
hold true and the matrix $g$ defined in \eqrefp{N29} is positive definite.

\noindent Proof. 
Let us prove that, for an algebraic solution, the triple $(x_1,x_2,x_3)$ must 
satisfy the inequality:
$$
x_1^2+x_2^2+x_3^2-x_1 x_2 x_3 < 4.\autoeqno{N28}
$$
Indeed, if $x_1^2+x_2^2+x_3^2-x_1 x_2 x_3  > 4$, then, according to the 
Algebraic Lemma the triple is not a good one. This contradicts the 
assumption that the solution is algebraic. If 
$x_1^2+x_2^2+x_3^2-x_1 x_2 x_3 = 4$, then $\mu={1\over2}+k$ with 
$k\in\interi$. This contradicts the basic assumption $2\mu\not\in\interi$. 
Then \eqrefp{N28} is satisfied and $\mu$ is a real number. Let us now prove 
\eqrefp{28.5}. If one of the coordinates, say $x_1$, is such that 
$x_1=\pm 2$, then
$$
x_1^2+x_2^2+x_3^2-x_1 x_2 x_3 = 4 + (x_2\mp x_3)^2,
$$
and, being $x_2$ $x_3$ real numbers, \eqrefp{N28} is violated. So,
$x_i\neq\pm 2$ for every $i$. Finally, applying the Sylvester criterion 
to the matrix $g$, we prove that $g$ is positive definite. In fact 
$$
\det G=8-2(x_1^2+x_2^2+x_3^2-x_1 x_2 x_3 )>0,
$$
and for any principal minor 
$$
\det\pmatrix{2&x_i\cr
x_i&2\cr}=4 - x_i^2>0.
$$
{\hfill QED}
\vskip 0.2 cm

\proclaim Lemma 1.17. For an algebraic solution of PVI$\mu$ the reflection 
group $G$ acts in the Euclidean space.

\noindent The proof immediately follows from the fact that the correspondent 
Gram matrix is positive definite.

\proclaim Corollary 1.5. For a good triple $(x_1,x_2,x_3)$ and for any braid 
$\beta\in B_3$, there exists three integer positive numbers $n_{12}^\beta$, 
$n_{13}^\beta$ and $n_{23}^\beta$ such that:
$$
\left(R_i^\beta R_j^\beta\right)^{n_{ij}^\beta}=\ID,\quad\hbox{for}
\quad i\neq j, \quad i,j=1,2,3.\autoeqno{N30}
$$

\noindent Proof. If $(e_1,e_2)=x_1=-2\cos\pi r$ with 
$r={m\over n}$, $m,n\in\interi$, then $R_1 R_2$ is a rotation by the angle 
$2\pi{m\over n}$. Hence:
$$
\left(R_1 R_2\right)^{n}=\ID.
$$
This holds true for any pair $R_i$ and $R_j$. Moreover, for any braid 
$\beta\in B_3$, the triple $\beta(x_1,x_2,x_3)$ is again good, then 
\eqrefp{N30} is proved.{\hfill QED}

\proclaim Corollary 1.6. The set of the solutions of the PVI$\mu$ equation 
with a real non resonant value of $\mu$ and real parameters $(x_1,x_2,x_3)$ 
satisfying 
$$
|x_i|<2,\quad i=1,2,3,
$$
is invariant with respect to the analytic continuation.

\noindent Proof. Applying the Sylvester criterion to 
the matrix $g$ defined in \eqrefp{N29}, we obtain that $g$ is positive 
definite. So the reflections $R_1,R_2,R_3$ can be realized in the Euclidean 
space.
After a transformation $(x_1,x_2,x_3)\mapsto(x_1^\beta,x_2^\beta,x_3^\beta)=
\beta(x_1,x_2,x_3)$, the new numbers $({x_1}^\beta,{x_2}^\beta,{x_3}^\beta)$ 
are the entries of the Gram matrix:
$$
g^\beta:=\pmatrix{2&x^\beta_1&x^\beta_3\cr
x_1^\beta&2&x^\beta_2\cr
x^\beta_3&x^\beta_2&2\cr},
$$
of the basis $(e^\beta_1,e^\beta_2,e^\beta_3)$, in the same Euclidean space. 
Then this matrix must be positive definite, namely $x_i^2<4$ as we wanted to 
prove. {\hfill QED}
\vskip 0.3 cm

In the second part of this paper, we will identify the set described in  
Corollary 1.6 with the class of solutions of PVI$\mu$ having asymptotic 
behaviour of algebraic type. This identification will be crucial in the 
computation of the five algebraic solutions of PVI$\mu$ we have classified.

As we have just shown, a good triple 
$$
(x_1,x_2,x_3)=\left(-2\cos\pi{ m_1\over n_1},-2\cos\pi{ m_2\over n_2},
-2\cos\pi{ m_3\over n_3}\right),
$$
corresponds to a representation of the Coxeter group generated by three 
reflections $R_1$, $R_2$, $R_3$ satisfying
$$
R_1^2=R_2^2=R_3^2=\ID,\quad 
(R_1 R_2)^{n_1}=(R_2 R_3)^{n_2}=(R_1 R_3)^{n_3}=\ID,\autoeqno{N31}
$$
in the three-dimensional Euclidean space. We denoted by $G$ the image of 
this representation. Moreover, for any braid $\beta\in B_3$, the matrices
$$
(R^\beta_1,R^\beta_2,R^\beta_3):=\beta(R_1,R_2,R_3),
$$
satisfy the same identities \eqrefp{N31}, with some new integers 
$n^\beta_1,n^\beta_2,n^\beta_3$. We stress that the reflections 
$R^\beta_1,R^\beta_2,R^\beta_3$ generate the same group $G$.

\proclaim Theorem 1.8. It follows from the above property that $G$ is an 
irreducible finite Coxeter group.
 
Let $n$ be the least common multiple of $n_1$, $n_2$ and $n_3$.  Put:
$$
\zeta=2 \cos{\pi\over n}.
$$

\proclaim Lemma 1.18. The numbers 
$$
x_i=-2\cos\pi {m_i\over n_i},\quad i=1,2,3,
$$
belong to the ring ${\cal K}_0$ of integers of the field 
${\cal K}:={\bf Q}[\zeta]$.

Recall (see [Wey]) that ${\cal K}$ is the normal extension of ${\bf Q}$ 
generated by $\zeta$ and ${\cal K}_0$ is the ring of all the algebraic 
integer numbers of ${\cal K}$, namely it consists of all the elements 
$x\in{\cal K}$ satisfying an algebraic equation of the form 
$$
x^k+a_1 x^{k-1} + \cdots + a_k=0, \quad\hbox{with}\quad
a_i\in\interi.
$$

\vskip 0.3 cm
\noindent Proof of Lemma 1.18.\quad Let $n=n_i m_i'$, then 
$$
\cos \pi{m_i\over n_i}=T_{m_im_i'}\left(\cos{\pi\over n}\right)
=T_{m_im_i'}\left({1\over2}\zeta\right),
$$
where
$$
T_k(x)=\cos(k\arccos x)=2^{k-1} x^k +\sum_{s=0}^{k-1}2^{s-1} a_{k s} x^s,
\autoeqno{N32}
$$
are the Tchebyscheff polynomials of the first kind (see [Bat]). Recall that 
all the coefficients $a_{k s}$ are integers, so $x_i=-2\cos \pi{m_i\over n_i}$ 
is a polynomial of $\zeta$ with integer coefficients. Moreover $\zeta$ is a 
root of the monic algebraic equation with integer coefficients
$$
2 T_n\left({\zeta\over2}\right) +2=
\zeta^n+\sum_{s=0}^{n-1} a_{n s}\zeta^s+2 =0.
$$
Hence $\zeta\in{\cal K}_0$ and 
$x_i=-2 T_{m_im_i'}\left({\zeta\over2}\right)\in{\cal K}_0$, as we wanted to 
prove. 
{\hfill QED}
\vskip 0.3 cm

\noindent Proof of Theorem 1.8.\quad
From the formulae \eqrefp{N20.5} it follows that the matrices  $R_1$, 
$R_2$ and $R_3$ are all defined over the same ring ${\cal K}_0$ of  
integers of ${\cal K}$:
$$
R_i\in{\rm Mat}({\cal K}_0,3).
$$
Moreover, these matrices are orthogonal with respect to $g$:
$$
R_i^{T}g R_i=g,\autoeqno{N33}
$$ 
where $g$ is defined in \eqrefp{N29}. Let 
$$
\Gamma:={\it Gal}({\cal K},{\bf Q})
$$
the Galois group of ${\cal K}$ over ${\bf Q}$, namely the group of all 
automorphisms
$$
\phi:{\cal K}\rightarrow {\cal K},
$$
identical on ${\bf Q}$. 

For any $\phi\in\Gamma$ we denote $\phi(R_i)$ and $\phi(g)$ the matrices 
obtained from $R_i$ and $g$ by the action 
$$
(x_1,x_2,x_3)\mapsto(\phi(x_1),\phi(x_2),\phi(x_3)).
\autoeqno{N34}
$$

\proclaim Lemma 1.19. For any $\phi\in\Gamma$ the following statements hold 
true:
\item{i)} $\det \phi(g)\neq 0$,
\item{ii)} The matrices $\phi(R_i)$ are orthogonal with respect to $\phi(g)$.
\item{iii)} For any $\beta\in B_3$ the matrices $\phi(R_i)^\beta$ satisfy 
the Coxeter relation \eqrefp{N30}.

\noindent The proof is obvious, due to the fact that any automorphism 
preserves all the algebraic relations.

\vskip 0.2 cm
From the above lemma, and from Algebraic Lemma, it follows that for 
any $\phi\in\Gamma$, the real symmetric matrix $\phi(g)$ must be positive 
definite. We will show that this implies that the group $G$ is finite. 
Let $N$ be the order of the Galois group $\Gamma$. We construct the 
block-diagonal matrices
$$
{\cal R}_i\in{\rm Mat}({\cal K}_0,3N), \quad i=1,2,3,
$$
as the matrices formed by $3\times 3$ blocks on the diagonal, such that the 
$j$-th block is $\phi_j(R_i)$, for $\phi_j\in\Gamma$, $j=1,2,\cdots, N$. The 
matrices ${\cal R}_i$ are orthogonal with respect to ${\cal G}$, that is 
the block-diagonal matrix having $\phi_j(G)$, for $\phi_j\in\Gamma$, 
$j=1,2,\cdots, N$, on the diagonal blocks.
We can apply Lemma 1.19 to the matrices ${\cal R}_i$ to show that 
they satisfy the Coxeter relation \eqrefp{N30}.
As a consequence we obtain a representation of our reflection group $G$ 
into the orthogonal group
$$
\eqalign{
G&\rightarrow{\cal O}\left({\cal K}^{3N},{\cal G}\right)\cr
R_i&\mapsto{\cal R}_i.\cr}\autoeqno{N35}
$$
By construction the matrices ${\cal R}_i$ preserve the sublattice
$$
{\cal K}_0^{3N}\subset{\cal K}^{3N}
$$
of the vectors the components of which are algebraic integers of the field 
${\cal K}$. We recall (see [Wey]) that the ring ${\cal K}_0$ of the 
algebraic integers of the field ${\cal K}$, is a finite-dimensional 
lattice. As a consequence, the image of the representation \eqrefp{N35} 
is a discrete subgroup of the orthogonal group. Since ${\cal G}$ is 
positive definite, the orthogonal group is compact and, hence, $G$ must be 
finite. The theorem is proved.{\hfill QED}
\vskip 0.3 cm

To complete the classification of the monodromy data related to the 
algebraic solutions it remains to classify the objects 
$$
(G,R_1,R_2,R_3),
$$
where $G$ is one of the Coxeter groups $A_3$, $B_3$ and $H_3$ and 
$(R_1,R_2,R_3)$ is a triple of generating reflections considered modulo 
the action \eqrefp{B19} of the braid group. This can be done by a 
straightforward computation of all the orbits of the triples of 
generating reflections. All of them were described and classified 
by Schwartz (see the introduction). We arrive again at the list of 
Theorem 1.6, where, as we already know, the triples \eqrefp{cl1} 
generate the group $W(A_3)$ of the symmetries of the tetrahedron, 
\eqrefp{cl2} generate the group $W(B_3)$ of the symmetries of the cube, 
while \eqrefp{cl3}, \eqrefp{cl4} and \eqrefp{cl5} correspond to three 
inequivalent triples of the generating reflections of the group $W(H_3)$. 
\vskip 0.3 cm







\semiautosez{2}{2. GLOBAL STRUCTURE OF THE SOLUTIONS OF 
PAINLEVE' VI${\bf \mu}$ HAVING CRITICAL BEHAVIOUR OF ALGEBRAIC TYPE}

In the first part of this paper, we found a class of solutions of 
PVI$\mu$ invariant with respect to the analytic continuation. For them, 
the reflection group $G$ acts in the three-dimensional {\it Euclidean}\/ 
space. Recall that the parameter $\mu$ must be real, the coordinates of 
the admissible triples $(x_1,x_2,x_3)$ must be real and satisfy the inequality
$$
-2<x_i<2,\quad i=1,2,3.
$$
In this second part, we prove that this class of solutions coincides with 
the class of the solutions of PVI$\mu$ having critical behaviour of the 
algebraic type
$$
y(x)= \left\{\eqalign{ a_0 x^{l_0} \left(1+{\cal O}(x^\eps)\right),
\qquad\hbox{as}\quad x\rightarrow 0,\cr
1-a_1(1-x)^{l_1} \left(1+{\cal O}((1-x)^\eps)\right),
\qquad\hbox{as}\quad x\rightarrow 1,\cr
 a_\infty x^{1-l_\infty} \left(1+{\cal O}(x^{-\eps})\right),
\qquad\hbox{as}\quad x\rightarrow \infty,\cr
}\right.\autoeqno{BB1}
$$
where $\eps>0$ is small enough, the indices $l_0,\, l_1,\, l_\infty$ 
are real and the coefficients $a_0,\, a_1,\, a_\infty$ are some complex 
numbers. We compute the behaviour of any branch of these solutions near 
the critical points. These results will be used to compute explicitly 
all the algebraic solutions classified in the first part.

First of all, we fix the notations. Let us choose: 
$$
u_1=0,\quad u_2=x,\quad u_3=1.
$$
Then the Fuchsian system \eqrefp{N1} reads
$$
\ddz Y= {\cal A}(z,x)Y=
\left({{\cal A}_1\over z} + {{\cal A}_2\over z-x} + {{\cal A}_3\over z-1}
\right)Y,
$$
and, putting
$$
{\cal A}_1:=A_0,\quad {\cal A}_2:=A_x,
\quad {\cal A}_3:=A_1,\quad {\cal A}_\infty=A_\infty,
$$
we obtain
$$
\ddz Y=
\left({A_0\over z} + {A_1\over z-1} + {A_x\over z-x}
\right)Y.\autoeqno{6}
$$
The branch cuts in $\overline\complessi$ are the same as in 
Section 1.1. We call now the basic loops $\gamma_0,\gamma_x,\gamma_1$. 
They are fixed as before, namely $\gamma_0,\gamma_x,\gamma_1$ play the 
role of the preceding $\gamma_1,\gamma_2,\gamma_3$ (see figure 1).
The Schlesinger equations read:
$$
\eqalign{
\ddx A_0(x)& =-{[A_0,A_x]\over x},\cr
\ddx A_1(x)& =-{[A_1,A_x]\over x-1},\cr
\ddx A_x(x)& = {[A_0,A_x]\over x}+{[A_1,A_x]\over x-1}.\cr}\autoeqno{7}
$$
The correspondent monodromy matrices are
$$
M_0, \quad M_x,\quad  M_1,
$$
which play the role of the preceding $M_1,M_2,M_3$ respectively. We recall 
that they satisfy
$$
M_\infty M_1 M_x M_0=\ID,\qquad
\det (M_i)=1, \quad {\rm Tr}(M_i)=2,\quad\hbox{for}\quad i=0,1,x,
\autoeqno{8}
$$
with
$$
M_\infty=\pmatrix{
\exp (2 i \pi \mu) & 0 \cr 0 &\exp (- 2 i  \pi \mu) \cr}.\autoeqno{9}
$$
With the above choice of $A_\infty$, $A_1$, $A_x$ and $A_0$, satisfying
$$
\det A_i=0, \quad {\rm Tr}(A_i)=0,\quad i=0,1,\infty,\autoeqno{8.5}
$$
the non-singular solution $A(z,x)$ of the Schlesinger equations turns 
out to be related to the solution of PVI$\mu$ in the following way (see 
[JMU]):
$$
[A(y,x)]_{12}=0,\quad\hbox{iff}\quad y(x)\,\, 
\hbox{solves PVI}\mu,\autoeqno{9.5}
$$
where $y$ is not identically equal to $0,1,x$.

We now state the first main theorem of this second part:

\proclaim Theorem 2.1. For any admissible triple $(x_0,x_1,x_\infty)$, 
$x_i\in\reali$, $|x_i|<2$ for $i=0,1,\infty$, there exists a unique 
branch $y(x;x_0,x_1,x_\infty)$ of a solution of PVI$\mu$, with the 
parameter $\mu$ satisfying the equation:
$$
4 \sin^2\pi\mu = x_0^2 + x_1^2 + x_\infty^2-x_0 x_1 x_\infty,\autoeqno{B6}
$$
with the asymptotic behaviour \eqrefp{BB1} near the critical points 
$0,1,\infty$. The indices are given by
$$
l_i={1\over\pi}\arccos(\cos 2\pi r_i)=\left\{\eqalign{
2r_i\quad\hbox{if}\quad 0<r_i\leq {1\over2}\cr
2-2r_i\quad\hbox{if}\quad {1\over2}\leq r_i<1\cr
}\right.\quad i=0,1,\infty,\autoeqno{B5}
$$
with
$$
x_i=-2\cos\pi r_i,\qquad i=0,1,\infty,
$$
and the leading coefficients $a_0,a_1,a_\infty$ are single-valued functions 
of the equivalence class of $x_0,x_1,x_\infty$ and of $\mu$. Namely, the 
coefficient $a_0$, for $x_0\neq 0$ , is given by:
$$
a_0={\exp(-i\pi\phi)\over4(2\mu+l_0-1)^2}
{\Gamma^2(1-l_0)\Gamma^2({1+l_0\over2})\Gamma({1+l_0\over2}+\mu)
\Gamma({1+l_0\over2}-\mu)\over\Gamma^2(l_0)\Gamma^2({1-l_0\over2})
\Gamma({1-l_0\over2}+\mu)\Gamma({1-l_0\over2}-\mu)}\autoeqno{A1}
$$
where 
$$
\exp(i\pi\phi)={x_0^2x_1^2-2 x_1^2-2 x_0 x_1 x_\infty+2x_\infty^2+i x_1
sign(x_0)\sqrt{4-x_0^2}
(2x_\infty-x_0 x_1)\over 2(x_1^2- x_0 x_1 x_\infty+x_\infty^2)}
\autoeqno{A2}
$$
and for $x_0=0$
$$
a_0={x_\infty^2\over x_1^2+x_\infty^2}.\autoeqno{A3}
$$  
The coefficient $a_1$ is given by the same formula with the substitution 
$x_0\leftrightarrow x_1$, $l_0\mapsto l_1$; $a_\infty$ is given by the same 
formula too, after the substitution 
$(x_0,x_1,x_\infty)\mapsto(x_\infty,-x_1,x_0-x_1x_\infty)$ and 
$l_0\mapsto l_\infty$.
Conversely any solution of the PVI$\mu$ equation, with a real value of $\mu$, 
having critical behaviour of algebraic type, can be obtained by the above 
construction.

\noindent{\bf Remark 2.1.}\quad The relation \eqrefp{B6} determines $\mu$ 
up to the transformations 
$$
\mu\mapsto\pm \mu + n,\quad n\in\interi.
$$ 
According to the results of Section 1.2, such an ambiguity can be absorbed
by the action of a symmetry on PVI$\mu$. Recall that these symmetries
preserve the class of the algebraic solutions.
\vskip 0.3 cm

Theorem 2.1 will be proved in Section 2.4.

\vskip 0.3 cm
\noindent{\bf 2.1. Local theory of the solutions of PVI$\mu$ having critical 
behaviour of algebraic type}
\vskip 0.3 cm

\noindent{\bf 2.1.1. Local asymptotic behaviour around 0.}\quad
In this section we characterize the local asymptotic behaviour of the 
solutions of PVI$\mu$ near the singular point $x=0$. First of all let us 
characterize the type of asymptotic behaviour that can be related to the 
algebraic solutions.

\proclaim Lemma 2.1. Let $y(x)$ be an algebraic solution of 
PVI$\mu$. Then the first term of its Puiseux series is
$$
y(x)\sim\,a_0 x^{1-\sigma_0}\quad\hbox{as}\quad x \rightarrow 0.
\autoeqno{15.5}
$$
for some constant $a_0\neq0$ and the rational number $\sigma_0$ must 
satisfy $0\leq\sigma_0<1$, with $a_0\neq1$ if $\sigma_0=0$.

\noindent Proof.\quad If $y(x)$ is an algebraic function, 
then it admits an expansion in Puiseux series around $0$
$$
y(x) = \sum_{k=k_0}^\infty a_k x^{k\over n}, \qquad k_0\in\interi, 
\qquad a_{k_0}\neq 0,
$$
where $n$ is some natural number. As a consequence, for $k_0\neq 0$, we have 
the following relation between the orders of the first and second derivative 
of $y$:
$$
{\cal O}(x^2 y'')={\cal O}(x y')={\cal O}(y)=
{\cal O}\left(x^{k_0\over n}\right).
\autoeqno{14.3}
$$
We now reduce to the common denominator the PVI$\mu$ equation and collect 
together all the terms of the same order in the numerator ${\cal N}$, using 
the rule \eqrefp{14.3}. The numerator is
$$
\eqalign{
{\cal N}= & 2 {y'}^2 x^4-{y'}^2 x^3 - {y'}^2 x^5 + 2 {y'}^2 x^2 y - 
2 {y'} x^3 y - 2 {y'}^2 x^3 y +  2 {y''} x^3 y + 2 {y'} x^4 y - 
2 {y'}^2 x^4 y - \cr
&- 4 {y''} x^4 y +2 {y'}^2 x^5 y  +2 {y''} x^5 y +   x y^2 - 2 {y'} x y^2 - 
2 x^2 y^2 + 4 \mu  x^2 y^2 - 4 \mu^2 x^2 y^2 + \cr
& + 6 {y'} x^2 y^2 - 3 {y'}^2 x^2 y^2 - 2 {y''} x^2 y^2 - 2 {y'} x^3 y^2  + 
6 {y'}^2 x^3 y^2 + 2 {y''} x^3 y^2 - 2 {y'} x^4 y^2 -\cr
& -3 {y'}^2 x^4 y^2 +2 {y''} x^4 y^2 - 2 {y''} x^5 y^2 - 8 \mu  x y^3 + 
8 \mu^2 x y^3 + 2 {y'} x y^3 + 4 x^2 y^3-\cr
&- 8 \mu  x^2 y^3 +8 \mu^2 x^2 y^3 - 6 {y'} x^2 y^3 + 2 {y''} x^2 y^3 + 
4 {y'} x^3 y^3 - 4 {y''} x^3 y^3  + 2 {y''} x^4 y^3-\cr
& -y^4 + 4 \mu  y^4 -4 \mu^2 y^4 - 3 x y^4 + 16 \mu  x y^4 - 16 \mu^2 x y^4 
- 2 x^2 y^4 + 4 \mu  x^2 y^4 - \cr
&-4 \mu^2 x^2 y^4 + 2 y^5  - 8 \mu  y^5 + 8 \mu^2 y^5 + 2 x y^5 -   
8 \mu  x y^5 + 8 \mu^2 x y^5 - y^6 + 4 \mu  y^6 - 4 \mu^2 y^6.\cr}
$$
The first term of the Puiseux series must be chosen in order to kill the 
lowest term in the numerator of the PVI$\mu$ equation.
If $k_0<0$, the lowest term is
$$
-y^6+4 \mu y^6 - 4\mu^2 y^6
$$
which, for $2\mu\not\in\interi$ cannot be zero for any choice of $a_0\neq0$. 
Then $k_0$ cannot be negative. If $n \geq k_0>0$, the lowest order term is
$$
2x^2 {y'}^2 y - 2 x y' y^2 -2 x^2 y'' y^2, 
$$
which is zero for any $y = a_{k_0} x^{k_0\over n}$.
For $k_0>n$, the lowest order term is
$$
-x^3 {y'}^2 + 2 x^3 y'' y + x y^2,
$$
which cannot be zero. Furthermore, for $k_0=0$, the lowest order term in the 
numerator ${\cal N}$ is
$$
-a_0^4(a_0-1)^2(2\mu-1)^2
$$
and, due to the assumptions $2\mu\not\in\interi$ and $a_0\neq0$, the only 
possible value of $a_0$ is $1$. Substituting $y=1+a_1 x^{k_1\over n}$, we 
obtain that the lowest order term in the numerator ${\cal N}$ is
$$
x^{2k_1\over n} a_1^2 ({k_1/n}+1-2\mu) ({k_1/n}-1+2\mu) 
$$
that is zero, for generic values of $\mu$, only if $a_1$ is $0$. If 
$\mu={1\over 2}\pm{k_1\over 2 n}$, we can again repeat the procedure. 
The numerator will be
$$
\eqalign{
\hat{\cal N}=&
-4 \hat\mu^2 \hat y^2 - 16 \hat\mu^2 \hat y^3 - 24 \hat\mu^2 \hat y^4 - 
16 \hat\mu^2 \hat y^5 - 4 \hat\mu^2 \hat y^6 +  \hat y^2 x + 
8 \hat\mu^2 \hat y^2 x + 2 \hat y^3 x + 24 \hat\mu^2 \hat y^3 x + \cr
&+\hat y^4 x + 24 \hat\mu^2 \hat y^4 x + 8 \hat\mu^2 \hat y^5 x + 
2 \hat y{y'} x + 4 \hat y^2{y'} x + 2 \hat y^3{y'} x - \hat y^2 x^2 - 
4 \hat\mu^2 \hat y^2 x^2 - 2 \hat y^3 x^2 -\cr
&- 8 \hat\mu^2 \hat y^3 x^2 - \hat y^4 x^2 - 4 \hat\mu^2 \hat y^4 x^2 - 
6 \hat y{y'} x^2 - 12 \hat y^2{y'} x^2 - 6 \hat y^3{y'} x^2 -{y'}^2 x^2 - 
4 \hat y{y'}^2 x^2 - \cr
&-3 \hat y^2{y'}^2 x^2 +2 \hat y {y''} x^2 + 4 \hat y^2 {y''} x^2 + 
2 \hat y^3 {y''} x^2 + 6 \hat y{y'} x^3 + 10 \hat y^2{y'} x^3 + 
4 \hat y^3{y'} x^3 +  \cr
& +3{y'}^2 x^3 + 10 \hat y{y'}^2 x^3 +  6 \hat y^2{y'}^2 x^3 -
6 \hat y {y''} x^3 - 10 \hat y^2 {y''} x^3 - 4 \hat y^3 {y''} x^3 - 
2 \hat y{y'} x^4 -  \cr
&-2 \hat y^2{y'} x^4 - 3{y'}^2 x^4  -3 \hat y^2{y'}^2 x^4 + 6 \hat y {y''} x^4 
+  8 \hat y^2 {y''} x^4 + 2 \hat y^3 {y''} x^4 +{y'}^2 x^5 + \cr
&+2 \hat y{y'}^2 x^5 - 2 \hat y {y''} x^5 - 2 \hat y^2 {y''} x^5
\cr}
$$
where $\hat\mu=\pm{k_1\over2n}$ and $\hat y= y-1$. Substituting 
$\hat y = a_1 x^{k_1\over n}$, the lowest order term in the numerator 
$\hat{\cal N}$ is automatically zero. Now, we want to eliminate the next 
lowest order term. Observe thet, now
$$
{\cal O}(x^2 y'')={\cal O}(x y')={\cal O}(\hat y)=
{\cal O}\left(x^{k_1\over n}\right).
$$
For the sake of definitess, suppose ${1\over 2}<\mu\leq 1$, i.e. 
$\hat\mu= {k_1\over2n}<{1\over 2}$ (the case $\hat\mu=-{k_1\over2n}$ is 
analogous). The next lowest order terms in the numerator $\hat{\cal N}$ are 
$$
-16\hat\mu^2 \hat y^3 + 4 x\hat y^2 y' - 4 x^2\hat y{y'}^2 + 
4x^2\hat y^2 y'' + x \hat y^2 + 8 x \hat\mu^2 \hat y^2 - 
6 x^2\hat y y' + 3 x^3{y'}^2 - 6 x^3\hat y y''
$$
To eliminate them, we substitute $y=1+a_1 x^{k_1\over n}+a_1 y^{k_2\over n}$, 
for some $k_2>k_1$. The above terms 
give
$$
-4 a_1^3 ({k_1\over n})^2 x^{3 {k_1\over n}}+{\cal O}(x^{1+2{k_1\over n}})
$$
that is zero if and only of $a_1=0$. So we obtain the forbidden solution 
$y(x)\equiv 1$. 
So, $k_0$ can not be zero, and $y(x)$ satisfies \eqrefp{15.5} with 
$0<l={k_0\over n}\leq 1$, namely $0\leq \sigma_0<1$. {\hfill QED}
\vskip 0.3 cm

In the above lemma we have seen the expected asymptotic behaviour of the 
algebraic solutions. We now state the main result of this section, which 
is more general, namely it holds also for non algebraic solutions.

\proclaim Theorem 2.2. For any pair of values $(a_0, \sigma_0)$, 
$0\leq\sigma_0<1$, there exists a unique branch of the solution of 
PVI$\mu$, for a fixed $\mu$, with the asymptotic behaviour 
$$
y(x)= a_0 x^{1-\sigma_0}(1+ x^\eps f(x))\quad\hbox{as}\quad 
x \rightarrow 0,\autoeqno{15}
$$
for some $\eps>0$ and $f(x)$ smooth function such that 
$\lim_{x\rightarrow0} f(x)={\rm const}$. 

In order that $x^{1-\sigma_0}$ is well defined, we have to make some cut 
in the complex plane. From now on, we cut along the line 
${\rm arg} x=\varphi$ for some $\varphi$.

\vskip 0.2 cm
\noindent{\bf Remark 2.2.}\quad Theorem 2.2 can be proved also for complex 
values of the index $\sigma_0$, provided that $0\leq\RE\sigma_0<1$. For 
algebraic solutions the  index $\sigma_0$ must be a rational number. Because 
of this, we consider only real indices.

\vskip 0.2 cm
\noindent{\bf 2.1.2. Proof of the existence.}\quad 
First of all we state the existence of solutions of the Schlesinger 
equations with a particular asymptotic behaviour. The following result 
will play an important role also in Section 2.3.

\proclaim Lemma 2.2 (Sato-Miwa-Jimbo). Given three constant matrices 
$A_i^0$, $i=0,1,x$ with zero eigenvalues such that $\Lambda=A_0^0+A_x^0$ 
has eigenvalues $\pm {\sigma\over2}$, $0\leq\sigma<1$, and 
$A_1^0=-\Lambda-A_\infty$, in any sector of $\overline\complessi$ 
containing none of the branch cuts, and sufficiently close to $0$, 
there exists a solution of the Schlesinger equations that satisfy
$$
|A_1(x)-A_1^0|\leq K |x|^{1-\sigma'} \quad  
\left|x^{-\Lambda}(A_1(x)-A_1^0) x^\Lambda\right| \leq K |x|^{1-\sigma'}
\autoeqno{17}
$$
$$ 
\left|x^{-\Lambda}A_0(x)x^\Lambda-A_0^0 \right|\leq K |x|^{1-\sigma'}
\quad
\left|x^{-\Lambda}A_x(x)x^\Lambda-A_x^0 \right|\leq K |x|^{1-\sigma'},
\autoeqno{18}
$$
where $K$ is some positive constant and $1>\sigma'>\sigma$.

We want to show that it is possible to choose $A_{0,1,x}$ and $\Lambda$ 
such that the corresponding solution $y(x)$ of the Painlev\'e VI equation 
obtained via \eqrefp{9.5} has the asymptotic behaviour \eqrefp{15}.
Let us consider an arbitrary constant matrix $\Lambda$ with eigenvalues 
$\pm {\sigma\over2}$; let $T$ be the diagonalizing matrix of $\Lambda$, 
namely
$$
\Lambda= T \pmatrix{{\sigma\over2}&0\cr 0&-{\sigma\over2}\cr} T^{-1}.
$$
Now, we choose $A_1^0=-A_\infty-\Lambda$ and $A_{0,x}^0$ such that 
$A_0^0+A_x^0=\Lambda$, namely
$$
A_0^0={1\over2}\Lambda + F, \qquad A_x^0={1\over2}\Lambda - F
$$
for some constant matrix $F$. Then:
$$
T^{-1} A_0^0 T = \pmatrix{{\sigma\over4}&0\cr 0&-{\sigma\over4}\cr} + E, 
\qquad 
T^{-1} A_x^0 T= \pmatrix{{\sigma\over4}&0\cr 0&-{\sigma\over4}\cr} - E, 
$$
where we can choose 
$E=\pmatrix{0&{b\sigma\over 4}\cr {-\sigma\over4 b}&0\cr}$, for some 
non-zero constant $b$. With this choice of $E$, $A_0$ and $A_x$ have zero 
eigenvalues. Using Lemma 2.2, we obtain that, as $x\rightarrow 0$:
$$
A_{0,x} \rightarrow T \,
{\pmatrix{x^{\sigma\over2}&0\cr 0&x^{-\sigma\over2}\cr}}
\left[ 
\pmatrix{{\sigma\over4}&0\cr 0&-{\sigma\over4}\cr} \pm 
\pmatrix{0&{b\sigma x^{\sigma}\over 4}\cr 
{-\sigma x^{-\sigma}\over4 b}&0\cr}
\right] {\pmatrix{x^{-\sigma\over2}&0\cr 0&x^{\sigma\over2}\cr}} 
T^{-1},
$$ 
and
$$
A_1\rightarrow-A_\infty-\Lambda.
$$
Substituting such asymptotic behaviors in the relation \eqrefp{9.5}, 
taking $T_{12},\, T_{11}\neq 0$ we obtain: 
$$
y(x)\sim -{T_{12} x^{1-\sigma}\over 4 b T_{11}};\autoeqno{p4}
$$
we are now free to choose the arbitrary constants 
$b,\,T_{11},\,T_{12},\,\sigma$ in such a way that 
$-{T_{12} \over 4 b T_{11}}=a_0$, $\sigma=\sigma_0$, for any fixed 
$a_0$ and $\sigma_0$.

\noindent{\bf Remark 2.3.} Other existence results for 
$\sigma\in\complessi\backslash\{]-\infty,0]\cup[1,+\infty[\}$ can be found 
in [IKSY] and [S1], [S2], [S3]. For indices with $\RE\sigma\not\in[0,1]$, the
asymptotics obtained in these papers are valid in more complicated 
domains near $0$.

\vskip 0.3 cm
\noindent{\bf 2.1.3. Proof of the uniqueness.}\quad 
Now we prove that the solution 
$y(x)$, $x\in {\cal B}(0,r)$, of Painlev\'e VI equation such that it 
satisfies \eqrefp{15} for some given constants $a_0$ and $\sigma_0\in [0,1)$, 
is uniquely determined by $a_0$ and $\sigma_0$. Here ${\cal B}(0,r)=
\left\{x\big | |x|\leq r,\, {\rm arg} x\neq \varphi,\, x\neq 0\right\}$.

\noindent The proof is based on the fact that Painlev\'e VI is equivalent to 
the following reduced Schlesinger equations \eqrefp{N13}:
$$
\dot q = {(q-1) q + 2 p (q-1) q (q-x)\over(x-1)x}, 
$$
$$ 
\dot p ={- p^2 (x - 2q - 2x q+3q^2)- p (2q-1)-(1-\mu)\mu\over(x-1)x},
$$
where:
$$
q = y, \quad p =  {x(x-1)\dot y-y(y-1)\over 2(y-x)y(y-1)}, \autoeqno{U2}
$$ 
and the dot means the derivative $\ddx$
We shall prove the local uniqueness of the solutions of the Hamiltonian 
system with the following asymptotic behaviour
$$
q(x)\sim  a x^l+ x^{l+\eps} \,f(x)\qquad
p(x)  \sim {l-1\over 2 a}{1\over x^l} + {x^\eps g(x)\over x^l }\autoeqno{U3}
$$ 
where $l=1-\sigma_0$, $a=a_0$, $\eps>0$ and $f(x)$ and $g(x)$ are some 
smooth functions in  ${\cal B}(0,r)$ which tend to zero as $x\rightarrow0$.

This is equivalent to show the theorem. In fact, from the uniqueness 
of $q$ it follows trivially the uniqueness of $y$. The following lemma 
holds true:

\proclaim Lemma 2.3. The estimates \eqrefp{U3} on the asymptotic behaviour 
of $(q(x),p(x)) $ are a consequence of \eqrefp{15}.

\noindent Proof.\quad
Since $q=y$, the assertion on $y$ is obvious due to the hypothesis 
\eqrefp{15}. Concerning $p$, we use its definition
$$
p =  {x(x-1)\dot y-y(y-1)\over 2(y-x)y(y-1)}
$$
and by a straightforward computation we show \eqrefp{U3} for p.{\hfill QED}

\vskip 0.3 cm

We now distinguish two cases: $0<l<1$, and $l=1$.
Let us consider the former case; it is convenient to introduce the new 
variables $(\tilde q,\tilde p)$       
$$
\tilde q = {y\over x^l}\quad \tilde p = x^l p; 
$$
which have a similar asymptotic behaviour
$$
\tilde q(x)= a+ x^\eps f(x)\qquad
\tilde p(x)= {l-1\over 2 a} + x^\eps g(x),
\autoeqno{U4}
$$
and the equations of the motion become
$$
\eqalign{
\dot{\tilde q} &= f_q(\tilde p, \tilde q, x, x^l),\cr
\dot{\tilde p} &= f_p(\tilde p, \tilde q, x, x^l),\cr}\autoeqno{U5}
$$
with
$$
f_q={-\tilde q(l-1-2\tilde p \tilde q)\over x}-
{\tilde q(1+2\tilde p \tilde q)\over x^{1-l}}+
{\tilde q(1+4\tilde p \tilde q)-
2 x^{-l}\tilde p \tilde q-x^l\tilde q^2(1+2\tilde p \tilde q)\over 1-x},
$$
and
$$
f_p={\tilde p(l-1-2\tilde p \tilde q)\over x}+{\mu-\mu^2+ 
2\tilde p \tilde q+ 3\tilde p^2 \tilde q^2\over x^{1-l}}+
{{\tilde p^2\over x^l}-\tilde p(1+4\tilde p \tilde q)+
x^l(\mu-\mu^2+ 2\tilde p \tilde q+ 3\tilde p^2 \tilde q^2) \over 1-x}.
$$
We want to prove the uniqueness of the solution $(\tilde q,\tilde p)$ 
of \eqrefp{U5}, satisfying \eqrefp{U4} for $x\in{\cal B}(0,r)$, 
in the ball $||\tilde p- {l-1\over 2 a}||,\, ||\tilde q-a|| \,\leq C_r$, 
for a constant $C_r$ vanishing when the radius $r\rightarrow 0$. Here  
$||f||= \sup_{{\cal B}(0,r)}|f(x)|$.
Let us suppose that there are two solutions $(\tilde q_1,\tilde p_1)$ and 
$(\tilde q_2,\tilde p_2)$ of the system \eqrefp{U5}, satisfying \eqrefp{U4}. 
Then, if we define ${\rm X}=\left(\matrix{\tilde q_1-\tilde q_2\cr 
\tilde p_1-\tilde p_2\cr}\right)$, we obtain, as a consequence of 
\eqrefp{U4}, that the following limits exist
$$
\lim_{|x|\rightarrow 0,\,{\rm arg}(x)=\vartheta} 
{{\rm X}^{(i)}(x)\over |x|^\eps}=0,\quad i=1,\,2,\autoeqno{U6}
$$
for some $0<\eps$, ${\rm X}^{(i)}$ being the i-th component of ${\rm X}$. 
Moreover, ${\rm X}$ satisfies the following
$$
{\rm X}'=\left(\matrix{
{[1-l+2\tilde p_1(\tilde q_1+\tilde q_2)](\tilde q_1-\tilde q_2)+
2\tilde q_2^2(\tilde p_1-\tilde p_2)\over x} + {\Delta 
{\rm Q}_1\over x^{1-l}} + {\Delta {\rm Q}_2\over x^l}+ 
\Delta {\rm Q}_3\cr
{[l-1-2\tilde q_2(\tilde p_1+\tilde p_2)](\tilde p_1-\tilde p_2)-
2\tilde p_1^2(\tilde q_1-\tilde q_2)\over x} + 
{\Delta {\rm P}_1\over x^{1-l}} + {\Delta {\rm P}_2\over x^l}+ 
\Delta {\rm P}_3\cr}\right),
$$
where 
$$
\Delta {\rm Q}_i
=Q_i(\tilde q_1,\tilde p_1,x)-Q_i(\tilde q_2,\tilde p_2,x),
\quad\hbox{and}\quad
\Delta {\rm P}_i=P_i(\tilde q_1,\tilde p_1,x)-P_i(\tilde q_2,\tilde p_2,x),
$$
$Q_1=\tilde q(1+2\tilde p\tilde q)$, $Q_2=-2\tilde p\tilde q$, 
$Q_3= {\tilde q(1+4\tilde p \tilde q)-
2 x^{1-l}\tilde p \tilde q-x^l\tilde q^2(1+2\tilde p \tilde q)\over 1-x}$, 
$P_1=\mu-\mu^2+2\tilde p\tilde q+3 \tilde p^2\tilde q^2$, $P_2=\tilde p^2$,
$P_3={p^2 x^{1-l}-\tilde p(1+4\tilde p \tilde q)+
x^l(\mu-\mu^2+ 2\tilde p \tilde q+ 3\tilde p^2 \tilde q^2) \over 1-x}$.

We want to prove that, under the hypothesis \eqrefp{U6}, 
${\rm X} \equiv 0$ (this is equivalent to prove our theorem). 
Performing the constant linear transformation 
${\rm X} = {\rm T}\, {\rm Z}$, where 
$$
{\rm T} = \left(\matrix{
1&0\cr {1-l\over 2 a^2}&{1\over 2 a^2}\cr}\right),
$$
we obtain 
$$
{\rm Z}'=\left(\matrix{
{[1-l+2\tilde p_1(\tilde q_1+\tilde q_2)](\tilde q_1-\tilde q_2)+
2\tilde q_2^2(\tilde p_1-\tilde p_2)\over x} 
+ {\Delta {\rm Q}_1\over x^{1-l}} + {\Delta {\rm Q}_2\over x^l}+ 
\Delta {\rm Q}_3\cr
{2 a^2 G(\tilde p_1,\tilde p_2, \tilde q_1, \tilde q_2)\over x}
+ {2 a^2\Delta {\rm P}_1 +(l-1) \Delta {\rm Q}_1\over x^{1-l}} + 
{2 a^2\Delta {\rm P}_2+(l-1) \Delta{\rm Q}_2\over x^l}+
2 a^2 \Delta {\rm P}_3+(l-1)\Delta{\rm Q}_3\cr}\right)
\autoeqno{U7}
$$
where 
$$
\eqalign{
G(\tilde p_1,\tilde p_2, \tilde q_1, \tilde q_2)=&
[l-1-2\tilde q_2(\tilde p_1+\tilde p_2)]
(\tilde p_1-\tilde p_2)-2\tilde p_1^2(\tilde q_1-\tilde q_2)\}+\cr
&+(l-1)\{[1-l+2\tilde p_1(\tilde q_1+\tilde q_2)](\tilde q_1-\tilde q_2)+
2\tilde q_2^2(\tilde p_1-\tilde p_2),\cr}
$$
and, from \eqrefp{U6}:
$$
\lim_{|x|\rightarrow 0,\,{\rm arg}(x)=\vartheta} 
{{\rm Z}^{(i)}(x)\over |x|^\eps}=0,\quad i=1,\,2.
\autoeqno{U8}
$$
In order to prove that ${\rm Z}\equiv 0$, we fix any direction in the 
complex plane ${\rm arg}(x)=\vartheta$ for some fixed $\vartheta$, and we 
consider the real variable $t=|x|$. Then we define:
$$
{\rm V}^{(i)}(t):=|Z^{(i)}(x)|.
$$
We want to prove that the assumption ${\rm V}^{(i)}(t_0)\neq 0$ for some 
$t_0>0$ leads to a cotradiction. To this aim we prove a differential 
inequality for the right derivative $D_{+}{\rm V}^{(i)}$ of $V^{(i)}(t)$. 
Since $D_{+}{\rm V}^{(i)}\leq |{\rm Z}^{(i)'}|$, to obtain such a differential 
inequality it is enough to estimate from above the modulus of the components 
of the right-hand-side of \eqrefp{U7}. To this aim we notice that all 
the polynomials $Q_i$, $P_i$ have the form:
$$
Q_i=\sum_{k,n=0}^3 a_{k,n}^{i}\tilde p^n\tilde q^k,\qquad 
P_i=\sum_{k,n=0}^3  b_{k,n}^{i}\tilde p^n\tilde q^k,
$$
with $a_{k,n}(x)$, $b_{k,n}(x)$ regular functions $x\in {\cal B}(0,r)$. 
As a consequence, we obtain, in the ball $||\tilde p- {l-1\over 2 a}||,\, 
||\tilde q-a|| \,\leq C_r$, the estimates:
$$
|\Delta {\rm Q}_i|, \quad |\Delta {\rm P}_i|
\leq c_1^i|{\rm Z}^{(1)}| + c_2^i|{\rm Z}^{(2)}|\autoeqno{U8,7}
$$
for some positive constants $c_1^i,c_2^i$. In fact
$$
\eqalign{
|\Delta {\rm Q}_i|&=
\big|\sum_{k, n} a_{k,n}[\tilde q_1^k(\tilde p_1^n-\tilde p_2^n) +
\tilde p_2^n(\tilde q_1^k-\tilde q_2^k)]\big|\leq 
\left\{
\sum_{k=0,1,2} C_r^{(1)^k}(||a_{k,1}||+2 C_r^{(2)} ||a_{k,2}||)\right\} \cr
& \cdot |\tilde p_1-\tilde p_2|
+\left\{  \sum_{n=0,1,2,3} C_r^{(2)^n}(||a_{1,n}||+2 C_r^{(1)} ||a_{2,n}||+ 
3  C_r^{(1)^2} ||a_{3,n}||)\right\}\cdot
|\tilde q_1-\tilde q_2|,\cr}
$$
where $C^{(1)}_r =C_r + 2 |a|$ and  $C^{(2)}_r =C_r + {1-l\over |a|}$. 
We obtain \eqrefp{U8,7} observing that $|\tilde q_1-\tilde q_2|$, 
$|\tilde p_1-\tilde p_2|$ are related to $|{\rm Z}^{(1)}|$, 
$|{\rm Z}^{(2)}|$ by the constant linear transformation $T$.

\vskip 0.2 cm\noindent
For the terms of order ${\cal O}({1\over x})$ in \eqrefp{U7} we have:
$$
\eqalign{
{|[1-l+2\tilde p_1(\tilde q_1+\tilde q_2)](\tilde q_1-\tilde q_2)+
2\tilde q_2^2(\tilde p_1-\tilde p_2)|\over |x|}\leq&
{ |-(1-l)(\tilde q_1-\tilde q_2)+2 a^2(\tilde p_1-\tilde p_2)|\over |x|}+\cr
&+{ C_r^{(3)}|\tilde q_1-\tilde q_2|\over |x^{1-\eps}|}  +
{ C_r^{(4)} |\tilde p_1-\tilde p_2|\over|x^{1-\eps}|}, \cr}
\autoeqno{Ui}
$$
and 
$$
\eqalign{
&{1\over x} |2 a^2\{[l-1-2\tilde q_2(\tilde p_1+\tilde p_2)]
(\tilde p_1-\tilde p_2)-2\tilde p_1^2(\tilde q_1-\tilde q_2)\}+
(l-1)\{[1-l+2\tilde p_1(\tilde q_1+\tilde q_2)]\cdot\cr
&(\tilde q_1-\tilde q_2)+2\tilde q_2^2(\tilde p_1-\tilde p_2)\}|
\leq { C_r^{(5)}\over|x^{1-\eps}|} |\tilde q_1-\tilde q_2| + 
{ C_r^{(6)}\over|x^{1-\eps}|} |\tilde p_1-\tilde p_2|, \cr} 
\autoeqno{Uii}
$$
for some positive constants $ C_r^{(3)},..., C_r^{(6)}$. Let us prove 
\eqrefp{Ui}:
$$
\eqalign{
&{|[1-l+2\tilde p_1(\tilde q_1+\tilde q_2)](\tilde q_1-\tilde q_2)+
2\tilde q_2^2(\tilde p_1-\tilde p_2)|\over |x|}\leq{ |-(1-l)
(\tilde q_1-\tilde q_2)+2 a^2(\tilde p_1-\tilde p_2)|\over |x|}+\cr
&   +{|2 a g_1(x)+{l-1\over a}(f_1(x)+f_2(x))+
x^\eps g_1(x)(f_1(x)+f_2(x))|\over |x^{1-\eps}|} |\tilde q_1-\tilde q_2|+\cr
&+{|2 f^2_2(x) x^\eps + 
4 a f_2(x)|\over |x^{1-\eps}|}|\tilde p_1-\tilde p_2|\leq 
{|-(1-l)(\tilde q_1-\tilde q_2)+2 a^2(\tilde p_1-\tilde p_2)|\over |x|}+\cr
&+{ C_r^{(3)}\over |x^{1-\eps}|} |\tilde q_1-\tilde q_2|+
 {C_r^{(4)}\over |x^{1-\eps}|} |\tilde p_1-\tilde p_2|,\cr}
$$
for some positive constants $ C_r^{(3)}$ and $ C_r^{(4)}$. The proof of 
\eqrefp{Uii} is analogous.
From the estimates \eqrefp{U8,7}, \eqrefp{Ui}, \eqrefp{Uii}, we obtain:
$$
\left(
\eqalign{
&|{\rm Z}^{(1)'}|\cr &|{\rm Z}^{(2)'}|\cr}\right)
\leq \left( {1\over|x|}\left(\matrix{0&1\cr 0&0\cr}\right) + 
{{\rm A}_1\over|x^{1-l}|}+{{\rm A}_2\over|x^l|}+
{{\rm A}_3\over|x^{1-\eps}|}+{\rm A}_4\right)\left( \eqalign{
&|{\rm Z}^{(1)}|\cr &|{\rm Z}^{(2)}|\cr}\right)\autoeqno{U8.22}
$$
for some constant matrices ${\rm A}_1$, ${\rm A}_2$, ${\rm A}_3$ and 
${\rm A}_4$ (Here we mean $\leq$ component by component). Finally, 
choosing $\tilde l=\max\{1-\eps,1-l,l\}$, we obtain from \eqrefp{U8.22}:
$$
\left(
\eqalign{
&D_{+,t}{\rm V}^{(1)}\cr & D_{+,t}{\rm V}^{(2)} \cr}\right)
\leq \left( {1\over t}\left(\matrix{0&1\cr 0&0\cr}\right) + 
{\tilde{\rm A}\over t^{\tilde l}}+{\rm A}_4\right)\left( \eqalign{
&{\rm V}^{(1)}\cr &{\rm V}^{(2)}\cr}\right),\autoeqno{U8.2}
$$
where $\tilde{\rm A}={\rm A}_1+{\rm A}_2+{\rm A}_3$ and $D_{+,t}$ is the right 
derivative w.r.t. $t$.

We perform the following change of variable $t^{1-\tilde l}=z$. The 
differential inequality for ${\rm V}$ in the new variable $z$ is
$$
D_{+,z}{\rm V}\leq\left(
{1\over z}\left(\matrix{0&1\cr 0&0\cr}\right){1\over 1- \tilde l} +
{\rm A}(z)\right){\rm V},\quad\hbox{with}\quad {\rm A}(z)=
{\tilde{\rm A}\over1-\tilde l}+{{\rm A}_4\over1-\tilde l}\, 
z^{{\tilde l\over 1-\tilde l}},
$$
where $D_{+,z}$ is the right derivative w.r.t. $z$. To show that ${\rm Z}=0$ 
we use the following:

\proclaim Comparison Theorem. Let us consider the following systems of 
$n$ first order ODEs in the real variable $z\in(0,a]$, for some $a>0$:
$$
D_{+}{\rm V}^{(i)} \leq F^{(i)}(z,{\rm V}), \quad{\rm V}^{(i)}(x_0)=
{\rm V}_0^{(i)},\quad i=1,..,n \autoeqno{U9} 
$$
$$
{{\rm d}{\rm U}^{(i)}\over{\rm d}z} = F^{(i)}(z,{\rm U}), 
\quad {\rm U}^{(i)}(x_0)={\rm U}_0^{(i)},\quad i=1,..,n 
\autoeqno{U10}
$$
where $F^{(i)}(z,{\rm U})$ are continuous functions in $z\in(0,a]$, 
$||{\rm U}-{\rm U}_0||<b$, non-decreasing in ${\rm U}^{(i)}$. If 
${\rm V}_0^{(i)}\geq {\rm U}_0^{(i)}$, for $i=1,..n$, then 
${\rm V}^{(i)}(z)\geq {\rm U}^{(i)}(z)$, for every $0<z\leq z_0$, $i=1,..n$.

\noindent For the proof see [Lak].\vskip 0.3 cm 

We now apply Comparison Theorem to show that the assumption $Z(t_0)\neq 0$ 
for some $t_0>0$ leads to a contradiction. Observe that by definition 
$\tilde l\geq 1$, then ${\rm V}$ satisfies \eqrefp{U9} with 
${\rm V}^{(i)}_0>0$ and $F$ linear in ${\rm V}$ given by:
$$
F(z,{\rm V})=\left(
{1\over z}\left(\matrix{0&1\cr 0&0\cr}\right){1\over 1- \tilde l} +
{\tilde{\rm A}\over1-\tilde l}+{{\rm A}_4\over1-\tilde l}\, z
\right){\rm V}.\autoeqno{N11.5}
$$
By Comparison Theorem, for any solution of \eqrefp{U10} with 
${\rm U}^{(i)}(z_0)={\rm V}_0^{(i)}$, and $F(z,{\rm U})$ of the form 
\eqrefp{N11.5}, we have 
${\rm V}^{(i)}(z)\geq {\rm U}^{(i)}(z)$, for every $0<z\leq z_0$, $i=1,2$. 
Moreover by standard arguments it is possible to take ${\rm U}$ in such 
a way that ${\rm U}^{(i)}(z)\geq0$ and to continue the functions ${\rm U}$,  
${\rm V}$ to $z=0$ preserving the relation:
$$
0\leq{\rm U}^{(i)}(z)\leq {\rm V}^{(i)}(z).
$$
Thus, by \eqrefp{U8} we obtain that ${\rm U}$ must satisfy
$$
\lim_{z\rightarrow 0} {{\rm U}^{(i)}(z)\over z^{\eps\over(1-\tilde l)}}=0,
\quad i=1,\,2.
\autoeqno{U11}
$$
Now, we use the following lemma:

\proclaim Lemma 2.4. The only solution ${\rm U}$ of \eqrefp{U10} with
$F(z,{\rm U})$ given by \eqrefp{N11.5} satisfying \eqrefp{U11} is 
${\rm U}\equiv 0$,

\noindent Proof. \quad Any non-zero solution of \eqrefp{U10} with
$F(z,{\rm U})$ of the form \eqrefp{N11.5} is given by
$$
{\rm U}(z)=T(z) z^{\pmatrix{0&{1\over 1-\tilde l}\cr 0&0\cr}}
\autoeqno{U12.5}
$$
where $T(z)$ is a holomorphic matrix function,
$T(z)= \pmatrix{1&0\cr 0&1\cr}+{\cal O}(z)$.
Now it is obvious that \eqrefp{U12.5} does not satisfy \eqrefp{U11}. Thus
${\rm U}\equiv 0$, as we wanted to prove. {\hfill QED}
\vskip 0.2 cm

Using the above lemma, we obtain ${\rm U}_0^{(i)}={\rm V}_0^{(i)}=0$, that
contradicts the assumption ${\rm V}_0^{(i)}=\neq0$. This concludes the proof 
of the uniqueness in the case $0<l<1$.

Let us briefly explain how to prove the uniqueness in the case $l=1$. Since 
the procedure is essentially the same as before, we shall skip the details. 
First of all we introduce the new variables $(\tilde q, \tilde p)$:
$$
\tilde q(x) = {q(x)\over x} \sim a +x^\eps  f(x)\qquad
\tilde p(x) = p(x) - \mu(1-\mu) \sim x^\eps g(x)
$$
which satisfy the equations of the motion:
$$
\eqalign{
\dot{\tilde q} & = Q_1(\tilde q,\tilde p) + 
{1\over x-1}  Q_2(\tilde q,\tilde p)\cr 
\dot{\tilde p} & = - {\tilde p\over x}+  P_1(\tilde q,\tilde p) + 
{1\over x-1}  P_2(\tilde q,\tilde p) \cr}
$$
where $Q_1(\tilde q,\tilde p) = 2(\mu-\mu^2+\tilde p)(\tilde q-1)\tilde q^2$, 
$ Q_2(\tilde q,\tilde p)= \tilde q(\tilde q-1)[1+(2\mu(1-\mu)+2\tilde p)
(\tilde q-1)]$, $P_1(\tilde q,\tilde p)=
(\mu-\mu^2+\tilde p)^2(2-3\tilde q)\tilde q-\mu(1-\mu)$, and 
$P_2(\tilde q,\tilde p)=\tilde p + (\mu-\mu^2+\tilde p)^2
(4\tilde q-3\tilde q^2-1)-2 (\mu-\mu^2+\tilde p)\tilde q$.
Then, if we define ${\rm X}$ as before we obtain
$$
{\rm X}'=\left(\matrix{
\Delta Q_1+{\Delta Q_2\over x-1}\cr
-{\tilde p_1-\tilde p_2\over x} + \Delta P_1+{\Delta P_2\over x-1}\cr}\right)
$$
that gives rise to the differential inequality:
$$
|{\rm X}'|\leq \left(
\left(\matrix{0&0\cr 0&1\cr}\right){1\over|x|}+{\rm A}_1 
+{{\rm A}_2\over|x-1|}\right)|{\rm X}|
$$
for some constant matrices ${\rm A}_1$ and ${\rm A}_2$. Obviously 
${\rm X}$ satisfies \eqrefp{U6} with any $0<\eps<1$. Again we apply 
Comparison Theorem to ${\rm V}:= \left(\eqalign{&|{\rm X}^{(1)}|\cr 
&|{\rm X}^{(2)}|\cr}\right)$
along any fixed direction on the complex plane. We take $x$ such 
that ${\rm arg}(x)=\vartheta$ for some fixed $\vartheta$ and define 
$t=|x|$. ${\rm V}$ satisfies \eqrefp{U9} with:
$$
F(t,{\rm V})=\left(
{1\over t}\left(\matrix{0&0\cr 0&1\cr}\right)+{\rm A}_1 +{{\rm A}_2\over t-1}
\right) {\rm V}.
$$
If ${\rm V}^{(i)}_0>0$ then, thanks to Comparison Theorem, it is 
possible to take a solution ${\rm U}$ of \eqrefp{U10}, with 
${\rm U}^{(i)}_0={\rm V}^{(i)}_0 $, $i=1,2$, such that
$$
0\leq{\rm U}^{(i)}(t)\leq {\rm V}^{(i)}(t),
$$
thus ${\rm U}^{(i)}$ satisfies \eqrefp{U11}. The general solution of 
\eqrefp{U10} is
$$
{\rm U} = {\rm U}_0\left(\left(\matrix{ 1&0\cr 0&t\cr}\right) + 
{\cal O}(t^2)\right)
$$ 
that satisfies \eqrefp{U11} iff $ {\rm U}_0=0$, namely $ {\rm U}\equiv 0$ 
that is absurd. This concludes the proof of the uniqueness.
{\hfill QED}

\vfill\eject
\noindent{\bf 2.1.4. Asymptotic behaviour of the solutions of the Schlesinger 
equations.}\quad
An important corollary to Theorem 2.2 is the following:
   
\proclaim Theorem 2.3. The solutions of the Schlesinger equations 
$A_{0,1,x}(x)$ corresponding to the solution of Painlev\'e VI equation 
with asymptotic behaviour \eqrefp{15} must satisfy the relations \eqrefp{17} 
and \eqrefp{18}.

\noindent Proof. Let us consider the solution $y(x)$ of 
Painlev\'e VI equation with asymptotic behaviour \eqrefp{15} and let us 
suppose that the corresponding solution of the Schlesinger 
equations $A_{0,1,x}(x)$ does not satisfy the relations \eqrefp{17} and 
\eqrefp{18}. As shown in the lemma 2.2, for any constant matrices 
$A_{0,1,x}^0$, $\Lambda$ such that $\Lambda=A_0^0+A_x^0$ has eigenvalues 
$\pm {\sigma\over2}$, $\sigma\in[0,1[$, and $A_1^0=-\Lambda-A_\infty$, there 
exists a solution $\hat A_{0,1,x}(x)$ of the Schlesinger equations that 
satisfy the relations \eqrefp{17} and \eqrefp{18}. Now, as shown in 
Section 2.2, we can choose $A_{0,1,x}^0$ in order that the corresponding 
solution $\hat y(x)$ of Painlev\'e VI equation has exactly the asymptotic 
behaviour \eqrefp{15}. Due to the uniqueness proved in Theorem 2.2, we 
have that $y(x)=\hat y(x)$, namely $A_{0,1,x}=\hat A_{0,1,x}$ up to 
conjugation by a constant diagonal matrix. This contradiction proves the 
theorem. {\hfill QED}

\vskip 0.3 cm
\noindent{\bf 2.1.5. Asymptotic behaviour of the PVI$\mu$ solution near $1$ 
and $\infty$.}\quad
We now state the analogues of Theorem 2.2 for the local asymptotic 
behaviour of the solutions of (PVI) near the singular points $x=1,\infty$:

\proclaim Theorem 2.2'. For any pair of values $(a_1, \sigma_1)$, 
$\sigma_1\in[0,1[$, there exists a unique branch of the solution of (PVI) 
with the asymptotic behaviour 
$$
y(x)\sim\,1-a_1 (1-x)^{1-\sigma_1}
\left(1+{\cal O}\left((1-x)^\eps\right)\right)
\quad\hbox{as}\quad x\rightarrow1,
\autoeqno{15i}
$$
for some $\eps>0$.

\noindent The proof of this theorem is analogous to the proof of theorem 1, 
namely one can state the analogous of the lemma 2.2 replacing 
$x\mapsto 1-x$, and then choose suitably $\Lambda$, $A_{0,1,x}^0$. The 
uniqueness is proved in the same way as the case $x\mapsto 0$.  

\proclaim Theorem 2.2''. For any pair of values $(a_\infty, \sigma_\infty)$, 
$\sigma_\infty\in[0,1[$, there exists a unique branch of the solution of 
(PVI) with the asymptotic behaviour 
$$
y(x)\sim\,a_\infty x^{\sigma_\infty}
\left(1+{\cal O}\left((x^{-\eps}\right)\right)
\quad\hbox{as}\quad x\rightarrow\infty,
\autoeqno{15ii}
$$
for some $\eps>0$.

\noindent
The proof of uniqueness is analogous to the one of Theorem 2.2. The proof of 
existence follows the same strategy as the one of Theorem 2.2, but with a 
different formulation of the lemma 2.2:

\proclaim Lemma 2.2'. Given some constant matrices $A_i^0$, $i=0,1,x$ with 
zero eigenvalues such that $\Lambda=A_0^0+A_t^0$ has eigenvalues 
$\pm {\sigma\over2}$, $0\leq \sigma<1$, in any sector of 
$\overline\complessi$ containing none of the branch cuts, and sufficiently 
close to $\infty$, there exists a solution of the Schlesinger equations
satisfying:
$$
|x^{A_\infty} A_x(x)x^{-A_\infty} -A_1^0|\leq K |x|^{\sigma'-1} \quad  
\left|x^{\Lambda}\left (x^{A_\infty} A_x(x)x^{-A_\infty}-A_1^0\right) 
x^{-\Lambda}\right|\leq K |x|^{\sigma'-1}
\autoeqno{17i}
$$
$$ 
\left|x^{\Lambda}x^{A_\infty}A_{0,1}(x)x^{-A_\infty}x{-^\Lambda}-
A_x^0 \right| \leq K |x|^{\sigma'-1},
\autoeqno{18i}
$$
where $K$ is some positive constant and $1>\sigma'>\sigma$.

\noindent Proof. Let us consider the Schlesinger equations 
\eqrefp{7} and perform the change of variable $x={1\over\hat x}$. 
Moreover we put:
$$
A_i(x):= x^{-A_\infty} \hat A_i(x) x^{A_\infty};
$$
Then we can apply Lemma 2.2 to the system:
$$
\eqalign{
\ddtx\hat A_0(\hat x)& =-{[\hat A_0,\hat A_1]\over\hat x},\cr
\ddtx\hat A_x(\hat x)& =-{[\hat A_x,\hat A_1]\over\hat x-1},\cr
\ddtx\hat A_1(\hat x)& ={[\hat A_0,\hat A_1]\over\hat x}+
{[\hat A_x,\hat A_1]\over\hat x-1},\cr}
$$
and obtain the estimates \eqrefp{17i} and \eqrefp{18i}. {\hfill QED}


\vskip 0.3 cm 
\noindent{\bf 2.2. The local asymptotic behaviour and the monodromy group 
of the Fuchsian system}
\vskip 0.3 cm

In this section we relate the local asymptotic behaviour of the solution 
$y(x)$ of PVI$\mu$ to the monodromy data of the associated Fuchsian system 
\eqrefp{6}. We essentially follow the same strategy of [Jim], even if we 
have to introduce some more tricks due to the 
fact that our matrices $A^0_{0,1,x}$ have eigenvalues all equal to zero. 
The main result of this section is the following:

\proclaim Theorem 2.4. For the solution $y(x)$ of PVI$\mu$, such that 
$y(x)\sim a_0 x^{1-\sigma_0}(1+{\cal O}(x^\eps))$, $0<\sigma_0<1$, the 
monodromy matrices of the Fuchsian system \eqrefp{6} have the form
$$
\eqalign {M_1=&{-i\over\sin\pi\thi}\cdot \cr
&\cdot\pmatrix{
\cos\pi\sigma_0-e^{-i\pi\thi}&
-2e^{-i\pi\thi}\sin{\pi(\thi+\sigma_0)\over2}\sin{\pi(\thi-\sigma_0)\over2}\cr
2e^{i\pi\thi}\sin{\pi(\thi+\sigma_0)\over2}\sin{\pi(\thi-\sigma_0)\over2}&
-\cos\pi\sigma_0+e^{i\pi\thi}\cr
}\cr}\autoeqno{m1}
$$
$$
C M_x C^{-1}={-i\over\sin\pi\sigma_0}\pmatrix{
e^{i\pi\sigma_0}-1&2 s e^{i\pi\sigma_0}\sin^2{\pi\sigma_0\over2}\cr
-{2\over s}e^{-i\pi\sigma_0} \sin^2{\pi\sigma_0\over2}&
1-e^{-i\pi\sigma_0}\cr}\autoeqno{mt}
$$
$$
C M_0 C^{-1}={-i\over\sin\pi\sigma_0}\pmatrix{
e^{i\pi\sigma_0}-1&-2s\sin^2{\pi\sigma_0\over2}\cr
{2\over s} \sin^2{\pi\sigma_0\over2}&
1-e^{-i\pi\sigma_0}\cr}\autoeqno{m0}
$$
where $\thi=2\mu$ and:
$$
{s\over r} ={1\over 4 a_0} {2\mu+\sigma_0\over 2\mu-\sigma_0}
{\Gamma^2(1+\sigma_0)\Gamma^2(1-{\sigma_0\over2})
\Gamma(1+\mu-{\sigma_0\over2})\Gamma(1-\mu-{\sigma_0\over2})\over
\Gamma^2(1-\sigma_0)\Gamma^2(1+{\sigma_0\over2})
\Gamma(1+\mu+{\sigma_0\over2})\Gamma(1-\mu+{\sigma_0\over2})}
\autoeqno{eqs}
$$
with an arbitrary complex number $r\neq 0$ and the matrix $C$ is:
$$
C=\pmatrix{
\sin{\pi(\thi-\sigma_0)\over2}& r\sin{\pi(\thi+\sigma_0)\over2}\cr
{1\over r}\sin{\pi(\thi+\sigma_0)\over2}& \sin{\pi(\thi-\sigma_0)\over2}\cr}.
\autoeqno{mc}
$$
In the case where $\sigma_0=0$ the monodromy matrices of the Fuchsian 
system \eqrefp{6} have the form
$$
M_1={1\over\cos{\pi\thi\over2}}
\pmatrix{e^{-i\pi{\thi\over2}}&i\pi e^{-i\pi{\thi\over2}}\cr
-{i\over\pi}\sin^2{\pi\thi\over2}e^{i\pi{\thi\over2}}&
e^{i\pi{\thi\over2}}\cr},\autoeqno{m10}
$$
$$
M_0=\pmatrix{
1-i s\tan{\pi\thi\over2}&
-i s\pi\exp(i\pi{\thi\over2})\sec{\pi\thi\over2}\cr
{i\over\pi} s\sin^2{\pi\thi\over2}\exp(-i\pi{\thi\over2})
\sec{\pi\thi\over2}&
1+i s\tan{\pi\thi\over2}\cr},\autoeqno{m00}
$$
$$
M_x=\pmatrix{
1-i(1-s)\tan{\pi\thi\over2}&
-i(1- s)\pi\exp(i\pi{\thi\over2})\sec{\pi\thi\over2}\cr
{i\over\pi}(1- s)\sin^2{\pi\thi\over2}\exp(-i\pi{\thi\over2})
\sec{\pi\thi\over2}&
1+i(1-s)\tan{\pi\thi\over2}\cr},\autoeqno{mx0}
$$
where $s=a_0$.

The main idea to prove this theorem is that, due to Theorem 2.3, the 
solutions of the Schlesinger equations corresponding to the PVI$\mu$ 
solution with the asymptotic behaviour \eqrefp{15} must satisfy the 
relations \eqrefp{17} and \eqrefp{18}. Using these relations, we obtain the 
monodromy matrices of the Fuchsian system \eqrefp{6} via the ones of two 
simpler systems, given in the following two lemmas (see [SMJ] and [Jim]):

\proclaim Lemma 2.5. Under the hypotheses \eqrefp{17}, \eqrefp{18}, the 
limit of the fundamental solution of the system \eqrefp{6}, normalized at 
infinity,
$\lim_{x\rightarrow0}\, Y_\infty(z,x)=\hat Y(z)$, exists, for 
$z\in\overline\complessi\backslash\{B_0\cup B_x\cup B_1\cup B_\infty\}$, 
$B_0,\, B_x,\, B_1$ and $B_\infty$ being balls around $0,\,x,\,1$ and 
$\infty$ respectively. This limit $\hat Y$ satisfies the differential equation:
$$
\ddz\hat Y=\left({A_1^0\over z-1}+{\Lambda\over z}\right)\hat Y;
\eqno{(\hat\Sigma)}
$$
and it has the following behaviour near the singularities of ${(\hat\Sigma)}$
$$
\eqalign{
\hat Y(z)=&
\left(1+{\cal O}({1\over z})\right)z^{-A_\infty}
\qquad z\rightarrow\infty \cr
&=\left(1+{\cal O}(z)\right)z^{\Lambda}\hat C_0 
\qquad z\rightarrow 0 \cr
& =\hat G_1\left(1+{\cal O}(z-1)\right)(z-1)^{J_1}\hat C_1
\qquad z\rightarrow1\cr}
\autoeqno{22.5}
$$
where $J_1$ is the Jordan normal forms of $A_1^0$, 
$\hat G_1 J_1\hat G^{-1}_1 = A_1^0$, 
$A_\infty=\pmatrix{\mu & 0\cr 0 & -\mu\cr}$. Here $\hat C_0$, $\hat C_1$ 
are the connection matrices of the system $(\hat\Sigma)$.

\noindent{\bf Remark 2.4.}\quad Observe that the matrix $\hat C_0$ is uniquely 
determined by the conditions \eqrefp{22.5}.

\proclaim Lemma 2.6. Under the hypotheses \eqrefp{17}, \eqrefp{18}, the 
limit of the fundamental solution of the system \eqrefp{6}, normalized 
around $\infty$,
$\lim_{x\rightarrow0}x^{-\Lambda}Y(xz,x)=\tilde Y(z)\hat C_0$ exists 
for $z\in\overline\complessi\backslash\{B_0\cup B_x\cup B_1\cup B_\infty\}$. 
It satisfies the system
$$
\ddz\tilde Y=\left({A_x^0\over z-1}+{A_0^0\over z}\right)\tilde Y;
\eqno{(\tilde\Sigma)}
$$
and it has the following behaviour near the singularities of
${(\tilde\Sigma)}$ 
$$
\eqalignno{
\tilde Y(z)=
&\left(1+{\cal O}({1\over z})\right)z^{\Lambda}& z\rightarrow\infty\cr
& =\tilde G_0 \left(1+{\cal O}(z)\right)z^{J_0}\tilde C_0 &z\rightarrow 0\cr
& =\tilde G_1\left(1+{\cal O}(z-1)\right)(z-1)^{J_1}\tilde C_1
&z\rightarrow 1\cr}
$$
where $J_{0,1}$ are the Jordan normal forms of $A_{0,x}^0$, $\tilde G_{0,1}$ 
are such that $\tilde G_{0,1} J_{0,1}\tilde G^{-1}_{0,1} = A_{0,x}^0$. We
denote $\tilde C_{0,1}$ the connection matrices of the system $(\tilde\Sigma)$.


As we have seen above, the matrices of the two systems have the following form:
$$
A_0^0={1\over2}\Lambda + F, \qquad A_x^0={1\over2}\Lambda - F,
\quad A_1^0=-A_\infty-\Lambda,
$$
for some constant matrix $F$, and for $\Lambda$ and $T$ such that
$$
\Lambda= T \pmatrix{{\sigma\over2}&0\cr 0&-{\sigma\over2}\cr} T^{-1}.
\autoeqno{ss1}
$$
Using the relations \eqrefp{8.5}, we have that
$$
F=T \pmatrix{0&{b\sigma\over 4}\cr {-\sigma\over4 b}&0\cr} T^{-1},
\autoeqno{ss2}
$$
for some parameter $b$. As a consequence the systems $(\hat\Sigma)$ and 
$(\tilde\Sigma)$ are determined, up to diagonal conjugation, by the four 
entries of the matrix $T$ and by $b$.  

Now, we explain how to compute the monodromy matrices of the original 
system \eqrefp{6} knowing the ones of the systems $(\hat\Sigma)$ and 
$(\tilde\Sigma)$. Later we will show how to compute the matrices 
$A_{0,x,1}^0$ and the monodromy matrices of $(\hat\Sigma)$ and 
$(\tilde\Sigma)$.

\proclaim Lemma 2.7. Let $\hat M_0,\hat M_1,\hat M_\infty=M_\infty$ be the 
monodromy matrices of the system $(\hat\Sigma)$ with respect to the 
fundamental matrix $\hat Y$ and the basis 
$\hat\gamma_0= \gamma_0\gamma_x,\,\gamma_1$ in 
$\pi_1\left(\overline\complessi\backslash\{0,1,\infty\}\right)$. Let 
$\tilde M_0,\tilde M_1,\tilde M_\infty=\exp(-2\pi i \Lambda)$ be the 
monodromy matrices of the system $(\tilde\Sigma)$ with respect to the 
fundamental matrix $\tilde Y$ and the basis 
$\tilde\gamma_0,\tilde\gamma_1=\gamma_x$. Then the monodromy matrices 
of the original system \eqrefp{6} are given by the formulae:
$$
M_0= \hat C_0^{-1}\tilde M_0 \hat C_0,\quad
M_x= \hat C_0^{-1}\tilde M_1 \hat C_0,\quad
M_1=\hat M_1,
\autoeqno{23}
$$
where $ \hat C_0$ is defined by \eqrefp{22.5}.

\noindent Proof. By the definition of $\hat Y$, the system $(\hat\Sigma)$ 
is obtained by merging of the singularities $0$ and $x$ of the system 
\eqrefp{6}. We can choose the loop $\hat\gamma_0$ to be homotopic to 
$\gamma_0\gamma_x$, with $\hat\gamma_0$ not crossing a ball 
the $B_0$ (see figure 7).
\midinsert
\centerline{\psfig{file=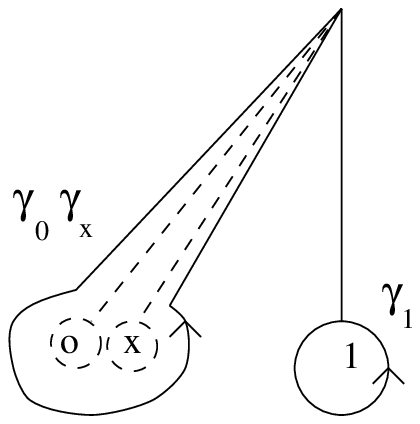,height=4cm}}
\vskip 0.6 cm
\centerline{{\bf Fig.7.} The paths $\gamma_x$ and $\gamma_0$ merge together 
as $x\to 0$. The} 
\centerline{homotopy class of $\gamma_0\gamma_x$ remains unchanged.}
\endinsert
%
%
%
As a consequence we obtain a relation between the monodromy matrices of 
the system \eqrefp{6} and the ones of the system $(\hat\Sigma)$
$$
\eqalign{\hat M_\infty = &\, M_\infty, \cr
\hat M_1 = &\, M_1, \cr
\hat M_0 = &\, M_x M_0. \cr}
$$ 
Similarly, by the definition of $\tilde Y$ the system $(\tilde\Sigma)$ is 
obtained by the merging (see figure 8) of the singularities $z'={1\over x}$ 
and $z'=\infty$ of the system for $Y'(z')$:
$$
{{\rm d}\over{\rm d} z'}  Y' = 
\left({A_0\over z'} + {A_1\over z'-{1\over x}}+{A_x\over z'-1} \right)Y'.
$$ 
\midinsert
\centerline{\psfig{file=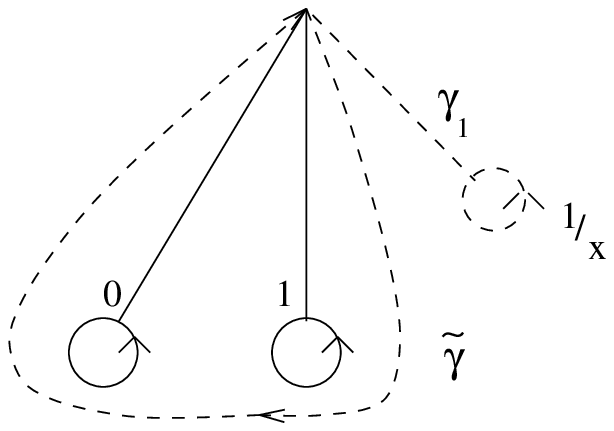,height=4cm}}
\vskip 0.7 cm
\centerline{{\bf Fig.8.} The paths $\gamma_1$ and $\gamma_\infty$ merge 
together as $x\to 0$. The homotopy }
\centerline{class of $\tilde\gamma\equiv\gamma_1\gamma_\infty$ coincides 
with the one of $(\gamma_0\gamma_x)^{-1}$}
\endinsert
%
%
So, in the basis $\hat Y$, the monodromy matrices of $(\tilde\Sigma)$ have 
the following form:
$$
\eqalign{\tilde M_\infty = &\, \hat C_0^{-1} M_\infty M_1 \hat C_0\cr
\tilde M_1 = &\, \hat C_0^{-1} M_x \hat C_0, \cr
\tilde M_0 = &\, \hat C_0^{-1} M_0 \hat C_0. \cr}
$$ 
The lemma is proved. {\hfill QED}

\vskip 0.3 cm
Now we want to compute the monodromy matrices $\hat M_i$ and $\tilde M_i$ 
and the connection matrix $\hat C_0$. To this aim we have to solve the 
systems $(\hat\Sigma)$ and $(\tilde\Sigma)$, namely we have to determine 
$T$ and $b$.
For $\sigma_0\neq 0$, this can be done introducing a suitable gauge 
transformation of $\hat Y$ and $\tilde Y$ such that the 
systems  $(\hat\Sigma)$ and $(\tilde\Sigma)$ are equivalent to a Gauss 
equation. The case $\sigma_0=0$ will be treated later.

\vskip 0.3 cm
\noindent{\bf 2.2.1. Reduction to the Gauss equation.}\quad
First of all let us notice that both the systems $(\hat\Sigma)$ and 
$(\tilde\Sigma)$ have similar form. We want to reduce them, via a 
suitable gauge transformation and a appropriate choice of the parameters 
$\alpha,\beta,\gamma$, to systems of the form:
$$
\ddz Y(z,\alpha,\beta,\gamma) = 
\left({B_0\over z}+{B_1\over z-1}\right) Y(z,\alpha,\beta,\gamma)
\autoeqno{24}
$$
where $B_0$, $B_1$ are some constant matrices with eigenvalues 
$1-\gamma,\,0$ and $\gamma-\alpha-\beta-1,\,0$ respectively and 
$B_0+B_1=
-\pmatrix{\alpha&0\cr 0&\beta\cr}$.

\proclaim Lemma 2.8. For $\alpha\neq\beta$, the system \eqrefp{24} is 
uniquely determined, up to a diagonal conjugation
$$
B_0\to T^{-1} B_0 T,\quad B_1\to T^{-1} B_1 T,\quad\hbox{with}\quad 
T=\pmatrix{1&0\cr 0&r\cr},\, r\neq0.\autoeqno{224}
$$
The entries $b_{ij}^0$ and $b_{ij}^1$ of the matrices $B_0$ and $B_1$ 
respectively, are given by the formulae
$$
b^0_{11}={\alpha(\gamma-1-\beta)\over\beta-\alpha},\quad
b^0_{22}={-\beta(\gamma-1-\alpha)\over\beta-\alpha},\quad
b^1_{11}={-\alpha(\gamma-1-\alpha)\over\beta-\alpha},\autoeqno{param1}
$$
$$
b^1_{22}={\beta(\gamma-1-\beta)\over\beta-\alpha},\quad
b^0_{12}b^0_{21}=b^1_{12}b^1_{21}=
{-\alpha\beta(\gamma-1-\beta)(\gamma-1-\alpha)\over(\beta-\alpha)^2}.
\autoeqno{param2}
$$

The system \eqrefp{24} can be solved using the Gauss hypergeometric function.
So, we can compute its connection matrices via the Kummer relations 
(see [Luke]) of the hypergeometric functions.

\proclaim Lemma 2.9. The solutions of \eqrefp{24} have the form 
$Y=\left(\eqalign{&y_1 \cr &y_2\cr} \right)$, with $y_1$ being an 
arbitrary solution of the following Gauss equation:
$$
z(1-z) y_1'' +[ c-(a+b+1)z] y_1' -a b y_1=0\autoeqno{29}
$$
where $a=\alpha$, $b=\beta+1$, $c=\gamma$
and $y_2$ given by:
$$
y_2(z)= r^{-1} {\beta-\alpha\over\beta(\gamma-\beta-1)} 
\left\{z(z-1)\ddz y_1(z)+\left[\alpha z +
{\alpha(\gamma-\beta-1)\over\beta-\alpha}\right]y_1(z)\right\}\autoeqno{30}
$$
where $r=-{(\beta-\alpha) b^0_{12}\over\beta(\gamma-\beta-1)}$. 

\noindent Proof. After the gauge transformation:
$$
Y(z,\alpha,\beta,\gamma)= 
z^{b^0_{11}}(1-z)^{b^1_{11}} U(z,\alpha,\beta,\gamma),
$$
one obtains from \eqrefp{24} the following Riemann equation for $u_1$
$$
u_1''+\left[{1 + b^0_{11} - b^0_{22}\over z} + 
{1 + b^1_{11} - b^1_{22}\over z-1} \right]u_1'
- {b^0_{11}b^0_{22}\over z^2(z-1)^2} u_1=0.
$$
Now $u_1$ is related with the solution $y_G$ of the Gauss equation 
\eqrefp{29}, with $a=-b^0_{11}-b^1_{11}$, $b=1-b^0_{22}-b^1_{22}$, 
$c=1-b^0_{11}-b^0_{22}$, via the relation 
$u_1=z^{-b^0_{11}}(1-z)^{-b^1_{11}} y_G$. As a consequence, thanks to 
\eqrefp{param1}, \eqrefp{param2}, we obtain that $y_1=y_G$ and $a=\alpha+1$,
$b=\beta$, $c=\gamma$.
$y_2$ it is given by:
$$
\left({b^0_{12}\over z}+{b^1_{12}\over z-1}\right) y_2 
= y_1' -\left({b^0_{11}\over z}+{b^1_{11}\over z-1}\right) y_1
$$
that gives the equation \eqrefp{30}. {\hfill QED}
\vskip 0.2 cm

To reduce the systems $(\hat\Sigma)$ and $(\tilde\Sigma)$ to the system 
\eqrefp{24} we need to diagonalize the matrices $A_1^0+\Lambda=-A_\infty$ 
and $\Lambda$ respectively and to perform a suitable gauge transform.
We need to introduce some notations. Denote 
$C_{0,1}^{\alpha,\beta,\gamma}$ the connection matrices of the 
system \eqrefp{24}. The matrices $J_{0,1}$ are the Jordan normal forms 
of $B_{0,1}$ and the matrices $G_{0,1}^{\alpha,\beta,\gamma}$ are such 
that $G_{0,1}^{\alpha,\beta,\gamma}J_{0,1}
\left(G_{0,1}^{\alpha,\beta,\gamma}\right)^{-1}=B_{0,1}$. Then for the 
asymptotic behaviour of an appropriate fundamental matrix 
$Y(z,\alpha,\beta,\gamma)$ of the system \eqrefp{24} we have
$$
\eqalignno{
Y(z,\alpha,\beta,\gamma)
& =\left(1+{\cal O}({1\over z})\right)
z^{-\pmatrix{\alpha&0\cr 0&\beta\cr} }& z\rightarrow\infty\cr
& =G_0^{\alpha,\beta,\gamma} \left(1+{\cal O}(z)\right)z^{J_0} 
C_0^{\alpha,\beta,\gamma} &z\rightarrow 0\cr
& = G_1^{\alpha,\beta,\gamma} \left(1+{\cal O}(z-1)\right)
(z-1)^{J_1} C_1^{\alpha,\beta,\gamma} &z\rightarrow 1.\cr}
$$
Some further remarks on the notations: from now on all the quantities 
with the hat are referred to the system $(\hat\Sigma)$ and all the 
quantities with the tilde to the system $(\tilde\Sigma)$. When we don't 
put any hat or tilde, the formulae are true for both the systems. In other words, they hold true for the generic system \eqrefp{24}; substituting all 
the quantities with the correspondent hat or tilde ones, the formulae hold 
true for the systems $(\hat\Sigma)$ or $(\tilde\Sigma)$ respectively.

We now choose the values of $\alpha,\beta,\gamma$ in relation with the 
eigenvalues of the matrices of the systems $(\hat\Sigma)$ and 
$(\tilde\Sigma)$. Namely, for $(\hat\Sigma)$ we take 
$$
\hat\alpha ={\thi-\sigma_0\over2},\qquad
\hat\beta =-{\thi+\sigma_0\over2},\qquad
\hat\gamma  =1-\sigma_0, \autoeqno{25}
$$
and for $(\tilde\Sigma)$ we take: 
$$
\tilde\alpha =-{\sigma_0\over2},\qquad
\tilde\beta ={\sigma_0\over2},\qquad
\tilde\gamma  =1.\autoeqno{26}
$$
With this choice of the values of $\alpha,\beta,\gamma$, one has:
$$
\hat J_0= \pmatrix { 1-\gamma&0\cr 0&0\cr},\quad 
\tilde J_{0}=\pmatrix{0&1\cr 0 &0\cr},\quad 
\hat J_1=\tilde J_1=\pmatrix{0&1\cr 0 &0\cr}.
$$
Now we can reduce the systems $(\hat\Sigma)$ and $(\tilde\Sigma)$ to the 
system \eqrefp{24} via the following gauge transformations:
$$
\hat Y = z^{\hat\alpha+\hat\beta\over 2} 
Y(z,\hat\alpha,\hat\beta,\hat\gamma),\qquad
\tilde Y=G_0^{\hat\alpha,\hat\beta,\hat\gamma} \,
Y(z,\tilde\alpha,\tilde\beta,\tilde\gamma),\autoeqno{22}
$$
where $G_0^{\hat\alpha,\hat\beta,\hat\gamma}$ is such that 
$$
\Lambda-{\hat\alpha+\hat\beta\over2}\Id= 
G_0^{\hat\alpha,\hat\beta,\hat\gamma}
\pmatrix{-\tilde\alpha&0\cr 0&-\tilde\beta\cr }
(G_0^{\hat\alpha,\hat\beta,\hat\gamma})^{-1}.
$$
As a consequence the connection matrices of \eqrefp{24} are related to  
the ones of $(\hat\Sigma)$ and $(\tilde\Sigma)$ by the following formulae:
$$
\hat G_1 = G_1^{\hat\alpha,\hat\beta,\hat\gamma},\quad
\hat C_1= C_1^{\hat\alpha,\hat\beta,\hat\gamma},\quad
\hat C_0 = G_0^{\hat\alpha,\hat\beta,\hat\gamma} 
C_0^{\hat\alpha,\hat\beta,\hat\gamma}
\autoeqno{27}
$$
$$
\tilde G_{0,1} = G_0^{\hat\alpha,\hat\beta,\hat\gamma} 
G_{0,1}^{\tilde\alpha,\tilde\beta,\tilde\gamma},\qquad
\tilde C_{0,1} = C_{0,1}^{\tilde\alpha,\tilde\beta,\tilde\gamma}
(G_0^{\hat\alpha,\hat\beta,\hat\gamma})^{-1}.\autoeqno{28}
$$

\vskip 0.3 cm
\noindent{\bf 2.2.2. Local behaviour of the solution of \eqrefp{24}.}\quad
The solutions of \eqrefp{29} around the singular points $0,1,\infty$ are 
known and one can compute $y_2$ by \eqrefp{30}. In this way one obtains the 
local behaviour of the fundamental solution $Y$ for $z\rightarrow0,1,\infty$, 
and one can compute the connection matrices by the Kummer relations (which 
are the connection formulae for the hypergeometric equation). The 
difference w.r.t. the situation of [Jim] is that in our case the Gauss 
equation is degenerate, namely:
$$
\hat c-\hat a-\hat b=0,\qquad\tilde c-\tilde a-\tilde b= 0,\qquad \tilde c=1.
$$
So, we have to consider the logarithmic solutions of the Gauss equation 
around $z=1$ for both the systems $(\hat\Sigma)$ and $(\tilde\Sigma)$, and 
around $z=0$ for $(\tilde\Sigma)$; moreover, we shall use the extension of 
the Kummer relations to this logarithmic case (see [Nor]).

In what follows we denote $F(a,b,c,z)$ the hypergeometric function and with 
$g(a,b,z)$ its logarithmic counterpart for $c=1$, namely:
$$
F(a,b,c,z)=\sum_{k=0}^\infty {(a)_k(b)_k\over k!(c)_k} \, z^k,
$$
$$
g(a,b,z)=\sum_{k=0}^\infty {(a)_k(b)_k\over k!} \, z^k
[\ln z + \psi(a+k)+\psi(b+k)-2\psi(k+1)],
$$
with the branchcut $|{\rm arg}(z)|<\pi$ (see figure 9). Here $\psi$ is 
the logarithmic derivative of the gamma function, and the expressions of 
the parameters $a,b,c$ via $\alpha,\beta,\gamma$ are given in the lemma 2.9. 

\midinsert
\centerline{\psfig{file=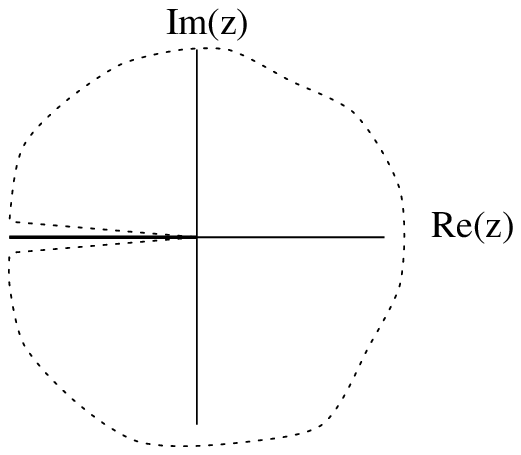,height=3cm}}
\vskip 0.5 cm
\centerline{{\bf Fig.9.} The branch cut $|{\rm arg}(z)|<\pi$.}
\endinsert


\noindent{\sl Fundamental solution near $\infty$}\quad
Since $a-b\neq 0$, the solutions of \eqrefp{29} around $\infty$ are not 
logarithmic. We obtain
$$
Y_\infty=\left(\eqalign{
&z^{-\alpha} F(\alpha,-\beta,\alpha-\beta,{1\over z})\qquad\quad
{-\alpha \beta\, z^{-\beta-1} r\over(\beta-\alpha)(\beta-\alpha+1)}
F(\beta+1,1-\alpha,\beta-\alpha+2,{1\over z})\cr
&{-\alpha \beta\, z^{-\alpha-1}\over r (\beta-\alpha)(\beta-\alpha-1)}
F(\alpha+1,1-\beta,\alpha-\beta+2,{1\over z})\qquad\quad
z^{-\beta} F(\beta,\alpha,\beta-\alpha,{1\over z})\cr
}\right)
$$
$$
Y_\infty\sim 
\left(1+{\cal O}({1\over z})\right)\, 
z^{-\pmatrix{\alpha&0\cr 0&\beta\cr} },\qquad z\rightarrow\infty.
$$
The monodromy around $\infty$ is 
$\pmatrix{\exp(2\pi i\alpha)&0\cr 0&\exp(2\pi i\beta)}$.

\noindent{\sl Fundamental solution near $1$.}\quad
Since $c-a-b=0$, the solutions are logarithmic:
$$
Y_1=\pmatrix{
F(\alpha,\beta+1,1,{1- z}) &  r\,g(\alpha,\beta+1,1,{1- z}) \cr
{1\over r} F(\alpha+1,\beta,1,{1- z}) & g(\alpha+1,\beta,1,{1- z}) \cr}.
$$
For $z\to1$
$$
Y_1\sim G_1^{\alpha,\beta}\,\,(1- z)^{\pmatrix{0&1\cr 0&0\cr}},
$$
with
$$
G_1^{\alpha,\beta}=\pmatrix{
1&r[\psi(\alpha)+\psi(1+\beta)-2\psi(1)]\cr
{1\over r}&\psi(1+\alpha)+\psi(\beta)-2\psi(1)
\cr}.
$$
The monodromy around $1$ is $\pmatrix{1&2 i\pi r\cr 0&1\cr}$.

\noindent{\sl Fundamental solution near $0$.}\quad
We have to distinguish the case $(\hat\Sigma)$, where the solutions of 
\eqrefp{29} around $0$ are not logarithmic, and the case $(\tilde\Sigma)$, 
where $c=1$ and the solutions are logarithmic. 

For $(\hat\Sigma)$ one has
$$
\hat Y_0=\pmatrix{
-{\hat\alpha\over\hat\beta-\hat\alpha} z^{-\hat\alpha-\hat\beta} 
F(-\hat\beta,1-\hat\alpha,1-\hat\alpha-\hat\beta,z)&
\hat r{\hat\beta\over\hat\beta-\hat\alpha} 
F(\hat\alpha,\hat\beta+1,\hat\alpha+\hat\beta+1,z)\cr
-{\hat\beta \over\hat r(\hat\beta-\hat\alpha)} z^{-\hat\alpha-\hat\beta} 
F(1-\hat\beta,-\hat\alpha,1-\hat\alpha-\hat\beta,z)&
{\hat\alpha\over\hat\beta-\hat\alpha}  
F(\hat\alpha+1,\hat\beta,\hat\alpha+\hat\beta+1,z)\cr }.
$$
For $z\to 0$ it behaves like
$$
\hat Y_0\sim G_0^{\hat\alpha,\hat\beta}\, 
z^{\pmatrix{-\hat\alpha-\hat\beta&0\cr 0&0\cr}}
$$
where 
$$
G_0^{\hat\alpha,\hat\beta}={1\over\hat\beta-\hat\alpha}\pmatrix{
-\hat\alpha&\hat r\,\hat\beta\cr -{\hat\beta\over \hat r}&\hat\alpha\cr}.
$$
The monodromy around $0$ is 
$\pmatrix{\exp(-2 i\pi(\hat\alpha+\hat\beta)&0 \cr 0&1\cr}$.

For $(\tilde\Sigma)$ one has
$$
\tilde Y_0=\pmatrix{
F(\tilde\alpha,1-\tilde\alpha,1,z) &\tilde r\, 
g(\tilde\alpha,1-\tilde\alpha,1,z) \cr
-{1\over\tilde r}F(\tilde\alpha+1,-\tilde\alpha,1,z)
&-g(\tilde\alpha+1,-\tilde\alpha,1,z)\cr},
$$
for $z\to 0$ it behaves like
$$
\tilde Y_0\sim G_0^{\tilde\alpha}\,\,z^{\pmatrix{0&1\cr 0&0\cr}}
$$
with
$$
G_0^{\tilde\alpha}=\pmatrix{
1&\tilde r[\psi(1-\tilde\alpha)+\psi(\tilde\alpha)-2\psi(1)]\cr
-{1\over\tilde r}&-\psi(1+\tilde\alpha)-\psi(-\tilde\alpha)+2\psi(1)
\cr}.
$$
The monodromy around $0$ is $\pmatrix{1&2 i\pi\tilde r\cr 0&1\cr}$.

\vskip 0.3 cm
\noindent{\bf 2.2.3. Connection formulae.}\quad 
In order to compute the connection matrices we write $Y_\infty$ in the form:
$$
Y_\infty=\left(\eqalign{
&\exp(-i\pi\alpha)\,U(\alpha,\beta+1,z)\qquad\quad
{-\alpha \beta\, \exp(-i\pi(\beta+1)) r\over(\beta-\alpha)(\beta-\alpha+1)}
U(\beta+1,\alpha,z)\cr
&{-\alpha \beta\, \exp(-i\pi(\alpha+1))\over r(\beta-\alpha)(\beta-\alpha-1)}
U(\alpha+1,\beta,z)\qquad\quad
\exp(-i\pi\beta)\,U(\beta,\alpha+1,z)\cr
}\right)
$$
where $U(a,b,z):=\left(z^{-1} e^{i\pi}\right)^{-a}F(a,1-b,1+a-b,{1\over z})$. 
For $z$ such that $|{\rm arg}(z)|<2\pi$, there are the following connection 
formulae:
$$
\eqalign{
U(a,b,z)\big|_{z\rightarrow1}=
&{-\exp(i\pi a)\Gamma(1+a-b)\over\Gamma(a)\Gamma(1-b)}
\big\{ [i\pi+\psi(1-b)-\psi(b)]F(a,b,1,1-z)+\cr
&g(a,b,1,1-z)\big\},\cr}
$$
$$
U(a,1-a,z)\big|_{z\rightarrow 0}=
{-\Gamma(2a)\over\Gamma(a)^2}
\left\{ [-i\pi+\psi(a)-\psi(1-a)]F(a,1-a,1,z)+g(a,1-a,1,z)
\right\},
$$
$$
\eqalign{
U(a,b,z)\big|_{z\rightarrow0}=&
{\Gamma(1+a-b)\Gamma(1-a-b)\over\Gamma(1-b)^2}
F(a,b,a+b,z)+{\Gamma(1+a-b)\Gamma(a+b-1)\over\Gamma(a)^2}\cdot \cr
& \cdot z^{1-a-b} F(1-b,1-a,2-a-b,z). \cr}
$$
Using these relations we obtain the analytic continuation of $Y_\infty$ 
around $0$ and $1$, and by the definition of the connection matrices
$$
Y_\infty\big|_{z\rightarrow 0,1}=Y_{0,1}C_{0,1}^{\alpha,\beta,\gamma},
$$
we obtain, by straightforward computations
$$
C_0^{\hat\alpha,\hat\beta}=\pmatrix{
{\exp(i\pi\hat\beta)\Gamma(\hat\alpha+\hat\beta)\Gamma(\hat\alpha-\hat\beta+1)
\over\Gamma(\hat\alpha)\Gamma(1+\hat\alpha)} &
-\hat r{\exp(-i\pi\hat\alpha)\Gamma(\hat\alpha+\hat\beta)
\Gamma(1-\hat\alpha+\hat\beta)\over\Gamma(\hat\beta)\Gamma(1+\hat\beta)}\cr
{1\over\hat r}{\exp(-i\pi\hat\alpha)\Gamma(-\hat\alpha-\hat\beta)
\Gamma(\hat\alpha-\hat\beta+1)
\over \Gamma(1-\hat\beta)\Gamma(-\hat\beta)}&
-{\exp(-i\pi\hat\beta)\Gamma(-\hat\alpha-\hat\beta)
\Gamma(1-\hat\alpha+\hat\beta)\over\Gamma(-\hat\alpha)\Gamma(1-\hat\alpha)}
\cr},
$$
$$
C_1^{\hat\alpha,\hat\beta}=\pmatrix{
-{\Gamma(\hat\alpha-\hat\beta)\over\Gamma(\hat\alpha)\Gamma(-\hat\beta)}
[i\pi+\pi\cot(\pi \hat\alpha)]&
-{\Gamma(-\hat\alpha+\hat\beta)\over\Gamma(-\hat\alpha)\Gamma(\hat\beta)}
[i\pi+\pi\cot(\pi \hat\beta)]\cr
-{1\over\hat r}
{\Gamma(\hat\alpha-\hat\beta)\over\Gamma(\hat\alpha)\Gamma(-\hat\beta)}&
-{\Gamma(-\hat\alpha+\hat\beta)\over\Gamma(-\hat\alpha)\Gamma(\hat\beta)}\cr},
\autoeqno{f1}
$$
$$
C_0^{\tilde\alpha}=\pmatrix{
{\Gamma(2\tilde\alpha)\over\Gamma^2(\tilde\alpha)}\exp(-i\pi\tilde\alpha)
[i\pi+\pi\cot(\pi\tilde\alpha)]&
\tilde r{\Gamma(-2\tilde\alpha)\over\Gamma^2(-\tilde\alpha)}
\exp(i\pi\tilde\alpha)
[\pi\cot(\pi \hat\alpha)-i\pi]\cr
-{1\over\tilde r}{\Gamma(2\tilde\alpha)\over\Gamma^2(\tilde\alpha)}
\exp(-i\pi\tilde\alpha)
&{\Gamma(-2\tilde\alpha)\over\Gamma^2(-\tilde\alpha)}\exp(i\pi\tilde\alpha)
\cr},\autoeqno{f2}
$$
$$
C_1^{\tilde\alpha}=\pmatrix{
-{\Gamma(2\tilde\alpha)\over\Gamma^2(\tilde\alpha)}
[i\pi-\pi\cot(\pi\tilde\alpha)]&
-\tilde r{\Gamma(-2\tilde\alpha)\over\Gamma^2(-\tilde\alpha)}
[\pi\cot(\pi \hat\alpha)+i\pi]\cr
-{1\over\tilde r}{\Gamma(2\tilde\alpha)\over\Gamma^2(\tilde\alpha)}
&-{\Gamma(-2\tilde\alpha)\over\Gamma^2(-\tilde\alpha)}\cr}.
$$
Now we have to compute the monodromy matrices in the basis $Y_\infty$. Using 
the formulae \eqrefp{23}, \eqrefp{27} and \eqrefp{28} we have
$$
M_1=(C_1^{\hat\alpha,\hat\beta})^{-1} 
\pmatrix{1&2\pi i \hat r\cr 0&1\cr} 
C_1^{\hat\alpha,\hat\beta},
\qquad M_{0,x}=(C_0^{\hat\alpha,\hat\beta})^{-1} 
(C_{0,1}^{\tilde\alpha})^{-1} \pmatrix{1&2\pi i\tilde r\cr 0&1\cr} 
C_{0,1}^{\tilde\alpha} C_0^{\hat\alpha,\hat\beta}.
$$
Now we put
$$
\hat r= {-\exp[-i \pi(\hat\alpha-\hat\beta)] 
\Gamma(\hat\alpha-\hat\beta)\Gamma(\hat\beta)\Gamma(-\hat\alpha)
\over \Gamma(\hat\beta-\hat\alpha)\Gamma(-\hat\beta)\Gamma(\hat\alpha)} r 
$$
and 
$$
\tilde r = {\exp(-2 i \pi\tilde\alpha) 
\Gamma(2\tilde\alpha)\Gamma(-\tilde\alpha)^2
\over\Gamma(-2\tilde\alpha)\Gamma(\tilde\alpha)^2} 
\tilde s.
$$
In this way we immediately obtain the formula \eqrefp{m1} for $M_1$ and
it turns out that
$$
C^{\hat\alpha,\hat\beta}_0= D^{\hat\alpha,\hat\beta}\cdot C
$$ 
where $C$ is given in the formula \eqrefp{mc} and
$$
D^{\hat\alpha,\hat\beta}:=
\pmatrix{{\exp(i\pi\beta)\Gamma(\hat\beta+\hat\alpha)
\Gamma(1-\hat\beta+\hat\alpha)\over
\hat\alpha\Gamma(\hat\alpha)^2\sin\pi\hat\alpha}&0\cr
0&{-\exp(-i\pi\beta)\Gamma(-\hat\beta-\hat\alpha)
\Gamma(1+\hat\beta-\hat\alpha)\over
\Gamma(1-\hat\alpha)\Gamma(-\hat\alpha)\sin\pi\hat\alpha}\cr}.
$$
As a consequence, one has
$$
C M_{0,x} C^{-1}= (D^{\hat\alpha,\hat\beta})^{-1}
(C_{0,1}^{\tilde\alpha})^{-1} \pmatrix{1&2\pi i\tilde r\cr 0&1\cr} 
C_{0,1}^{\tilde\alpha} D^{\hat\alpha,\hat\beta}.
$$
By straightforward computations one can easily check that, for
$$
\tilde s= - {\exp(2i\pi\beta)\Gamma(\hat\beta+\hat\alpha)
\Gamma(1-\hat\beta+\hat\alpha)
\Gamma(1-\hat\alpha)\Gamma(-\hat\alpha)
\over \Gamma(-\hat\beta-\hat\alpha)
\Gamma(1+\hat\beta-\hat\alpha)
\hat\alpha\Gamma(\hat\alpha)^2} s
$$
i.e. for
$$
{\tilde r\over\hat r} = - {\Gamma(1-\sigma_0)^2\Gamma({\sigma_0\over2})^2\over
\Gamma(1+\sigma_0)^2\Gamma(-{\sigma_0\over2})^2} 
{\Gamma(1+{\thi+\sigma_0\over2})\Gamma(1+{-\thi+\sigma_0\over2})
\over \Gamma(1+{\thi-\sigma_0\over2})\Gamma(1-{\thi+\sigma_0\over2})}
{s\over r},
$$
the formulae \eqrefp{mt}, \eqrefp{m0} hold true.

To conclude the proof we have to prove the relation \eqrefp{eqs}, namely we 
want to prove that 
${\tilde r\over\hat r}=-{1\over4a_0} {2\mu+\sigma_0\over2\mu-\sigma_0}$.
To this aim we compute the matrices $A_{0,1,x}^0$ and $\Lambda$ and then 
the asymptotic behaviour of $y$ in terms of $\sigma_0$ and $\tilde r$.
To compute the matrices $A_{0,1,x}^0$ and $\Lambda$ we observe that, thanks 
to the gauges \eqrefp{22},
$$
A_1^0=\hat B_1,\qquad\Lambda=\hat B_0 + {\hat\alpha+\hat\beta\over2}\ID,
\qquad
A_{0,x}^0=G_0^{\hat\alpha,\hat\beta}\tilde B_{0,1} 
(G_0^{\hat\alpha,\hat\beta})^{-1}.
$$
First of all one has to compute the $B_{0,1}$:
$$
B_0= G_0^{\alpha,\beta,\gamma} J_0 \,(G_0^{\alpha,\beta,\gamma})^{-1},
\qquad
B_1= G_1^{\alpha,\beta,\gamma} J_1 \, (G_1^{\alpha,\beta,\gamma})^{-1}
$$
then  
$$
\hat B_0 = {1\over\hat\beta-\hat\alpha}\pmatrix{
\alpha^2&-\hat r\hat\alpha\beta\cr {\hat\alpha\beta\over\hat r}&-\hat\beta^2},
\qquad
\hat B_1 = {\hat\alpha\hat\beta\over\hat\beta-\hat\alpha}\pmatrix{
-1&\hat r\cr -{1\over\hat r}&1\cr},
$$
and
$$
\tilde B_1 ={\tilde\alpha\over2}\pmatrix{
-1&\tilde r\cr -{1\over\tilde r}&1\cr},\qquad
\tilde B_0=-{\tilde\alpha\over2}
\pmatrix{1&\tilde r\cr -{1\over\tilde r}&-1\cr}.
$$
It is then obvious that, referring to \eqrefp{ss1} and \eqrefp{ss2}, 
$b=\tilde r$, $T=G_0^{\hat\alpha,\hat\beta}$. Using the formula 
\eqrefp{p4},
$$
y(x)\sim -{\hat r(\sigma_0+2\mu)\over 4 \tilde r(2\mu-\sigma_0)}x^{1-\sigma_0}.
$$
This proves the formula \eqrefp{eqs} and concludes the proof of the theorem, 
in the case $\sigma_0\neq 0$.

For completeness we write here the result for the matrices $A_{0,1,x}^0$ and 
$\Lambda$:
$$
\Lambda={1\over4\thi} \pmatrix{
-\thi^2-\sigma_0^2&(-\thi^2+\sigma_0^2)\hat r\cr
{-\sigma_0^2+\thi^2\over\hat r}&\thi^2+\sigma_0^2\cr},\qquad
A_1^0={\thi^2-\sigma_0^2\over 4\thi}\pmatrix{
-1&\hat r\cr -{1\over \hat r}&1\cr},
$$

\leftline{$A_0^0={1\over8\theta}\left(\eqalign{
&{\theta^2-\sigma_0^2\over2}({\tilde r\over\hat r}+
{\hat r\over\tilde r})-\theta^2-\sigma_0^2 \cr
&{1\over\hat r}\left(-{(\theta-\sigma_0)^2\hat r\over2\tilde r}-
{(\theta+\sigma_0)^2\tilde r\over2\hat r}+\theta^2-\sigma_0^2\right)\cr}\right.
$}
 
\rightline{$\left.
\eqalign{\hat r\left({(\theta-\sigma_0)^2\tilde r\over2\hat r}+
{(\theta+\sigma_0)^2\hat r\over2\tilde r}+\sigma_0^2-\theta^2 \right)\cr
-{\theta^2-\sigma_0^2\over2}({\tilde r\over\hat r}+
{\hat r\over\tilde r})+\theta^2+\sigma_0^2\cr}\right)$}

\leftline{$A_x^0={1\over8\thi}\left(\eqalign{
&{\sigma_0^2-\thi^2\over2}({\tilde r\over\hat r}+
{\hat r\over\tilde r})-\thi^2-\sigma_0^2\cr
&{1\over\hat r}\left({(\thi-\sigma_0)^2\hat r\over2\tilde r}+
{(\thi+\sigma_0)^2\tilde r\over2\hat r}+ \thi^2-\sigma_0^2\right)\cr}
\right.$}

\rightline{$\left.
\eqalign{\hat r\left(
-{(\thi-\sigma_0)^2\tilde r\over2\hat r}-
{(\thi+\sigma_0)^2\hat r\over2\tilde r}+\sigma_0^2-\thi^2\right)\cr
{\thi^2-\sigma_0^2\over2}({\tilde r\over\hat r}+{\hat r\over\tilde r})
+\thi^2+\sigma_0^2\cr}\right).$}

\vskip 0.3 cm
\noindent{\bf 2.2.4. Case $\sigma_0=0$.}\quad
In this case the solution of the system $(\hat\Sigma)$ has logarithmic 
behaviour around $0$. Moreover, as seen before, it has a logarithmic 
behaviour around $1$. For this system we can use all the formulae 
derived for $(\tilde\Sigma)$, substituting $\tilde \alpha$ by $\hat\alpha$. 
The treatment of the $(\tilde\Sigma)$, is even easier. Ideed in this 
case $\Lambda$ has zero eigenvalues and it is straightforward to solve the 
system \eqrefp{24} exactly. In fact in this case we have
$$
\tilde B_0+\tilde B_1
=\pmatrix{0&1\cr 0& 0\cr},\quad 
\det\tilde B_i={\rm Tr}\tilde B_i=0,\quad i=0,1.
$$
Then the matrices $\tilde B_0$ and $\tilde B_1$ are uniquely determined up 
to an arbitrary parameter $s$:
$$
\tilde B_0=\pmatrix{0&s\cr 0& 0\cr},\quad
\tilde B_1 = \pmatrix{0&1-s\cr 0& 0\cr},
$$
and we can solve the differential equation \eqrefp{24} explicitly:
$$
\tilde Y=\pmatrix{1&s \log z +(1-s)\log(z-1)\cr 0&1\cr}.
$$
The solution $\tilde Y$ has the following asymptotic behaviour near the 
singular points:
$$
\eqalign{
\tilde Y&=\left(\ID+{\cal O}\left({1\over z}\right) \right)
z^{J},\qquad\hbox{as}\qquad z\to\infty,\cr
\,&=\tilde G_0\left(\ID+{\cal O}(z) \right)
z^J\tilde C_0,\qquad\hbox{as}\qquad z\to 0,\cr
\,&=\tilde G_1\left(\ID+{\cal O}(z-1) \right)
(z-1)^J\tilde C_1\qquad\hbox{as}\qquad z\to 1,\cr}
$$
where $J=\pmatrix{0&1\cr 0& 0\cr}$. It's easy to verify that 
$$
\tilde C_0=\pmatrix{1&0\cr 0&s\cr},\quad 
\tilde C_1=\pmatrix{1&0\cr 0&1-s\cr},\quad 
\tilde G_0=\pmatrix{1&0\cr 0&{1\over s}\cr},\quad 
\tilde G_1=\pmatrix{1&0\cr 0&{1\over 1-s}\cr}.
$$ 
As a consequence the monodromy matrices of the system \eqrefp{24} are
$$
\tilde M_0=\pmatrix{1&2\pi i s\cr 0& 1\cr},
\quad\tilde M_1=\pmatrix{1&2\pi i(1-s)\cr 0& 1\cr}.\autoeqno{31}
$$
The correspondent monodromy matrices of the full system \eqrefp{6} are 
given by:
$$
M_1=(C_1^{\hat\alpha})^{-1} \pmatrix{1&2i\pi\hat r\cr 0&1\cr}
C_1^{\hat\alpha},\quad\hbox{and}\quad
M_{0, x}= (C_0^{\hat\alpha})^{-1} \tilde M_{0, 1} C_0^{\hat\alpha},
$$
where
$$
C_0^{\hat\alpha}=\pmatrix{
{\pi\Gamma(2\hat\alpha)\over\Gamma^2(\hat\alpha)\sin\pi\hat\alpha}&
\hat r{\pi\Gamma(-2\hat\alpha)\over\Gamma^2(-\hat\alpha)\sin\pi\hat\alpha}\cr
-{1\over\hat r}{\Gamma(2\hat\alpha)\exp(-i\pi\hat\alpha)
\over\Gamma^2(\hat\alpha)}
&{\Gamma(-2\hat\alpha)\exp(i\pi\hat\alpha)\over\Gamma^2(-\hat\alpha)}\cr},
$$
and
$$
C_1^{\hat\alpha}=\pmatrix{
{\pi\Gamma(2\hat\alpha)\exp(-i\pi\hat\alpha)\over\Gamma^2(\hat\alpha)
\sin\pi\hat\alpha}
& -\hat r{\pi\Gamma(-2\hat\alpha)\exp(i\pi\hat\alpha)\over
\Gamma^2(-\hat\alpha)\sin\pi\hat\alpha}\cr
-{1\over\hat r}{\Gamma(2\hat\alpha)\over\Gamma^2(\hat\alpha)}
&-{\Gamma(-2\hat\alpha)\over\Gamma^2(-\hat\alpha)}\cr}.
$$
We observe that
$$
\hat C_{0,1}^{\hat\alpha}= C^{0,1} D^{\hat\alpha},
$$
where
$$
D^{\hat\alpha}=\pmatrix{
{\pi\Gamma(2\hat\alpha)\exp(-i\pi\hat\alpha)
\over\Gamma(\hat\alpha)^2\sin\pi\hat\alpha}&0\cr
0&{-\Gamma(-2\hat\alpha)\over\Gamma(-\hat\alpha)^2}\cr},
$$ 
and
$$
C^0=\pmatrix{\exp(i\pi\hat\alpha)&
{-\hat r\pi\over\sin\pi\hat\alpha}\cr
{-\sin\pi\hat\alpha\over\pi\hat r}&\exp(i\pi\hat\alpha)\cr},\qquad
C^1=\pmatrix{1&
{\hat r\pi\exp(i\pi\hat\alpha)\over\sin\pi\hat\alpha}\cr
{-\sin\pi\hat\alpha\exp(i\pi\hat\alpha)\over\pi\hat r}&1\cr}.
$$
We can factor out the diagonal matrix $D^{\hat\alpha}$ in \eqrefp{31}, and 
take $\hat r=1$. In this way, we obtain the formulae \eqrefp{m10}, 
\eqrefp{m00}, \eqrefp{mx0}.
The asymptotic behaviour of $y(x)$ can be computed as before. For 
$\sigma_0=0$ we obtain:
$$
y\sim a_0 x\qquad\hbox{for}\quad a_0=s.
$$
This concludes the proof of the theorem.

\vskip 0.3 cm
\noindent{\bf 2.2.5. The asymptotic behaviour near $1,\infty$ and the monodromy 
data.}\quad
We can prove the analogue of Theorem 2.2 near $1$ and $\infty$. 
Namely, for any pair of values $(a_1,\sigma_1)$ there exists a unique 
branch of the solution of PVI$\mu$ with the asymptotic behaviour 
$$
y(x)\sim\,1-a_1 x^{1-\sigma_1}\quad\hbox{as}\quad x \rightarrow 1.
\autoeqno{15'}
$$
It is possible to parameterize the monodromy matrices as in 
Theorem 2.2 substituting $\sigma_0$ with $\sigma_1$ and $M_0$ with $M_1$ 
and vice-versa.
Analogously, for any pair of values $(a_\infty,\sigma_\infty)$ there exists 
a unique branch of the solution of (PVI) with the asymptotic behaviour 
$$
y(x)\sim\,a_\infty x^{\sigma_\infty}\quad\hbox{as}\quad x \rightarrow \infty,
\autoeqno{15''}
$$
and it is possible to parameterize the monodromy matrices as before, 
substituting $\sigma_0$ with $\sigma_\infty$ and applying the braid 
$\beta_2$ to the monodromy matrices.


\vskip 0.3 cm
\noindent{\bf 2.3. From the local asymptotic behaviour to the global one.}
\vskip 0.3 cm

In this section we prove Theorem 2.1 which gives the asymptotic 
behaviour of the branches of the solutions in terms of the triplets
$(x_0, x_1,x_\infty)$.  

\proclaim Lemma 2.10. For the solution $y^{(0)}(x)$ of PVI$\mu$ behaving as 
$$
y^{(0)}(x)= a_0 x^{1-\sigma_0}\left(1+{\cal O}(x^\eps)\right)
\quad\hbox{as}\quad 
x \rightarrow 0,
$$
with $0\leq\sigma_0<1$ and $a_0\neq0$, $a_0\neq1$ for $\sigma_0=0$, the 
canonical form \eqrefp{N15} of the monodromy matrices 
$M_0^{(0)},M_x^{(0)},M_1^{(0)}$ given by \eqrefp{m0},  \eqrefp{mt}, 
\eqrefp{m1}, or \eqrefp{m00},  \eqrefp{mx0}, \eqrefp{m10} for 
$\sigma_0=0$, is the following:
$$
M_0=\pmatrix{1&-x^{(0)}_0\cr 0&1\cr},\quad M_x
=\pmatrix{1&0\cr x^{(0)}_0&1\cr},
\quad M_1=\pmatrix{1+ {x^{(0)}_1 x^{(0)}_\infty\over x^{(0)}_0}&
 -{(x^{(0)}_1)^2\over x^{(0)}_0}\cr 
{(x^{(0)}_\infty)^2\over x^{(0)}_0}&
1-{x^{(0)}_1 x^{(0)}_\infty\over x^{(0)}_0}\cr},
$$
where the triple $(x^{(0)}_0,x^{(0)}_1, x^{(0)}_\infty)$ is defined, up 
to equivalence, by the following formulae, for $\sigma_0\neq0$:
$$
\eqalign{
x_0^{(0)}&=-2\sin{\pi\sigma_0\over2},\cr
x_1^{(0)}&=-\sqrt{2(\cos\pi\sigma_0-\cos2\pi\mu)}
{\sin{\pi\phi\over2}\over\sin{\pi\sigma_0\over2}},\cr
x_\infty^{(0)}&=-\sqrt{2(\cos\pi\sigma_0-\cos2\pi\mu)}
{\cos{\pi(\sigma_0+\phi)\over2}\over\sin{\pi\sigma_0\over2}},\cr}
\autoeqno{N50}
$$
with $\phi$ given by
$$
e^{i\pi\phi}={1\over4a_0}{\sigma_0+2\mu\over\sigma_0-2\mu}
{\Gamma(1+\sigma_0)^2 \Gamma(1-{\sigma_0\over2})^2
\Gamma(1+\mu-{\sigma_0\over2})\Gamma(1-\mu-{\sigma_0\over2})
\over\Gamma(1-\sigma_0)^2 \Gamma(1+{\sigma_0\over2})^2
\Gamma(1+\mu+{\sigma_0\over2})\Gamma(1-\mu+{\sigma_0\over2})},
\autoeqno{N51}
$$
and for $\sigma_0=0$:
$$
\eqalign{
x_0^{(0)}&=0,\cr
x_1^{(0)}&=-|\sin\pi\mu| \sqrt{1-a_0}\cr
x_\infty^{(0)}&=-|\sin\pi\mu| \sqrt{a_0}.\cr}\autoeqno{N52}
$$

The proof of this lemma can be obtained by straightforward computations, 
using the algorithm of lemma 1.5. Similar formulae for the parameters 
$(x^{(1)}_0,x^{(1)}_1, x^{(1)}_\infty)$ and 
$(x^{(\infty)}_0,x^{(\infty)}_1, x^{(\infty)}_\infty)$ can be obtained 
respectively starting from a solution $y^{(1)}(x)$ of PVI$\mu$ behaving as 
$$
y^{(1)}(x)= 1-a_1(1- x)^{1-\sigma_1}\left(1+{\cal O}(1-x)^\eps\right)
\quad\hbox{as}\quad 
x \rightarrow 1,
$$
or from another solution $y^{(\infty)}(x)$ of PVI$\mu$ behaving as 
$$
y^{(\infty)}(x)= a_\infty x^{\sigma_\infty}
\left(1+{\cal O}({1\over x^\eps})\right)
\quad\hbox{as}\quad x \rightarrow \infty.
$$
So, given an admissible triple $(x_0, x_1,x_\infty)$, with $x_i\in\reali$, 
$|x_i|<2$ for $i=0,1,\infty$, we choose the parameters $\mu$,  
$(a_0,\sigma_0)$, $(a_1,\sigma_1)$ and $(a_\infty,\sigma_\infty)$ in such 
a way that \eqrefp{B6} is satisfied and
$$
x_i^{(0)}=x_i^{(1)}=x_i^{(\infty)}=x_i,\quad\hbox{for}\quad i=0,1,\infty.
$$
Using the explicit formulae \eqrefp{N50}, \eqrefp{N51} for $x_0\neq 0$, 
we derive the expressions \eqrefp{A1}, \eqrefp{A2}. Similarly, using 
\eqrefp{N52} for $x_0=0$ we derive the expression \eqrefp{A3}. In the same 
way, we derive the analogous expressions for $(a_1,\sigma_1)$ and 
$(a_\infty,\sigma_\infty)$. The three correspondent branches 
$y^{(0)}(x),y^{(1)}(x), y^{(\infty)}(x)$ of solutions of PVI$\mu$, 
with $\mu$ given by \eqrefp{B6} must coincide. In fact, the associated 
auxiliary Fuchsian systems have the same, modulo diagonal conjugation, 
monodromy matrices. This proves the existence of a solution of PVI$\mu$ 
with the asymptotic behaviour \eqrefp{BB1}, with the indices given by 
\eqrefp{B5} and the coefficients specified as above, for any admissible 
triple $(x_0, x_1,x_\infty)$, with $x_i\in\reali$, $|x_i|<2$ for 
$i=0,1,\infty$. The uniqueness of such a branch follows from theorem 1.3.

Conversely, for any such a solution we obtain an admissible triple 
$(x_0, x_1,x_\infty)=(x^{(0)}_0, x^{(0)}_1,x^{(0)}_\infty)
=(x^{(1)}_0, x^{(1)}_1,x^{(1)}_\infty)
=(x^{(\infty)}_0, x^{(\infty)}_1,x^{(\infty)}_\infty)$, 
using the formulae \eqrefp{N50}, \eqrefp{N51} or \eqrefp{N52} and their 
analogies. Let us prove that the numbers $(x_0, x_1,x_\infty)$ are real 
and satisfy $|x_i|<2$ for $i=0,1,\infty$. Indeed, from the definition of 
the parameters, it follows:
$$
(x^{(0)}_0)^2=4 \sin^2\pi\sigma_0,
( x^{(1)}_1)^2=4 \sin^2\pi\sigma_1,
(x^{(\infty)}_\infty)^2=4 \sin^2\pi\sigma_\infty.
$$
This proves that our construction covers, for real $\mu$, all the solutions 
of PVI$\mu$ with critical behaviour of algebraic type.

Finally, using corollary 1.5, we infer that the class of solutions of 
PVI$\mu$, with real $\mu$, having  critical behaviour of algebraic type 
is invariant with respect to the analytic continuation. The law of 
transformation of the critical indices $l_0,l_1,l_\infty$ of the 
expansions \eqrefp{BB1}, is described by theorem 1.5.

\vfill\eject


\noindent{\bf 2.4. The complete list of algebraic solutions}
\vskip 0.3 cm

We summarize the results of this paper in the following

\proclaim Classification Theorem. Any algebraic solution of the equation 
PVI$\mu$ with $2\mu\not\in\interi$ is equivalent, in the sense of symmetries
\eqrefp{s1}, \eqrefp{s2}, \eqrefp{B10} to one of the five solutions 
$(A_3)$, $(B_3)$, $(H_3)$, $(H_3)'$, $(H_3)''$ below.

We already know that the classes of equivalent algebraic solutions are 
labelled by the five regular polyhedra and star-polyhedra in the 
three-dimensional space. We will construct representatives in these classes 
for the following values of the parameter $\mu$
$$
\mu=-{1\over4},\quad-{1\over3},\quad-{2\over5},\quad-{1\over5},
\quad-{1\over3}.
$$
The correspondent algebraic solutions will have $4$, $3$, $10$, $10$, $18$
branches respectively. Recall that these are the lenghts of the orbits 
\eqrefp{cl1}, \eqrefp{cl2}, \eqrefp{cl3}, \eqrefp{cl4}, \eqrefp{cl5} 
respectively with respect to the action of the pure braid group (see 
remark 1.11 above). We give now the explicit formulae for the solutions with 
brief explanations of the derivations of them.

\vskip 0.2 cm
\item{} {\bf Tetrahedron.}
\quad We have $(x_0,x_1,x_\infty) = (-1,0,-1)$, then $\mu=-{1\over4}$ and
$$
\eqalign{
&y= {(s-1)^2(1+3s)(9s^2-5)^2\over(1+s)(25-207 s^2+1539 s^4+243 s^6)},\cr
&x={(s-1)^3(1+3s)\over(s+1)^3(1-3s)}.\cr}\eqno{(A_3)}
$$
(We present the solution in the parametric form). The monodromy matrices, 
in the canonical form \eqrefp{N15}, are:
$$
M_0=\pmatrix{1&1\cr 0&1\cr}\quad M_x=\pmatrix{1&0\cr -1&1\cr}\quad 
M_1=\pmatrix{1&0\cr 1&1\cr}. 
$$
This solution was found in [Dub] in the implicit form (E.29). This was also 
obtained, independently, by N. Hitchin (see [Hit2]). To reduce (E.29) to the 
above form, we have to solve the cubic equation (E.29 b) with the 
substitution:
$$
t={32(1-18 s^2+81 s^4)\over 27(1+9s^2+27 s^4 + 27 s^6)}.
$$
Then the three roots of (E.29 b) are:
$$
\eqalign{
\omega_1&={13-66 s^2-27 s^4\over3(1+3 s^2)^2},\cr
\omega_{2,3}&={-5+42 s^2\pm 144 s^3 + 27 s^4\over3(1+3 s^2)^2}.\cr}
$$

\vskip 0.2 cm
\item{}{\bf Cube.}\quad 
We have $(x_0,x_1,x_\infty) =  (-1,0,-\sqrt 2)$ and $\mu=-{1\over3}$. The 
solution
$$
\eqalign{
&y={(2-s)^2(1+s)\over(2+s)(5s^4-10s^2+9)}, \cr
&x={(2-s)^2(1+s)\over(2+s)^2(1-s)},\cr}\eqno{(B_3)}
$$
was obtained in [Dub]. The canonical form for the monodromy matrices is:
$$
M_0=\pmatrix{1&1\cr 0&1\cr}\quad M_x=\pmatrix{1&0\cr -1&1\cr}\quad 
M_1=\pmatrix{1&0\cr -2&1\cr}. 
$$

\vskip 0.3 cm

\noindent{\bf Coxeter group ${\bf W(H_3)}$, of symmetries of 
icosahedron.}\quad We have three possible choices of the point  
$(x_0,x_1,x_\infty)$ which lead to three different solutions. 
\vskip 0.2 cm
\item{}{\bf Icosahedron.}\quad The orbit \eqrefp{cl3} corresponds to the 
standard triple of reflections for the icosahedron. 
$(x_0,x_1,x_\infty)=(0,-1,-{1+\sqrt 5\over 2})$, then $\mu=-{2\over5}$ and
$$
\eqalign{
&y={(s-1)^2(1+3s)^2(-1+4 s+s^2)(7-108 s^2+314 s^4-588 s^6+119 s^8)^2\over
 (1+s)^3(-1+3s) P(s)}\cr
&x={{{{\left( -1 + s \right) }^5}\,{{\left( 1 + 3\,s \right) }^3}\,
     \left( -1 + 4\,s + {s^2} \right) }\over 
   {{{\left( 1 + s \right) }^5}\,{{\left( -1 + 3\,s \right) }^3}\,
     \left( -1 - 4\,s + {s^2} \right) }},\cr}\eqno{(H_3)}
$$
with
$$
\eqalign{
P(s)=&49-2133 s^2+34308 s^4-259044 s^6+16422878 s^8-7616646 s^{10}+
13758708 s^{12}\cr
&+5963724 s^{14}-719271 s^{16}+42483 s^{18}.\cr}
$$
The canonical form for the monodromy matrices is:
$$
M_0=\pmatrix{1&1\cr 0&1\cr}\quad M_x=\pmatrix{1&0\cr -1&1\cr}\quad 
M_1=\pmatrix{1&0\cr {-3-\sqrt5\over2}&1\cr}. 
$$
The above solution was already obtained in [Dub] in the implicit form 
(E.33). The above explicit formula can be obtained solving (E.33 b) in 
the form:
$$
t={(1-4s-s^2)(-1-4s+s^2)(-1+5s^2)\over (1+3s^2)^3}
$$
$$
\eqalign{\omega_1 &=
{{25 - 585\,{s^2} + 3530\,{s^4} - 6690\,{s^6} - 3955\,{s^8} + 
     507\,{s^{10}}}\over {{{\left( 1 + 3\,{s^2} \right) }^5}}}\cr
\omega_2 &=
{{-7 + 215\,{s^2} - 1910\,{s^4} - 4096\,{s^5} + 5150\,{s^6} + 20480\,{s^7} + 
     6125\,{s^8} - 357\,{s^{10}}}\over {{{\left( 1 + 3\,{s^2} \right) 
}^5}}}\cr
\omega_3 &=
{{-7 + 215\,{s^2} - 1910\,{s^4} + 4096\,{s^5} + 5150\,{s^6} - 20480\,{s^7} + 
     6125\,{s^8} - 357\,{s^{10}}}\over {{{\left( 1 + 3\,{s^2} \right) }^5}}.}
\cr}
$$

The last two solutions for the orbits \eqrefp{cl4} and \eqrefp{cl5}, with the 
icosahedral symmetry are new. They correspond to great icosahedron and great 
dodecahedron respectively. To compute 
them we use the following algorithm. The leading terms of the Puiseux 
expansions near the ramification points $0,1,\infty$ of each branch can 
be computed by the formulae \eqrefp{B5}, \eqrefp{A1}, \eqrefp{A2} and 
\eqrefp{A3}. From this the genus of the algebraic curve $F(y,x)=0$ is 
easily computed. Namely, the genus of \eqrefp{cl4} is $0$ and the genus 
of \eqrefp{cl5} is $1$. Since the symmetries of PVI$\mu$ preserve the 
indices $l_0,l_1,l_\infty$ (up to permutations), they preserve the genus too. 

We observe that the appearance of genus $1$ in the last 
solution related to the great dodecahedron could seem less surprising 
if we recall that the topology of this immersed two-dimensional surface 
is different from the topology of all the other polyhedra and star-polyhedra. 
In fact, this is a surface of genus $4$, while all the others have genus 
$0$ (see [Cox]).

Let us now list the last two solutions.

\vskip 0.2 cm
\item{}{\bf Great Icosahedron.}\quad 
$(x_0,x_1,x_\infty)= (-1,0,{1-\sqrt 5\over 2})$, then $\mu=-{1\over5}$ and
$$
\eqalign{
&y={{{{\left( -1 + s \right) }^4}\,{{\left( 1 + 3\,s \right) }^2}\,
     \left( -1 + 4\,s + {s^2} \right) \,
     {{\left( 3 - 30\,{s^2} + 11\,{s^4} \right) }^2}}\over 
   {\left( 1 + s \right) \,\left( -1 + 3\,s \right) \,
     \left( 1 + 3\,{s^2} \right)\, P(s)}},\cr
&x={{{{\left( -1 + s \right) }^5}\,{{\left( 1 + 3\,s \right) }^3}\,
     \left( -1 + 4\,s + {s^2} \right) }\over 
   {{{\left( 1 + s \right) }^5}\,{{\left( -1 + 3\,s \right) }^3}\,
     \left( -1 - 4\,s + {s^2} \right) }}.\cr}\eqno{(H_3)'}
$$
with
$$
P(s)=
\left( 9 - 342\,{s^2} + 4855\,{s^4} - 28852\,{s^6} + 63015\,{s^8} - 
       1942\,{s^{10}} + 121\,{s^{12}} \right).
$$ 
The canonical form for the monodromy matrices is
$$
M_0=\pmatrix{1&1\cr 0&1\cr}\quad M_x=\pmatrix{1&0\cr -1&1\cr}\quad 
M_1=\pmatrix{1&0\cr {-3+\sqrt5\over2}&1\cr}. 
$$

\vskip 0.2 cm
\item{}{\bf Great Dodecahedron.}\quad 
$(x_0,x_1,x_\infty)=(-1,-1,{1-\sqrt{5}\over2})$, $\mu=-{1\over3}$. The 
canonical form for the monodromy matrices is
$$
M_0=\pmatrix{1&1\cr 0&1\cr}\quad 
M_x=\pmatrix{1&0\cr -1&1\cr}\quad 
M_1=\pmatrix{{3-\sqrt{5}\over 2} &1\cr 
{-3+\sqrt{5}\over 2} &{1+\sqrt{5}\over 2 }\cr}. 
$$
This is the most complicated solution and we will briefly explain how did we 
obtain it. As we already said, it is an algebraic function with $18$ branches. 
It has two branch points of order $5$, two of order $3$ and two regular 
branches, over every ramification point $0,1,\infty$. The 
branches $y_1(x),\cdots,y_{18}(x)$ near $x=0$ have the form:
$$
\eqalign{ 
y_k(x)&= e^{2\pi i k\over5}\left({7\over13}\right)^2 6^{-2\over5}x^{4\over5} +
{\cal O}(x), \qquad k=1,\cdots,5\cr
y_{k+5}(x)&= e^{2\pi i k\over5} {6^{4\over5}\over 19^2} x^{2\over5} +
{\cal O}(x^{4\over5}), \qquad k=1,\cdots,5\cr
y_{10+k}(x)&= e^{2\pi i k\over3} {2^{2\over3}\over 18} {1+i\sqrt{15}\over4}  
x^{2\over3} + {\cal O}(x), \qquad k=1,\cdots,3\cr
y_{13+k}(x)&= e^{2\pi i k\over3} {2^{2\over3}\over 18} {1-i\sqrt{15}\over4} 
x^{2\over3} + {\cal O}(x), \qquad k=1,\cdots,3\cr
y_{17,18}(x)&= {3\pm\sqrt{5}\over 6} x+ {\cal O}(x^2).  \cr
}$$
The Puiseux expansions near $x=1$ and $x=\infty$ can be obtained from these 
formulae applying the symmetries \eqrefp{s1} and \eqrefp{s2} respectively. 
Using these formulae, one can compute any term of the 
Puiseux expansions of all the branches. Due to computer difficulties, at the 
moment, we do not manage to produce the explicit elliptic parameterization 
of the algebraic curve. We give this in the form of an algebraic curve of 
degree $36$. 
$$
x^{15} F(x,y,t)=0 \eqno{(H_3)''}
$$
where
$$
t=x+{1\over x}
$$
and
$$
F(x,y,t)=\left( 11423613917539180989 - 57169813730203944\,t - 
        13869163074392577\,{t^2} \right.
$$
$$
\left. + 1307302991918736\,{t^3} - 
        31962210377\,{t^4} - 556854952\,{t^5} + 282475249\,{t^6} \right) ^2
x^9 + 
$$
$$
+ 9\,\left( -42194267411458338799378785573556538817 
\right.
$$
$$
-58759262104428568315429822622247510492\,t +
$$
$$+10095266581644469686796601774497789110\,{t^2} 
$$
$$- 
     969805106597038829472153249647160780\,{t^3} 
$$
$$+ 
     13082239583395373581545441399627177\,{t^4} 
$$
$$- 
     77058446549850745165440956773416\,{t^5} 
$$
$$- 
     2150599531632473735225276196788\,{t^6} 
$$
$$+ 
     5521397776112060589691860200\,{t^7} 
$$
$$+ 
     34431689430132242698256649\,{t^8} - 4868379539328005204126748\,{t^9} 
$$
$$+ 
     543298990997997546590\,{t^{10}} - 5420393254540081020\,{t^{11}} 
$$
$$\left. - 
     79792266297612001\,{t^{12}} \right) \,{x^8}\,\left( 1 + x \right) \,y 
$$
$$+ 3\,\left( 1958540422461728364360419152068949457061 
\right.
$$
$$- 
     2209393132972329408615760780334959957197\,t 
$$
$$+ 
     3730909713982979160856372675878664915614\,{t^2} 
$$
$$- 
     998456940153012666787445238400320871842\,{t^3} 
$$
$$+ 
     127635631698846877473536225998246411623\,{t^4} 
$$
$$- 
     1482254497297246657880332587875322027\,{t^5} 
$$
$$+ 
     8847488219466307166390055942913100\,{t^6} 
$$
$$+ 
     115442527212405524632938663371228\,{t^7} 
$$
$$- 
     1429542565820801871766655325509\,{t^8} 
$$
$$+ 
     7631245862019705473945545029\,{t^9} 
$$
$$+ 
     215859931310548879418863190\,{t^{10}} 
$$
$$- 
     61369918621176379581242\,{t^{11}} 
$$
$$+ 443864075501101414537\,{t^{12}} 
$$
$$\left. + 
     239376798892836003\,{t^{13}} \right) \,{x^8}\,{y^2} 
$$
$$+ 8\,\left( -14176454404869485146893293606614356611775 
\right. $$
$$+ 
     28055900200716580071561317517835005937542\,t 
$$
$$- 
     22805214425386778204880821526748969126548\,{t^2} 
$$
$$+ 
     7203824414325261553629527716660045847442\,{t^3} 
$$
$$- 
     941980967212371956151450075677431178469\,{t^4} 
$$
$$+ 
     8038149376577574092957422850094452382\,{t^5} 
$$
$$- 
     22951481031624768208910085230213950\,{t^6} 
$$
$$- 
     1567265072303229457696028497735362\,{t^7} 
$$
$$+ 
     17710721122570843039502588273105\,{t^8} 
$$
$$- 
     100631536944878626686735036764\,{t^9} 
$$
$$- 
     513222963217440801801106854\,{t^{10}} 
$$
$$+ 
     346131811374226777113368\,{t^{11}} 
$$
$$\left. - 1385596068253512936373\,{t^{12}}
      \right) \,{x^7}\,\left( 1 + x \right) \,{y^3} 
$$
$$+6\,\left( 372245126038285018174621123839906354129684 
\right.
$$
$$- 
     538123164074188598920246979212739683700019\,t 
$$
$$+ 
     148097514409311992531878851512796911164392\,{t^2} 
$$
$$+ 
     131183971273631981690818920653550831952470\,{t^3} 
$$
$$- 
     81255871706326021946002239342890704787124\,{t^4} 
$$
$$+ 
     13496764871847929345618085983834973142347\,{t^5} 
$$
$$- 
     77522788387610267024773351181536807752\,{t^6} 
$$
$$- 
     15773276759106249673670395386220524\,{t^7} 
$$
$$+ 
     43306225979803766149728090033557588\,{t^8} 
$$
$$- 
     398429371182004923439978792917533\,{t^9} 
$$
$$+ 
     1762391142835398582269106251136\,{t^{10}} 
$$
$$+ 
     2388970787440139757429804470\,{t^{11}} 
$$
$$- 
     3439831965306734927530580\,{t^{12}} 
$$
$$\left. + 6227088023782555081845\,{t^{13}}
      \right) \,{x^7}\,{y^4} 
$$
$$+ 6\,\left( -1723593607166532080927038166074395946035900 
\right.$$
$$+ 
     2040599178888751273951239802196511605636025\,t 
$$
$$+ 
     233403452055579960658996997839662220718652\,{t^2} 
$$
$$- 
     1266443503829326206558682114543567079698662\,{t^3} 
$$
$$+ 
     601472177228487289796456276889218593842200\,{t^4} 
$$
$$- 
     90437995682218432847116536290883275511849\,{t^5} 
$$
$$+ 
     312726298388429638276816726813440288752\,{t^6} 
$$
$$+ 
     2550025722128829156352333577788952604\,{t^7} 
$$
$$- 
     575104019542617989271403929336386084\,{t^8} 
$$
$$+ 
     4033434111720024901257683928676919\,{t^9} 
$$
$$- 
     12117519775974826925705049828908\,{t^{10}} 
$$
$$- 
     6556200250126220600244219830\,{t^{11}} 
$$
$$+ 
     14241314959493049268655304\,{t^{12}} 
$$
$$\left. - 9850799841945252913527\,{t^{13}}
      \right) \,{x^6}\,\left( 1 + x \right) \,{y^5} 
$$
$$+2\,\left( 21644465131825357400382632120971734649857264 
\right.
$$
$$+ 
     1609281090760898028021044308156991213714259\,t 
$$
$$- 
     38819578264561008264384235926474662382590097\,{t^2} 
$$
$$+ 
     19461702669444598007359556763777303231431634\,{t^3} 
$$
$$+ 
     6809228919729818789589103927528063675916654\,{t^4} 
$$
$$- 
     6502005871163632747227392004794985000703103\,{t^5} 
$$
$$+ 
     1176880351448958476049654953342289787508545\,{t^6} 
$$
$$- 
     2305666451065939258487468463173953769764\,{t^7} 
$$
$$- 
     66480299243646643579685794767909280748\,{t^8} 
$$
$$+ 
     15241200180683229312383060139285323669\,{t^9} 
$$
$$- 
     78386104021612006912908252437015263\,{t^{10}} 
$$
$$+ 
     148880803913659745854808576095474\,{t^{11}} 
$$
$$+ 
     57461755817453559455002717038\,{t^{12}} 
$$
$$- 
     103978838096755580164874073\,{t^{13}} 
$$
$$\left. + 22079893897296692642367\,{t^{14}}
      \right) \,{x^6}\,{y^6} 
$$
$$+ 12\,\left( -3846090734682844350682028823106084687428444 
\right.
$$
$$- 
     6334186503545807261781842904989156398103475\,t 
$$
$$+ 
     13149149531444644750567430723972640952739601\,{t^2} 
$$
$$- 
     833478059474840476986699901359926517619482\,{t^3} 
$$
$$- 
     6957568559854570720866520753764257206265698\,{t^4} 
$$
$$+ 
     3704223179468379920674190362976879513996499\,{t^5} 
$$
$$- 
     576109607175627912516018268715191958362489\,{t^6} 
$$
$$+ 
     449528751470785721434087640482632944192\,{t^7} 
$$
$$+ 
     97703132604175842587572789513022822852\,{t^8} 
$$
$$- 
     16011119124352462144696375400709770865\,{t^9} 
$$
$$+ 
     57470033726311385832789992330535883\,{t^{10}} 
$$
$$- 
     64186326557416461737056084783958\,{t^{11}} 
$$
$$- 
     26094692352698228396829172726\,{t^{12}} 
$$
$$+ 
     25456979816863844482106097\,{t^{13}} 
$$
$$\left. - 1051423518918890125827\,{t^{14}}
      \right) \,{x^5}\,\left( 1 + x \right) \,{y^7} 
$$
$$+ 3\,\left( 10610560214390717981593236262575236159801442 
\right.
$$
$$+ 
     111851425974428655491946571184648885748780846\,t 
$$
$$- 
     46715489836280837819492546495478177251546881\,{t^2} 
$$
$$- 
     139030514021835704781706562226557058780091068\,{t^3} 
$$
$$+ 
     91919803393449431080407081029445741045516160\,{t^4} 
$$
$$+ 
     17286168898942635534811373112001322602942242\,{t^5} 
$$
$$- 
     24220067368877156108293128216367359106738667\,{t^6} 
$$
$$+ 
     4637848726272650440828320863095758089676600\,{t^7} 
$$
$$- 
     768878602349792937649689712426759282462\,{t^8} 
$$
$$- 
     1741702208080601219941076960674171058190\,{t^9} 
$$
$$+ 
     291817792658031512381548036217335383393\,{t^{10}} 
$$
$$- 
     703979549232767832584811827025318204\,{t^{11}} 
$$
$$+ 
     431869052449224832223544125220252\,{t^{12}} 
$$
$$+ 
     211026301709949401515207349982\,{t^{13}} 
$$
$$\left. - 
     88672767468343194407328645\,{t^{14}} \right) \,{x^5}\,{y^8} 
$$
$$+\left( 29516093700614758561713684532397059536112914 
\right.
$$
$$- 
     250427363718702318175164789552712875171021018\,t 
$$
$$- 
     276026230387316277405019906885146213095909789\,{t^2} 
$$
$$+ 
     798207668051094993868151955635553085940964248\,{t^3} 
$$
$$- 
     344467523172443772269951999752021715572380244\,{t^4} 
$$
$$- 
     120301086904813944940738235485125501727252554\,{t^5} 
$$
$$+ 
     109511261517293051818507102282993534830623481\,{t^6} 
$$
$$- 
     19194080365623650515069714757340701033233072\,{t^7} 
$$
$$- 
     51032394918091643798626435084286655402518\,{t^8} 
$$
$$+ 
     23520399409116636731567242035209841281394\,{t^9} 
$$
$$- 
     2981630391452674856104059990243312091171\,{t^{10}} 
$$
$$+ 
     4599349671236813290177817558881705144\,{t^{11}} 
$$
$$- 
     1416023622121443571963532624008680\,{t^{12}} 
$$
$$- 
     884464824710945312016748613390\,{t^{13}} 
$$
$$\left. + 
     126470477973158698785100695\,{t^{14}} \right) \,{x^4}\,
   \left( 1 + x \right) \,{y^9} 
$$
$$+3\,\left( -10911639347758887707980330460476045164474906 
\right.
$$
$$- 
     58595420545899882777261835025201788418411880\,t 
$$
$$+ 
     510148671470006459933001332369933260295921507\,{t^2} 
$$
$$- 
     421411310692898883818707378130970229534690273\,{t^3} 
$$
$$- 
     167152535194408490298020642999503515413037404\,{t^4} 
$$
$$+ 
     259430867087315198025819011038004422708530022\,{t^5} 
$$
$$- 
     46562493376530130089912033636689480134093531\,{t^6} 
$$
$$- 
     21317412351831533402564088853511437604716591\,{t^7} 
$$
$$+ 
     5998321070650622801147765820299240777440702\,{t^8} 
$$
$$+ 
     42318724644239526363080933921378429784628\,{t^9} 
$$
$$- 
     20025402467700581457477461594500061379211\,{t^{10}} 
$$
$$+ 
     2543243908066676381203685445622311376505\,{t^{11}} 
$$
$$- 
     2382710944926354080442398446589552744\,{t^{12}} 
$$
$$+ 
     326767340108407733203130499898126\,{t^{13}} 
$$
$$+ 
     276491632581630260746191363123\,{t^{14}} 
$$
$$\left. - 
     8431365198210579919006713\,{t^{15}} \right) \,{x^4}\,{y^{10}} 
$$
$$+12\,\left( -8733106988589310312964985955814462819079774 
\right.
$$
$$+ 
     64016502559572741731505856556403573520814046\,t 
$$
$$- 
     144735875200699985710584864792377620711111027\,{t^2} 
$$
$$+ 
     66707682241185015299381142924382044312819100\,{t^3} 
$$
$$+ 
     73079504596828994404469155513764487920535332\,{t^4} 
$$
$$- 
     75764221368814908722894464652763782420369670\,{t^5} 
$$
$$+ 
     18144429124108777614777763808360976174204563\,{t^6} 
$$
$$+ 
     1935278311931502378271609920389021155764328\,{t^7} 
$$
$$- 
     842229927380965714885405569219684268010166\,{t^8} 
$$
$$- 
     43576284807293907660657241391770657495742\,{t^9} 
$$
$$+ 
     12197493655076249609592887507484958732619\,{t^{10}} 
$$
$$- 
     1236906021392508491488029828633988502116\,{t^{11}} 
$$
$$+ 
     656542960338462752927669545849000176\,{t^{12}} 
$$
$$- 
     33306604264547610372382676649434\,{t^{13}} 
$$
$$\left. - 
     41617303783672454910317327355\,{t^{14}} \right) \,{x^3}\,
   \left( 1 + x \right) \,{y^{11}} 
$$
$$+2\,\left( 187554318014748213259275472412247634165779260 
\right.
$$
$$- 
     762602529289184570716467664651505597594612711\,t 
$$
$$+ 
     781137281074558579579982041608982979326338521\,{t^2} 
$$
$$+ 
     609916709974782293997005108542677011211842354\,{t^3} 
$$
$$- 
     1300509731487820027375772452301766352965003310\,{t^4} 
$$
$$+ 
     495021556835245340779707559257082787835603199\,{t^5} 
$$
$$+ 
     115214340026017185187174532535570472667840767\,{t^6} 
$$
$$- 
     93151202230256791649185399844658098150754584\,{t^7} 
$$
$$+ 
     7202343146876999296118334182261935450026080\,{t^8} 
$$
$$+ 
     1570014403002611359784317587694467895480399\,{t^9} 
$$
$$+ 
     371422610102654444897802704145750294582991\,{t^{10}} 
$$
$$- 
     107776513895628817503345944747365689465454\,{t^{11}} 
$$
$$+ 
     10941070332070115448154615812894796283346\,{t^{12}} 
$$
$$- 
     3005301616251624310474021462023025751\,{t^{13}} 
$$
$$+ 
     41932897601387821985740913526761\,{t^{14}} 
$$
$$\left. + 
     87168052345572726428101766292\,{t^{15}} \right) \,{x^3}\,{y^{12}} 
$$
$$+6\,\left( -49769698820686378083743112964571467299790844 
\right.
$$
$$+ 
     144208424244836145109900273916270122534480379\,t 
$$
$$- 
     31061701177281181816778329394184658654578125\,{t^2} 
$$
$$- 
     316115877280896274823885457480018099381675270\,{t^3} 
$$
$$+ 
     413355491574586200245516861333774692564805758\,{t^4} 
$$
$$- 
     189323246544195194714568800051070616915145871\,{t^5} 
$$
$$+ 
     16358105769812410561563772201680559507177601\,{t^6} 
$$
$$+ 
     10699064629034472787428641330738981639504716\,{t^7} 
$$
$$- 
     2410111354335748769859016926649172233001216\,{t^8} 
$$
$$+ 
     411171528032730820725099829153187243925077\,{t^9} 
$$
$$- 
     226789949801107593485509158636442297272707\,{t^{10}} 
$$
$$+ 
     47843156588060602804351104035547186862994\,{t^{11}} 
$$
$$- 
     4050842729096500910242692387785686802242\,{t^{12}} 
$$
$$+ 
     503056665812777017090193080924335791\,{t^{13}} 
$$
$$- 
     1358120473689033753999092690625\,{t^{14}} 
$$
$$\left. - 
     4470156530542191098877013656\,{t^{15}} \right) \,{x^2}\,
   \left( 1 + x \right) \,{y^{13}} 
$$
$$+6\,\left( 40840375974497844675523416705994772589613816 
\right.
$$
$$- 
     31065278771086485580585566101107070670979537\,t 
$$
$$- 
     256058709793754221748465779158452450846012816\,{t^2} 
$$
$$+ 
     559068508134217391316137877579221486889948317\,{t^3} 
$$
$$- 
     405356674135473476033247305191846353227354396\,{t^4} 
$$
$$+ 
     56661024338479719866449660305901198060970619\,{t^5} 
$$
$$+ 
     64125926679095609696574105582575413252584004\,{t^6} 
$$
$$- 
     28890521413485617192612359619257982301215223\,{t^7} 
$$
$$+ 
     2839341861289460260556040613235064311172368\,{t^8} 
$$
$$+ 
     465313317652398071330581978147640308517821\,{t^9} 
$$
$$- 
     351597845897601623389761113030041346811064\,{t^{10}} 
$$
$$+ 
     179321665024113431049966683538039741491991\,{t^{11}} 
$$
$$- 
     39080763389941286051140782872898013637420\,{t^{12}} 
$$
$$+ 
     3329327789375495478533443619894892092425\,{t^{13}} 
$$
$$- 
     149553232132910627041355116122361020\,{t^{14}} 
$$
$$\left. + 
     205192601853360299747062918515\,{t^{15}} \right) \,{x^2}\,{y^{14}} 
$$
$$+4\,\left( -5628131126860710389875035169789649454420876 
\right.
$$
$$- 
     78804574894375995596657603789677394324984101\,t 
$$
$$+ 
   352860279249727357643015437757245513637138672\,{t^2} 
$$
$$- 
     564941261129218367227412931528430799748606546\,{t^3} 
$$
$$+ 
     430939089346630550106464017811753376858010640\,{t^4} 
$$
$$- 
     148284188482685309994113608598302973692022739\,{t^5} 
$$
$$+ 
     4873973943047627527185771161170310871611116\,{t^6} 
$$
$$+ 
     10181461701693419091610691800269773692583228\,{t^7} 
$$
$$- 
     1226684412907984419281022032089194096771900\,{t^8} 
$$
$$- 
     1114701349894370233505605371103641055314707\,{t^9} 
$$
$$+ 
     706698148832598485833137372995728746006888\,{t^{10}} 
$$
$$- 
     230885597278675059768074093486733449982986\,{t^{11}} 
$$
$$+ 
     40110760213781966306595755424591426952408\,{t^{12}} 
$$
$$- 
     2944406938738808019234484282441173992613\,{t^{13}} 
$$
$$\left. + 
     29909989810256194655311832623132956\,{t^{14}} \right) \,x\,
   \left( 1 + x \right) \,{y^{15}} 
$$
$$+3\,\left( -19345311524103689299806429866595584344434933 
\right.
$$
$$+ 
     165880840018062517894524148661179410853072546\,t 
$$
$$- 
     433975351186661527899190510419861031681577223\,{t^2} 
$$
$$+ 
     515516306674309051714096086492072331808918060\,{t^3} 
$$
$$- 
     283562876761607595979024343783955270990852289\,{t^4} 
$$
$$+ 
     35089717870652037166528865782071242284918734\,{t^5} 
$$
$$+ 
     33297928990127187049831304457387943687578909\,{t^6} 
$$
$$- 
     12917764244851664872827620472556082803226856\,{t^7} 
$$
$$- 
     266713623245328356955979252488258143292463\,{t^8} 
$$
$$+ 
     555900198844440351814987030522263162652334\,{t^9} 
$$
$$+ 
     344809125199575823496923125385565831315595\,{t^{10}} 
$$
$$- 
     325689072459807008457121908075371991483716\,{t^{11}} 
$$
$$+ 
     117388439783020206894897144460070846332949\,{t^{12}} 
$$
$$- 
     21123688072686368568170196496753937437182\,{t^{13}} 
$$
$$\left. + 
     1569161588742434760282235480090100082255\,{t^{14}} \right) \,x\,{y^{16}} 
$$
$$+ 3\,\left( 9783299760488948030219433006083570296689357 
\right.
$$
$$- 
     59321119347918543659930676521984384042169430\,t 
$$
$$+ 
     141416477837529651726686264572772822193430055\,{t^2} 
$$
$$- 
     177096809878289456793903796377476455257673500\,{t^3} 
$$
$$+ 
     127907586479651422318564410835908192786763365\,{t^4} 
$$
$$- 
     54372658309139640733439296021048049726746698\,{t^5} 
$$
$$+ 
     13488394375983259178386269031077826541323679\,{t^6} 
$$
$$- 
     2113244794203694376441274534558687456380488\,{t^7} 
$$
$$+ 
     141785257824097311610019381070069013792095\,{t^8} 
$$
$$+ 
     316821091130460893567937727374441119017078\,{t^9} 
$$
$$- 
     305480360931755555721215775431316200256739\,{t^{10}} 
$$
$$+ 
     142141070595224470100170760768533902542116\,{t^{11}} 
$$
$$- 
     38457837145846954116338584809621985652097\,{t^{12}} 
$$
$$+ 
     5806436836462494743682658146129894324810\,{t^{13}} 
$$
$$\left. - 
     380650359326333515862984779019865187923\,{t^{14}} \right) \,
   \left( 1 + x \right) \,{y^{17}} 
$$
$$+ \left( -10524240525647109159259219575804205851284241 
\right.
$$
$$+ 
     54473554541130948618895564303322618142908415\,t 
$$
$$- 
     116506041709060140591221129409148683447467665\,{t^2} 
$$
$$+ 
     134582692171688928175407226575844795641654895\,{t^3} 
$$
$$- 
     91188532083602481766920967127844408926436405\,{t^4} 
$$
$$+ 
     36353592052137018784330858422674753425698363\,{t^5} 
$$
$$- 
     7639385747972665348009342010125607135270845\,{t^6} 
$$
$$+ 
     93749989416978017153766638267237965058515\,{t^7} 
$$
$$+ 
     613514003480165484061014972915970847589645\,{t^8} 
$$
$$- 
     277973572971202497026511206431230055671555\,{t^9} 
$$
$$+ 
     4186799108745525715968930085758736947197\,{t^{10}} 
$$
$$+ 
     68475809505229552273919737578535787115805\,{t^{11}} 
$$
$$- 
     41025357958210023316522198089867194215575\,{t^{12}} 
$$
$$+ 
     12129433061687109251202289166065827342585\,{t^{13}} 
$$
$$- 
     1903251796631667579314923895099325939615\,{t^{14}} 
$$
$$\left. + 
     126883453108777838620994926339955062641\,{t^{15}} \right) \,{y^{18}}.
$$

\vfill\eject 

\baselineskip=12pt
\noindent{\bf BIBLIOGRAPHY }
\vskip 0.3 cm

\item{[AgLak]}
R.P. Agarwal and V. Lakshmikantham, {\it Uniqueness and Non-uniqueness 
Criteria for Ordinary Differential Equations,}\/ World Scientific (1993).
\vskip 0.2 cm

\item{[Bat]}
H. Bateman and A. Erdelyi, {\it Higher Transcendental Functions,}\/
Krieger Publishing, Florida (1981).
\vskip 0.2 cm

\item{[Bir]}
J.S. Birman, {\it Braids, Links, and Mapping Class Groups,}\/ Ann. Math. 
Stud. Princeton University (1975).
\vskip 0.2 cm

\item{[Cox]}
H.S.M. Coxeter, {\it Regular Polytopes,}\/ Macmillan Company New York 
(1963).
\vskip 0.2 cm

\item{[Cro]}
W.J.R. Crosby, Problems and Solutions, {\it Mounthly Amer. Math.}\/ (1946) 
104-107.
\vskip 0.2 cm

\item{[Dek]}
W. Dekkers, The Matrix of a Connection having Regular Singularities 
on a Vector Bundle of Rank $2$ on $P^1(\complessi)$, {\it Springer Lect. Notes 
Math.}\/ {\bf 712} (1979) 33-43.
\vskip 0.2 cm

\item{[Dub]} 
B. Dubrovin, Geometry of 2D Topological Field Theories, {\it Springer Lect. 
Notes Math.}\/ {\bf 1620}, (1995) 120-348.
\vskip 0.2 cm

\item{[FlN]}
H. Flashka and A.C. Newell, Monodromy and Spectrum Preserving 
Deformations, {\it Comm. Math. Phys.}\/  {\bf 76} (1980) 67-116.
\vskip 0.2 cm

\item{[Fuchs]}
R. Fuchs, \"Uber Lineare Homogene Differentialgleichungen Zweiter 
Ordnung mit im drei im Endrichen Gelegene Wesentlich Singul\"aren 
Stellen, {\it Math. Ann.}\/ {\bf 63} (1907) 301-321.
\vskip 0.2 cm

\item{[Gamb]}   
B. Gambier, Sur les Equations Differentielles du Second Ordre et du 
Primier Degr\`e dont l'Integrale est a Points Critiques Fixes, {\it Acta 
Math.}\/  {\bf 33}, (1910) 1-55.
\vskip 0.2 cm

\item{[Gar]}
R. Garnier, Sur des Equations Differentielles du Troisieme Ordre dont 
l'Integrale Generale est Uniforme et sur una Classe d'Equations nouvelles 
d'Ordre Superieur dont  l'Integrale Generale a ses Points Critiques Fixes,
{\it Ann. Sci. Ecole Norm. Sup.}\/  {\bf 3} No.29 (1912) 1-126.
\vskip 0.2 cm

\item{[Gor]}
P. Gordan, Uber endliche Gruppen linearer Transformationen einer 
Veranderlichen, {\it Math. Ann.}\/  {\bf 12}, 1877,  23-46.
\vskip 0.2 cm

\item{[Hit]}
N.J. Hitchin, Twistor Spaces, Einstein Metrics and Isomonodromic 
Deformations, {\it J. Differential Geometry}\/ {\bf 42}, No.1 July 1995 30-112.
\vskip 0.2 cm

\item{[Hit1]}
N.J. Hitchin, Poncelet Polygons and the Painlev\'e Transcendents, 
{\it Geometry and Analysis,}\/ edited by Ramanan, Oxford University Press 
(1995) 151-185.
\vskip 0.2 cm

\item{[Hit2]}
N.J. Hitchin, Talk at I.C.T.P. (Trieste) April 1993.
\vskip 0.2 cm

\item{[IKSY]}
K. Iwasaki, H. Kimura, S. Shimomura, M. Yoshida, {\it From Gauss to 
Painlev\'e, a Modern Theory for Special Functions,}\/  Chapter 4,
Aspects of Mathematics, Vieweg (1991).
\vskip 0.2 cm

\item{[Ince]}
E.L. Ince, {\it Ordinary Differential Equations,}\/
Dover Publications, New York (1956).
\vskip 0.2 cm

\item{[ItN]} 
A.R. Its and V.Yu. Novokshenov, The Isomonodromic Deformation Method 
in the Theory of Painlev\'e Equations, {\it Springer Lect. Notes Math.}\/ 
{\bf 1191} (1986).
\vskip 0.2 cm

\item{[Jim]}
M. Jimbo, Monodromy Problem and the Boundary Condition for some 
Painlev\'e Equations, {\it Publ. RIMS}\/ {\bf 18} (1982), 1137-1161.
\vskip 0.2 cm

\item{[JMU]}
M. Jimbo, T. Miwa and K. Ueno, Monodromy Preserving Deformation of 
the Linear Ordinary Differential Equations with Rational Coefficients I, 
II, {\it Publ. RIMS}\/ {\bf 14} (1978) 223-267, {\bf 15} (1979) 201-278.
\vskip 0.2 cm

\item{[Kit]}
A.V. Kitaev, Rational Solutions of the Fifth Painlev\'e Equation,
{\it Differential and Integral Equations}\/ {\bf 7} (1994) no.3-4, 967-1000.
\vskip 0.2 cm

\item{[Kr]}
L. Kronecker, {\it J. de Math.}\/ {\bf 19} (1854) 177-192.
\vskip 0.2 cm

\item{[Lak]}
V. Lakshmikantham and S. Leela, {\it Differential and Integral 
Inequalities,}\/ Academic Press New York (1969). 
\vskip 0.2 cm

\item{[Luka]}
N.A. Lukashevich, On the Theory of the Third Painlev\'e Equation, 
{\it Differ. Uravn.}\/ {\bf 3} (1967) no.11, 1913-1923.
\vskip 0.2 cm

\item{[Luke]}
L. Yu. Luke, {\it Mathematical Functions and Their Approximations,}\/
Academic Press (1975).
\vskip 0.2 cm

\item{[Mag]}
W. Magnus, Ring of Fricke Characters and Automorphism Groups of Free 
Groups, {\it Math. Z.}\/ {\bf 170}, (1980) 91-103.
\vskip 0.2 cm

\item{[Mal]}
B. Malgrange, Sur les Deformations Isomonodromiques I, Singularit\'es 
R\'eguli\`eres, Seminaire de l'Ecole Normale Superieure 1979-1982, 
{\it Progress in Mathematics}\/ {\bf 37}, Birkh\"auser, Boston (1983) 401-426.
\vskip 0.2 cm

\item{[Man]}
Yu.I. Manin, Sixth Painlev\'e Equation, Universal Elliptic Curve and 
Mirror of ${\bf P}^2$, preprint al-geom$/9605010$ (1996).
\vskip 0.2 cm

\item{[Ma]}
M. Mazzocco, Picard and Chazy Solutions to PVI Equation, preprint 1998.
\vskip 0.2 cm

\item{[Miwa]}
T. Miwa, Painlev\'e Property of Monodromy Preserving Deformation 
Equations and the Analyticity of $\tau$-function, {\it Publ. RIMS, Kyoto 
Univ.}\/ {\bf 17} (1981) 703-721.
\vskip 0.2 cm

\item{[Mur]}
Y. Murata, Classical Solutions of the Third Painlev\'e Equation, {\it 
Nagoya Math. J.}\/ {\bf 139} (1995) 37-65.
\vskip 0.2 cm

\item{[Mur1]}
Y. Murata, Rational Solutions of the Second and Fourth Painlev\'e 
equations, {\it Funk. Ekvacioj}\/ {\bf  28}, no.1 (1995).
\vskip 0.2 cm

\item{[Nor]}
N.E. N\"orlund, The Logarithmic Solutions of the Hypergeometric 
equation, {\it Mat. Fys. Skr. Dan. Vid. Selsk.}\/ {\bf 2}, no.5 (1963) 1-58.
\vskip 0.2 cm

\item{[Ok]}
K. Okamoto, On the $\tau$-Function of the Painlev\'e Equations, 
{\it Physica}\/ {\sl 2D}\/ (1981) 525-535.
\vskip 0.2 cm

\item{[Ok1]}
K. Okamoto, Studies on the Painlev\'e Equations I, Sixth Painlev\'e 
Equation, {\it Ann. Mat. Pura Appl.}\/ {\bf 148} (1987) 337-381.
\vskip 0.2 cm

\item{[Ok2]}
K. Okamoto, Painlev\'e Equations and Dynkin Diagrams, {\it Painlev\'e 
Transcendents,}\/ edited by D. Levi and P. Winternitz, Plenum Press, 
New York (1992) 299-313.
\vskip 0.2 cm

\item{[Ok3]}
K. Okamoto,  Sur les Feuilletages Associ\'es aux Equation du Second 
Ordre \`a Points Critiques Fixes de Painlev\'e, Espaces des Conditions 
Initiales, {\it Japan. J. Math.}\/ {\bf 5} (1979), 1-79. 
\vskip 0.2 cm

\item{[Pain]}
P. Painlev\'e, Sur les Equations Differentielles du Second Ordre et
d'Ordre Superieur, dont l'Interable Generale est Uniforme, {\it Acta Math.}\/ 
{\bf 25} (1902) 1-86.
\vskip 0.2 cm

\item{[Pic]}
E. Picard, M\'emoire sur la Th\'eorie des Functions Alg\'ebriques de deux 
Varables, {\it Journal de Liouville}\/ {\bf 5} (1889), 135-319.
\vskip 0.2 cm

\item{[Seg]}
J. Segert, Frobenius Manifolds from Yang-Mills Instantons, preprint 
October (1997) dg-ga$/9710031$.
\vskip 0.2 cm

\item{[SMJ]}
M. Sato, T. Miwa and M. Jimbo, Holonomic Quantum Fields II, {\it
Publ. RIMS}\/ {\bf 15} (1979) 201-278. 
\vskip 0.2 cm

\item{[Sch]}
L. Schlesinger, \"Uber eine Klasse von Differentsial System 
Beliebliger Ordnung mit Festen Kritischer Punkten, {\it J. fur Math.}\/
{\bf 141} (1912), 96-145.
\vskip 0.2 cm

\item{[Schw]}
H.A. Schwartz, \"Uber Diejenigen F\"alle in Welchen die Gaussische 
Hypergeometrische Reihe einer Algebraische Funktion iheres vierten 
Elementes Darstellit, {\it Crelle J.}\/ {\bf 75} (1873) 292-335.
\vskip 0.2 cm

\item{[S1]} 
S. Shimomura, Painlev\'e Transcendents in the Neighbourhood of a Fixed 
Singular Point, {\it Funkcial. Ekvac.}\/ {\bf 25} (1982) 163-184. 
\vskip 0.2 cm

\item{[S2]} 
S. Shimomura, Series Expansions for Painlev\'e Transcendents in the 
Neighbourhood of a Fixed Singular Point, {\it Funkcial. Ekvac.}\/ {\bf 25} 
(1982) 185-197. 
\vskip 0.2 cm

\item{[S3]} 
S. Shimomura, Supplement to ``Series Expansions for Painlev\'e 
Transcendents in the Neighbourhood of a Fixed Singular Point'', {\it 
Funkcial. Ekvac.}\/ {\bf 25} (1982) 363-371. 
\vskip 0.2 cm

\item{[Sib]}
Y. Sibuya, {\it Linear Differential Equations in the Complex Domain: Problems 
of Analytic Continuation,}\/ AMS TMM {\bf 82} (1990).
\vskip 0.2 cm

\item{[Tod]}
K.P. Tod, Self-dual Einstein Metrics from the Painlev\'e VI Equation, {\it 
Phys. Lett. A,}\/ {\bf 190} (1994) 221-224.
\vskip 0.2 cm

\item{[Um]}
H. Umemura, Birational Automorphism Groups and Differential 
Equations, {\it Nagoya Math. J.,}\/ {\bf 119} (1990) 1-80.
\vskip 0.2 cm

\item{[Um1]}
H. Umemura,  Irreducebility of the First Differential Equation of 
Painlev\'e, {\it Nagoya Math. J.,}\/ {\bf 117} (1990) 231-252.
\vskip 0.2 cm

\item{[Um2]}
H. Umemura, H. Watanabe, Solutions of the Second and Fourth Differential 
Equations of Painlev\'e, to appear in {\it Nagoya Math. J.,}\/ (1997).
\vskip 0.2 cm

\item{[Um3]}
H. Umemura, H. Watanabe, Solutions of the Third Differential Equation of 
Painlev\'e, preprint (1997).
\vskip 0.2 cm

\item{[Wat]}
H. Watanabe, Birational Canonical Transformations and Classical Solutions 
of the Sixth Painlev\'e Equation, preprint of the Kyushu University, 
Japan (1997).
\vskip 0.2 cm

\item{[Wat1]}
H. Watanabe, Solutions of the Fifth Painlev\'e Equation I, {\it Hokkaido 
Math. J.}\/ XXIV, {\bf 2} (1995) p. 53-89.
\vskip 0.2 cm

\item{[Wey]}
H. Weyl, {\it Algebraic Theory of Numbers,}\/ Ann. Math. Stud. {\bf 1}, 
Princeton University Press (1940).
\vskip 0.2 cm

\bye